\newif\ifdraft		
\newif\iflualatex	
\newif\ifpictures 
\draftfalse
\lualatexfalse
\picturestrue

\documentclass[
	11pt,
	a4paper,
	\ifdraft draft\else final\fi,
	oneside%
]{article}

\iflualatex \usepackage{fontspec} \else \usepackage[T1]{fontenc} \fi
\usepackage{lmodern}
\usepackage[main=english, ngerman]{babel} 
\usepackage{fullpage}
\usepackage{tocloft}
\usepackage{caption}
\usepackage{xcolor} 
	\definecolor{darkblue}{RGB}{0,35,102}
	\definecolor{darkgreen}{RGB}{15,140,15}
\usepackage[final]{graphicx} 
	\ifdraft\ifpictures\else \setkeys{Gin}{draft} \fi\fi 
\usepackage{tabularray} 
	\UseTblrLibrary{amsmath}
\usepackage{csquotes} 
\usepackage{paralist} 

\usepackage{mathrsfs} 
\usepackage{amsmath, amsfonts, amssymb, amsthm} 
\usepackage{MnSymbol} 
	\DeclareSymbolFont{Symbols}{OMS}{cmsy}{m}{n}
	\SetSymbolFont{Symbols}{bold}{OMS}{cmsy}{b}{n}
	\DeclareMathSymbol{\setminus}{\mathbin}{Symbols}{"6E} 
\usepackage{mathtools} 
\usepackage{subdepth} 
\usepackage{dsfont} 
\usepackage{xfrac} 
\usepackage{esint} 

\usepackage[color,notcite,notref]{showkeys}
	\colorlet{labelkey}{orange}
	\colorlet{refkey}{orange}
	
\usepackage{todonotes} 
	\ifdraft\else\renewcommand{\todo}[2][]{} \fi 
\usepackage{changes} 
	\definechangesauthor[name=Kai, color=orange]{KR}

\usepackage{placeins} 
\usepackage{seqsplit} 
\usepackage{caption, subcaption} 
\usepackage{fontawesome5} 

\usepackage[final]{hyperref} 
	\hypersetup{
		colorlinks, 
		linkcolor=blue, 
		urlcolor=red, 
		filecolor=green, 
		citecolor=blue}
\usepackage[capitalise]{cleveref} 
\usepackage[open, openlevel=2, numbered]{bookmark}
\usepackage[backend=biber, style=alphabetic, giveninits=true]{biblatex} 
\usepackage{aliascnt}
\theoremstyle{plain} 
\newtheorem{theorem}{Theorem}[section]

\newaliascnt{lemma}{theorem}
\newtheorem{lemma}[lemma]{Lemma}
\aliascntresetthe{lemma}

\newaliascnt{proposition}{theorem}
\newtheorem{proposition}[proposition]{Proposition}
\aliascntresetthe{proposition}

\newaliascnt{corollary}{theorem}
\newtheorem{corollary}[corollary]{Corollary}
\aliascntresetthe{corollary}

\newaliascnt{remark}{theorem}
\newtheorem{remark}[remark]{Remark}
\aliascntresetthe{remark}

\theoremstyle{definition} 
\newaliascnt{definition}{theorem}
\newtheorem{definition}[definition]{Definition}
\aliascntresetthe{definition}

\newaliascnt{assumption}{theorem}
\newtheorem{assumption}[assumption]{Assumption}
\aliascntresetthe{assumption}

\newaliascnt{example}{theorem}

\aliascntresetthe{example}

\crefname{theorem}{Theorem}{Theorems}
\crefname{lemma}{Lemma}{Lemmas}
\crefname{proposition}{Proposition}{Propositions}
\crefname{corollary}{Corollary}{Corollaries}
\crefname{remark}{Remark}{Remarks}
\crefname{definition}{Definition}{Definitions}
\crefname{assumption}{Assumption}{Assumptions}
\crefname{example}{Example}{Examples}

\DeclareMathOperator{\SO}{SO}
\DeclareMathOperator{\Gl}{GL}
\DeclareMathOperator{\dist}{dist}

\DeclareMathOperator{\sym}{sym}
\DeclareMathOperator{\skewMat}{skew}
\newcommand{\el}{\mathrm{el}} 
\newcommand{\vi}{\mathrm{vi}} 
\newcommand{\limSymbol}{\boxtimes}
\newcommand{\R}{\mathbb{R}}

\newcommand{\N}{\mathbb{N}}
\newcommand{\Q}{\mathbb{Q}}
\newcommand{\Z}{\mathbb{Z}}
\newcommand{\eps}{\varepsilon}
\newcommand{\dd}{\mathrm{d}}
\newcommand{\wto}{\rightharpoonup}
\newcommand{\transp}{\top} 
\newcommand{\placeholder}{\makebox[1ex]{$\boldsymbol{\cdot}$}}
\newcommand{\grad}[1][]{\if\relax\detokenize{#1}\relax \nabla \else \nabla_{\!#1}\fi}
\newcommand{\dif}{\mathrm{D}}

\DeclarePairedDelimiter{\abs}{\lvert}{\rvert} 
\DeclarePairedDelimiter{\norm}{\lVert}{\rVert}
\DeclarePairedDelimiter{\of}{\lparen}{\rparen}
\DeclarePairedDelimiterX{\compose}[2]{\lparen}{\rparen}{#1 \:\delimsize\vert\allowbreak\:\mathopen{} #2}
\newcommand\setSymbol[1][]{\nonscript\:#1\vert\allowbreak\nonscript\:\mathopen{}} 
\newcommand{\given}{}
\DeclarePairedDelimiterX{\set}[1]{\lbrace}{\rbrace}{\renewcommand{\given}{\setSymbol[\delimsize]}#1} 
\DeclarePairedDelimiterX{\sprod}[2]{\langle}{\rangle}{#1\,,\,\mathopen{} #2}

\makeatletter
\def\fourvdots{\raisebox{-2pt}{\vbox{\baselineskip3\p@ \lineskiplimit\z@%
  \hbox{.}\hbox{.}\hbox{.}\hbox{.}}}}
\makeatother

\DeclareMathOperator{\Cont}{C}
\DeclareMathOperator{\SobH}{H}
\DeclareMathOperator{\SobW}{W}
\DeclareMathOperator{\Leb}{L}
\DeclareMathOperator{\AC}{AC}

\newcommand{\stepemph}[1]{\textsc{#1}}

\makeatletter
\newcommand{\matheqmargin}[1]{\setlength{\@mathmargin}{#1}}
\newcommand{\namedlabel}[2][]{\def\@currentlabelname{#1}\label{#2}}
\makeatother

	{\par\xdef\tpd{\the\prevdepth}\egroup\prevdepth=\tpd}
\makeatletter
	{\@fleqnfalse}%
	{\@fleqntrue\ignorespacesafterend}
\makeatother

\counterwithin{equation}{section} 

\allowdisplaybreaks 

\author{Stefan Neukamm, Kai Richter}
\title{Derivation of a bending model for rods from finite viscoelasticity}

\makeatletter
\hypersetup{
	pdftitle={\@title},
	pdfsubject={Elasticity},
	pdfauthor={Kai Richter},
	pdfkeywords={KEYWORD1, KEYWORD2},
	pdfdisplaydoctitle}
\makeatother

\DeclareDelimFormat[bib,biblist]{nametitledelim}{\addcolon\space}
\DeclareFieldFormat*{title}{\mkbibemph{#1}\isdot}
\DeclareFieldFormat*{journaltitle}{#1}

\AtEveryBibitem{%
  \clearfield{url}%
  \clearfield{issn}%
  \clearfield{isbn}%
  \clearfield{month}}
\begin{document}
\selectlanguage{english}%
\makeatletter%
\begin{flushleft}
	\def\refTUDsymbol{\textasteriskcentered}%
	\newcommand{\refTUD}{\hyperref[Title:TUD]{\textsuperscript{\refTUDsymbol}}}%
	\newcommand{\mail}[2]{\faEnvelope[regular] \hypersetup{hidelinks}\href{mailto:#1}{#1} (#2)}%
	\LARGE\bfseries
	\parbox[t]{0.8\linewidth}{\flushleft\@title}
	\par \bigskip

  \large\normalfont
  Stefan Neukamm\refTUD\ $\cdot$
  Kai Richter\refTUD
  \par

  \normalsize
	\mbox{} \hfill \today
	\par \bigskip
	
	\footnotesize 
	\textsuperscript{\refTUDsymbol}Faculty of Mathematics, Technische Universit{\"a}t Dresden, 01062 Dresden, Germany \label{Title:TUD}
	\par \medskip
	
	\mail{stefan.neukamm@tu-dresden.de}{S.~Neukamm} \\
	\mail{kai.richter@tu-dresden.de}{K.~Richter}
\end{flushleft}
\makeatother
\vspace{-2\baselineskip} \bigskip\bigskip

\paragraph*{Abstract.} 
Starting from a fully three dimensional model of finite viscoelasticity in the Kelvin--Voigt rheology, we derive a bending model for thin, viscoelastic rods. For this, we utilize a formulation of the model as a time-dependent metric gradient flow and pass to the limit by establishing convergence of curves of maximal slope using $\Gamma$-convergence.
A posteriori, we relate the metric gradient flow limits to gradient flows in Hilbert spaces and discuss their properties.
In contrast to previous results, where limit deformations are infinitesimal, in the bending regime, we deal with large rotations. This prevents us from using an additive ansatz for the construction of recovery sequences. Instead, we adapt from a previous result by the authors a multiplicative ansatz for elastoplastic Cosserat rods to the viscoelastic setting.
For regularization purposes, we consider strain gradient terms in the three dimensional model that vanish in the limit.
\bigskip

\begin{footnotesize}
\begin{compactdesc} 
	\item[\textbf{Keywords:}] viscoelasticity $\cdot$ metric gradient flow $\cdot$ $\Gamma$-convergence $\cdot$ dimension reduction $\cdot$ bending $\cdot$ rod $\cdot$ strain gradient elasticity $\cdot$ third gradients.
	\item[\textbf{MSC-2020:}] 74D10 (primary) $\cdot$ 74K10 $\cdot$ 74Q15 $\cdot$ 49J45 $\cdot$ 35A15 $\cdot$ 35Q74.
\end{compactdesc}
\end{footnotesize}

\paragraph{Acknowledgments.}
Stefan Neukamm and Kai Richter received support from the German Research Foundation (DFG) via the research unit FOR 3013, “Vector- and tensor-valued surface PDEs” (project number 417223351).

{\hypersetup{hidelinks} \tableofcontents}
\ifdraft\listoftodos\fi


\section{Introduction}
\label{Sec:Introduction}

In finite elasticity one usually needs to deal with nonlinear, nonconvex energies. Such models come naturally with many difficulties for analytic studies or numerical computations. An important branch of continuum mechanics is thus the derivation of effective models adapted to special situations.
In this paper we focus on a model for finite viscoelasticity in the presence of a special geometry, namely thin rods. Recently, efforts have been made to rigorously derive effective models for viscoelasticity upon others in terms of linearization \cite{FK18} and von~K{\'a}rm{\'a}n plate \cite{FK20}, ribbon \cite{FM22} and rod \cite{FM23} models. The von~K{\'a}rm{\'a}n regime assumes small displacements leading to a linear state space and at least partly linearized equations which is heavily used in the mentioned.
In this paper, we shall build on the tools developed in these works, but we leave the situation of infinitesimal deformations and study the bending regime.
Our main result is the rigorous justification of a viscoelastic bending model for rods as the evolutionary $\Gamma$-limit of three-dimensional finite viscoelasticity.
The authors recently derived an effective bending-torsion rod model for elastoplasticity in \cite{NR26}. The main achievement of said paper was the construction of a mutual recovery sequence for the energy and dissipation simultaneously. Adapting this construction to viscoelasticity is a major point in the present paper.
\medskip

With this dimension reduction result we add another work to the long list of $\Gamma$-convergence-like results for finite elasticity that started in \cite{FJM02} with the development of the geometric rigidity estimate. For bending rods many results have been achieved in static elasticity, see \cite{MM03,Neu12,CRS17,BNS20,BGKN22-PP,FKZ24}. However, for evolutionary problems, \cite{NR26} is to the authors knowledge the first result in this direction for the bending regime. The complexity of the bending regime is the geometrically nonlinear state space that is expected in the limit.
Extending the results for the bending regime to viscoelasticity is a natural next step. In fact, this theory has gained a lot of attention recently from the mathematical community, see e.g.\ \cite{MOS14,KR20,FKV21,Sen21,DJKV22,Gah24,Sen24,Mac25,BKK26-PP}. The list is by far not exhaustive. Furthermore, viscoelasticity is an important building block for many more advanced material models such as thermo- and poroviscoelasticity, see e.g.\ \cite{MR20,BFK23,BFM24,OL24,RS23}, so that we can expect that our result can be adapted to many interesting problems.
\medskip

To specify the situation, we consider the thin rod $\Omega^h := \omega \times hS$ where $\omega := (0,\ell)$ is the center line, $S \subset \R^2$ denotes the profile and $0 < h \ll 1$ is the thickness of the rod. The evolution of a viscoelastic body in Kelvin--Voigt rheology is governed by the force balance equation
\begin{equation} \label{Eq:Intro:PDE}
		\dif_y \mathcal E^h\of{t,y(t)} + \dif_{\dot y} \mathcal R^h\of[\big]{y(t),\dot y(t)} = 0 
		\qquad\text{for a.e.\ } t \in [0,T],
\end{equation}
where $y(t)$ is a deformation, $\mathcal E^h$ is the total energy of the system and $\mathcal R^h$ is the viscous dissipation potential. Throughout this paper we are interested in weak PDE solutions and gradient flow solutions that are based on a reformulation of the equation as a metric gradient flow, see \cref{Def:Intro:Metric-gradient-flow} below.
The total energy and the dissipation potential are defined as follows. After rescaling the domain to unit thickness (that is to $\Omega := \Omega^1$), we consider as a state space for the deformation $y$, the set
\begin{equation}
	\mathcal Q^h_{\rm bc} := \set*{y \in \SobW^{3,p_*}(\Omega, \R^3) \given \det \grad[h]y > 0 \text{ a.e.\ in } \Omega,\, y(0,\placeholder) = v_{\rm bc} + hR_{\rm bc}\mathbf{\bar{x}} \text{ a.e.\ in } S}.
\end{equation}
Here, we fix $p_* > 3$ and $(v_\mathrm{bc}, R_\mathrm{bc}) \in \R^3 \times \SO(3)$. We note that $\grad[h] := \compose*{\partial_1}{h^{-1}\bar\grad}$ with $\bar\grad := (\partial_2,\partial_3)$ is the rescaled gradient that appears due to the rescaling of the domain to unit thickness. In the state space we include one-sided clamped boundary conditions. The total energy $\mathcal E^h: [0,T] \times \mathcal Q^h_{\rm bc} \to \R$, is given as

\begin{subequations} \label{Eq:Intro:Energy}
\begin{align}
	\mathcal{E}^h(t,y) &:= \mathcal{E}^h_\el(y) + \mathcal{E}^h_{\rm sg}(y) - \sprod{l(t)}{y}, \text{ where} \\
	\mathcal{E}^h_\el(y) &:= h^{-2}\int_\Omega W_\el\of*{\grad[h]y} \quad\text{and} \\
	\mathcal{E}^h_{\rm sg}(y) &:= h^{-\alpha_C p_*}\int_\Omega H_\el^C\of*{\grad C_y} 
		 + h^{\alpha_{F} p_*}\int_\Omega H_\el^F\of*{\grad \grad[h]y}
 		 + h^{\alpha_\mathbb{F} p_*}\int_\Omega H_\el^\mathbb{F}\of*{\grad^2 \grad[h]y}.
\end{align}
\end{subequations}
As presented, the energy consists of three parts. $\mathcal E^h_\el$ is a standard elastic energy with scaling $h^{-2}$ according to the bending regime. We suppose that the elastic energy density $W_\el$ is frame-indifferent, non-degenerate and minimized at identity. The precise assumptions are presented in \cref{Sec:Assumptions}. 
Via the term $\sprod{l(t)}{y}$ we allow for a time-dependent external force where we assume for the loading $l \in \AC([0,T]; \SobH^1(\Omega, \R^3)')$.
Finally, $\mathcal E^h_{\rm sg}$ is a strain gradient term that is similar to the term considered in \cite{NR26}. The term is adapted to the bending regime as we propose a stronger regularization of the right Cauchy-Green strain tensor $C_y := \grad[h]y^\transp \grad[h]y$ compared to the full deformation gradient. The scaling coefficients $0 < \alpha_C, \alpha_F < 1$ and $\alpha_{\mathbb F} > 0$ satisfy compatibility conditions, see \cref{Ass:ElasticMaterialLaw}. The precise form of the strain gradient term is motivated by the goal to achieve that it vanishes in the limit while still providing sufficient regularity, since we include $\mathcal E^h_{\rm sg}$ mainly for regularization purposes as in many related works in the field, see \cite{FK18,FK20,FM23}. It is an open problem how to rigorously justify the derivations when no strain gradient terms are considered. We note that strain gradients have been present in the physics literature for a long time, see \cite{Tou62,Tou64,Min64}. In fact, already in \cite{Min64} a split of strain gradients into gradients of the Cauchy-Green strain tensor and gradients of the rotational part is considered.
\medskip

For our purpose it is convenient to model the dissipation by a dissipation distance $\mathcal{D}^h: \mathcal Q^h_{\rm bc} \times \mathcal Q^h_{\rm bc} \to [0,\infty)$, given by
\begin{equation} \label{Eq:Intro:Dissipation}
	\mathcal{D}^h(y_1, y_2) := \left( h^{-2} \int_\Omega \mathscr{D}^2_\vi(\grad[h]y_1, \grad[h]y_2) \right)^{1/2}.
\end{equation}
Here $\mathscr{D}_\vi$ is a pseudo-distance with precise assumptions given in \cref{Sec:Assumptions}. From the dissipation distance we can recover the dissipation potential via the formula
\begin{equation} \label{Eq:Intro:Dissipation-distance}
	\mathcal R^h(y,\dot y) := \lim_{\eps \to 0} \tfrac{\eps^{-2}}{2} \mathcal D^h(y + \eps \dot{y}, y)^2 = \tfrac{1}{4} \dif_{y_1}^2 \mathcal D^h(y,y)^2[\dot y, \dot y].
\end{equation}
This definition states that the stress depends on the deformation gradient and its velocity, which is a model for nonlinear viscoelasticity introduced in \cite{RE55}. We note that there are different models for nonlinear viscoelasticity, e.g.\ considering a multiplicative split of the deformation gradient into an elastic and a viscous part, see e.g.\ \cite{SY24,CKS24}.

\paragraph{The limiting models.}
The main goal of this paper is to relate the system \eqref{Eq:Intro:PDE} for small $h \ll 1$ to a limiting system by studying the limit $h \to 0$. Indeed, we obtain convergence of solutions to the solution of a limit system that has again a gradient flow structure given by
\begin{equation} \label{Eq:Intro:PDE-limit}
		\dif_q \mathcal E^0\of[\big]{t,q(t)} + \dif_{\dot q} \mathcal R^0\of[\big]{\dot q(t)} = 0 
		\qquad\text{for a.e.\ } t \in [0,T],
\end{equation}
where $q = (v,R,a,\varphi)$ denotes a limiting state and has the following form: $(v,R)$ is an inextensible Cosserat rod configuration, that is $(v,R) \in \mathcal A_{\rm rod}$ with
\begin{equation}
	\mathcal A_{\rm rod} := \set*{(v,R) \in \SobH^2(\omega, \R^3) \times \SobH^1(\omega, \R^{3\times 3}) \given \partial_1 v = R e_1, R \in \SO(3) \text{ a.e.}}.
\end{equation}
Such rod configurations have already been derived for bending problems in static elasticity, see \cite{MM03,Neu12,BNS20}. They describe the deformation of the rod by the deformation of the center line $v$ and and an orientation frame for the profile $R$, and thus allow for bending and twisting of the rod. $a \in \Leb^2(\omega)$ and $\varphi \in \Leb^2(\omega, \SobH^1_{\rm av}(S, \R^3))$ with
\begin{equation*}
	\SobH^1_{\rm av}(S, \R^3) := \set*{\varphi \in \SobH^1(S, \R^3) \given {\textstyle \fint_S} \varphi = {\textstyle \fint_S} \partial_3\varphi_2 - \partial_2\varphi_3 = 0}
\end{equation*}
are first order correctors for such deformations and lead to the state space
\begin{equation}
	\mathcal{Q}^0_{\rm bc} := \set*{(v,R,a,\varphi) \in \mathcal{A}_{\rm rod} \times \Leb^2(\omega) \times \Leb^2(\omega, \SobH^1_{\rm av}(S, \R^3)) \given v(0) = v_{\rm bc}, R(0) = R_{\rm bc}}.
\end{equation}
One has to be careful in how to interpret \cref{Eq:Intro:PDE-limit}, as $\mathcal A_{\rm rod}$ is not a vector space. A rigorous way to understand the equation is through the metric gradient flow formulation presented below. In \cref{Sec:Rod-models} we study the obtained viscoelastic rod models and derive a refined interpretation of \eqref{Eq:Intro:PDE-limit}. In fact, it turns out that the derivatives in \eqref{Eq:Intro:PDE-limit} can be defined with respect to an (infinite dimensional) manifold structure on $\mathcal Q^0_{\rm bc}$, see \cref{Rem:Manifold-structure-on-bending-deformations}.
The structure of the correctors $(a,\varphi)$ is already known from derivations in static elasticity, see \cite{Neu12}. For this static problem it turns out that the correctors are determined by corrector equations only depending on the rod configuration $(v,R)$ and the material law of the elastic energy. This observation was extended to elastoplastic evolutions in \cite{NR26}. This is convenient as it allows to define a reduced energy functional only depending on $(v,R)$ and fully reducing the equations to equations on the center line $\omega$. In contrast, for the viscoelastic setting studied in the present paper, the evolution of the correctors can in general not be separated from the evolution of the rod as the correctors influence the elastic as well as the dissipative processes. However, we provide in this paper a compatibility condition (see \eqref{Ass:C} in \Cref{Sec:Assumptions}), under which it is possible to show convergence to a reduced system
\begin{equation} \label{Eq:Intro:PDE-reduced-limit}
		\dif_q \mathcal E^{\rm red}\big(t, q(t)\big) + \dif_{\dot q} \mathcal R^{\rm red}\big(\dot q(t)\big) = 0 
		\qquad\text{for a.e.\ } t \in [0,T],
\end{equation}
defined on the reduced state spaces
\begin{equation}
	\mathcal Q^{\rm red}_{\rm bc} := \set*{(v,R) \in \mathcal A_{\rm rod} \given v(0) = v_{\rm bc}, R(0) = R_{\rm bc}}.
\end{equation}
Before we define the limiting and reduced energy and dissipation functionals, let us introduce the following notion of convergence, with which we can relate the limiting states to the states of the 3D system.

\begin{definition}
\label{Def:Intro:Convergence-rod-configurations}
Let $(y^h) \subset \SobH^1(\Omega, \R^3)$ and $q=(v,R,a,\varphi) \in \mathcal{A}_{\rm rod} \times \Leb^2(\omega) \times \Leb^2(\omega, \SobH^1_{\rm av}(S, \R^3))$. We say that \emph{$(y^h)$ converges to $(v,R)$ as reduced rod configuration} and write $y^h \xrightarrow{\rm red\, rod} (v,R)$ if
\begin{subequations}
\begin{align}
	y^h &\to \mathrlap{v}\phantom{R}
	\qquad \text{in } \Leb^2(\Omega, \R^3), 
	\label{Eq:Intro:Convergence-rod-configurations:a} \\
	\grad[h]y^h &\to R 
	\qquad \text{in } \Leb^2(\Omega, \R^{3\times 3}).
	\label{Eq:Intro:Convergence-rod-configurations:b}
\intertext{Furthermore, if additionally}
	h^{-1}\left(\sqrt{C_{y^h}} - I\right) &\wto \mathrlap{\sym\compose*{R^\transp\partial_1 R\mathbf{\bar x} + ae_1}{\bar\grad\varphi}} 
	&& \text{weakly in } \Leb^2(\Omega, \R^{3\times 3}),
	\label{Eq:Intro:Convergence-rod-configurations:c}
\end{align}
\end{subequations}
then we say \emph{$(y^h)$ converges to $q$ as rod configuration} and write $y^h \xrightarrow{\rm rod} q$. The limit $q$ is unique, see \cref{Lem:Non-degeneracy-of-Qvi}.
\end{definition}

We note that boundary conditions are preserved along sequences with bounded energy, i.e.\ $(y^h) \subset \mathcal{Q}^h_{\rm bc}$, $y^h \xrightarrow{\rm red\, rod} (v,R)$ and $\limsup_{h \to 0} \mathcal{E}^h_\el(y^h) < \infty$ imply $(v,R) \in \mathcal{Q}^{\rm red}_{\rm bc}$.
The limiting and reduced energies $\mathcal E^0: [0,T] \times \mathcal Q^0_{\rm bc} \to \R$ and $\mathcal E^{\rm red}: [0,T] \times \mathcal Q^{\rm red}_{\rm bc} \to \R$ are given as
\begin{subequations}
\begin{alignat}{3}
	\mathcal{E}^0(t,q) &:= \mathcal{E}^0_\el(q) - \sprod{l(t)}{v}, &&\text{ where }  &
	\mathcal{E}^0_\el(q) &:= \int_\Omega Q_\el\compose*{K_R\mathbf{\bar{x}} + a e_1}{\bar{\grad}\varphi}, 
		\label{Eq:Intro:limit-energy}\\
	\mathcal{E}^{\rm red}(t,v,R) &:= \mathcal{E}^{\rm red}_\el(v,R) - \sprod{l(t)}{v}, &&\text{ where } &
	\mathcal{E}^{\rm red}_\el(v,R) &:= \int_\omega Q_\el^{\rm red}(K_R)
		\label{Eq:Intro:reduced-energy}
\end{alignat}
\end{subequations}
for $q = (v,R,a,\varphi) \in \mathcal Q^0_{\rm bc}$, where we define the bending-torsion strain $K_R := R^\transp\partial_1 R$, the map $\mathbf{\bar x}(x) := x_2e_2 + x_3e_3$ and the linearized elastic energy density $Q_\el(\placeholder) := \frac12\dif_F^2 W_\el(I)[\placeholder,\placeholder]$. The reduced linearized elastic energy density $Q_\el^{\rm red}$ is defined in \eqref{Eq:Reduced-energy}.
The energy consists of an elastic term $\mathcal E^0_\el$ (resp.\ $\mathcal E^{\rm red}_\el$) and a loading term and coincides with the bending energy in the static theory of finite elasticity, see \cite{Neu12}, with potentially explicit dependence on the correctors. The energy is still nonconvex as $K_R$ depends nonlinearly on $R$ but the complexity is highly reduced as the linearized elastic energy densities $Q_\el$ and $Q^{\rm red}_\el$ are convex quadratic forms. 
For the dissipation distances $\mathcal D^0: \mathcal Q^0_{\rm bc} \times \mathcal Q^0_{\rm bc} \to \R_+$ and $\mathcal D^{\rm red}: \mathcal Q^{\rm red}_{\rm bc} \times \mathcal Q^{\rm red}_{\rm bc} \to \R_+$ we obtain
\begin{subequations}
\begin{align}
	\mathcal{D}^0(q, \hat{q}) &:= \left( \int_\Omega Q_\vi\compose*{K_R\mathbf{\bar{x}} - K_{\hat{R}}\mathbf{\bar{x}} + a e_1 - \hat{a}e_1}{\bar{\grad}\varphi - \bar{\grad}\hat{\varphi}} \right)^{1/2}, \\
	\mathcal{D}^{\rm red}\of*{(v,R), (\hat{v},\hat{R})} &:= \left( \int_\omega Q_\vi^{\rm red}\of*{K_R - K_{\hat{R}}} \right)^{1/2},
\end{align}
\end{subequations}
where, $Q_\vi(\placeholder) := \tfrac{1}{2} \dif_{F_1}^2 \mathscr{D}_\vi^2(I,I)[\placeholder, \placeholder]$ and $Q^{\rm red}_\vi$ is defined in \eqref{Eq:Reduced-dissipation}.

\paragraph{Reformulation as a metric gradient flow and main result.}
For the rigorous justification of the limit model we reformulate the viscoelastic model \eqref{Eq:Intro:PDE} as a metric gradient flow, see \cite{AGS05}. This approach was proposed in \cite{MOS14} and used successfully in previous works to derive limiting models using $\Gamma$-convergence techniques, see \cite{FK18,FK20,FM23}. As the limits are metric gradient flows as well, we also provide a reformulation of \eqref{Eq:Intro:PDE-limit} and \eqref{Eq:Intro:PDE-reduced-limit}.
We show in \cref{Sec:Rod-models,Sec:Characterization-3D} that in fact the metric reformulations are equivalent to \eqref{Eq:Intro:PDE}, \eqref{Eq:Intro:PDE-limit} and \eqref{Eq:Intro:PDE-reduced-limit}.
In contrast to the works cited above, we allow for a time-dependent energy and therefore use a generalized notion of metric gradient flows, proposed in \cite{RMS08,FG17}. To introduce this concept, consider a metric space $Q$ with metric $D$ and a map $E: [0,T]\times Q \to (-\infty,\infty)$. Let us introduce the necessary concepts.

\begin{definition}[cf.\ {\cite[Chap.~1]{AGS05}}, {\cite[Sec.~1]{RMS08}} and {\cite[Def.~2.5]{FG17}}]~
\label{Def:Intro:Metric-gradient-flow}
\begin{enumerate}[(a)]
	\item We say a curve $q: [0,T] \to Q$ is \emph{absolutely continuous} w.r.t.\ $D$ and write $q \in \AC([0,T]; Q)$, if there exists a map $m \in \Leb^1(0,T)$, such that
	\begin{equation}
		D\big(q(s), q(t)\big) \leq \int_s^t m(r)\,\dd r, \qquad \text{for all } 0 \leq s \leq t \leq T.
	\end{equation}
	We denote the smallest map $m$ with this property by $\abs*{q'}_{D}$ and call it the \emph{metric derivative} of $q$. One can show that $\abs*{q'}_{D}$ exists (see \cite[Thm.~1.1.2]{AGS05}) and is given by
	\begin{equation}
		\abs*{q'}_{D}(t) = \lim_{s \to t} \frac{D\big(q(s), q(t)\big)}{\abs*{t - s}} \qquad \text{ for a.e.\ } t \in [0,T].
	\end{equation}
	Moreover, if $Q$ is a reflexive Banach space and $D$ induced by the norm, then $q \in \SobW^{1,1}(0,T; Q)$ is the continuous representant, a.e.\ differentiable and $\abs*{q'}_{D} = \norm{\dot q}_Q$ (see \cite[Rem.~1.1.3]{AGS05}). 
	
	\item For $t \in [0,T]$, we denote by $\abs*{\partial E(t)}_{D}: Q \to [0,\infty]$ the \emph{local slope}, given by
	\begin{equation}
		\abs*{\partial E(t)}_{D}(q) := \limsup_{\hat q \to q} \frac{\big(E(t,q) - E(t,\hat q)\big)^+}{D(q,\hat q)}
	\end{equation}
	
	\item A curve $q: [0,T] \to Q$ is called a \emph{curve of maximal slope} for $(Q,D,E)$ (w.r.t.\ the local slope), if $q$ is absolutely continuous, $t \mapsto E(t,q(t))$ is absolutely continuous, the partial derivative $\partial_t E(t,q(t))$ exists for a.e.\ $t \in [0,T]$ and
	\begin{equation}
	\label{Eq:Intro:Curve-of-max-slope}
		\tfrac{\dd}{\dd t} E(t,q(t)) \leq \partial_t E(t,q(t)) - \tfrac12\abs*{q'}_{D}^2(t) - \tfrac12 \abs*{\partial E(t)}_{D}^2(q(t)) \qquad
		\text{for a.e.\ } t \in [0,T].
	\end{equation}
\end{enumerate}
\end{definition}

Curves of maximal slope take the role of solutions to the metric gradient flow specified by $(Q, D, E)$. To transfer these concepts to our setting, we note that indeed $(\mathcal Q^h_{\rm bc}, \mathcal D^h)$ is a metric space and $(\mathcal Q^0_{\rm bc}, \mathcal D^0)$ and $(\mathcal Q^{\rm red}_{\rm bc}, \mathcal D^{\rm red})$ are even complete metric spaces, see \cref{Sec:Rod-models}.
As the viscoelastic dissipation distances necessarily satisfy frame-indifference to be physical, see \cite{Ant98,MOS14}, definiteness of the distances is a priori not clear, but follows from Korn's inequality. For the finite-thickness case we use a nonlinear version of Korn's inequality, presented in \cite{CM15}. To prove existence of curves of maximal slope in the finite-thickness case, it is convenient to introduce the states of bounded energy,
\begin{equation} \label{Eq:Intro:state-space-bounded-energy}
	\mathcal Q^{h,M}_{\rm bc} := \set[\Big]{y \in \mathcal Q^h_{\rm bc} \given \sup_{t \in [0,T]} \mathcal E^h(t,y) \leq M} \qquad \text{for } M > 0.
\end{equation}
We show in \cref{Thm:existence-bulk} that curves of maximal slope exist for $h \ll 1$ with a delicate smallness assumption on the thickness $h$ that depends on the initial energy. The existence result is a consequence of \cite{RMS08}, where curves of maximal slope are approximated by the time-discretized minimizing movements scheme,
\begin{equation*}
	\left\{\begin{aligned}
		t_i &:= t_{i - 1} + \tau, \\
		y_i &\in \operatorname{Arg\,min}_{y \in \mathcal Q^h_{\rm bc}} \left(\mathcal E^h(t_i,y) + \tfrac{1}{2\tau} \mathcal D^h\of[\big]{y_{i-1},y}^2\right),
	\end{aligned}\right.
\end{equation*}
given $t_0 = 0$ and $y_0 \in \mathcal Q^h_{\rm bc}$. We note that an existence result for large strains has recently been obtained in \cite{Mac25}.
Existence and uniqueness of curves of maximal slope for the limiting system is discussed in \cref{Sec:Rod-models}. Let us now state our main results.

\begin{theorem}[Convergence of curves of maximal slope]
\label{Thm:Main-result:convergence-to-limit}
Consider \cref{Ass:ElasticMaterialLaw,Ass:ViscousMaterialLaw}. For $h > 0$, let $y^h: [0,T] \to \mathcal Q^h_{\rm bc}$ be a curve of maximal slope for $(\mathcal Q^h_{\rm bc}, \mathcal D^h, \mathcal E^h)$. Assume that the initial value $y^h_0 := y^h(0)$ satisfies $y^h_0 \xrightarrow{\rm rod} q_0$ and $\mathcal E^h(0,y^h_0) \xrightarrow{h \to 0} \mathcal E^0(0,q_0)$ for some $q_0 \in \mathcal Q^0_{\rm bc}$.
Then, there exists $q: [0,T] \to \mathcal Q^0_{\rm bc}$, such that $y^h(t) \xrightarrow{\rm rod} q(t)$ for all $t \in [0,T]$ and
\begin{equation}
	\mathcal E^h\big(t,y^h(t)\big) \to \mathcal E^0\big(t,q(t)\big).
\end{equation}
The limit $q$ is the unique curve of maximal slope for $(\mathcal Q^0_{\rm bc}, \mathcal D^0, \mathcal E^0)$ which satisfies $q(0) = q_0$.
\end{theorem}

\begin{theorem}[Convergence of curves of maximal slope to the reduced system]
\label{Thm:Main-result:convergence-to-reduced}
Consider \cref{Ass:ElasticMaterialLaw,Ass:ViscousMaterialLaw} and the compatibility condition \eqref{Ass:C}. For $h > 0$, let $y^h: [0,T] \to \mathcal Q^h_{\rm bc}$ be a curve of maximal slope for $(\mathcal Q^h_{\rm bc}, \mathcal D^h, \mathcal E^h)$. Assume that $y^h_0 := y^h(0)$ satisfies $y^h_0 \xrightarrow{\rm red\, rod} (v_0,R_0)$ and $\mathcal E^h(0,y^h_0) \xrightarrow{h \to 0} \mathcal E^{\rm red}(0,v_0,R_0)$ for some $(v_0, R_0) \in \mathcal Q^{\rm red}_{\rm bc}$.
Then, there exists $(v,R): [0,T] \to \mathcal Q^{\rm red}_{\rm bc}$, such that $y^h(t) \xrightarrow{\rm red\,rod} (v(t),R(t))$ for all $t \in [0,T]$ and
\begin{equation}
	\mathcal E^h\big(t,y^h(t)\big) \to \mathcal E^{\rm red}\big(t,v(t),R(t)\big).
\end{equation}
The limit $(v,R)$ is the unique curve of maximal slope for $(\mathcal Q^{\rm red}_{\rm bc}, \mathcal D^{\rm red}, \mathcal E^{\rm red})$ which satisfies $v(0) = v_0$, $R(0) = R_0$.
\end{theorem}

The proof of \cref{Thm:Main-result:convergence-to-limit,Thm:Main-result:convergence-to-reduced} is due to a general convergence result for metric gradient flows which we present in \cref{Sec:Convergence-metric-gradient-flow}. This convergence result is an adaptation of \cite{Ser11,FK18}, see also \cite{FK20,FM23}, extended to the case of time-dependent energies.

\paragraph{Structure of the paper.} 
In \cref{Sec:Rod-models} we study the limiting viscoelastic rod models in a self-contained passage.
In \cref{Sec:Assumptions}, we introduce the assumptions on the 3D model. Afterwards in \cref{Sec:Compactness} we establish a compactness statement for sequences with bounded energy.
In \cref{Sec:3D-model}, we study the 3D model and show existence and establish a characterization of curves of maximal slope. Our main results \cref{Thm:Main-result:convergence-to-limit,Thm:Main-result:convergence-to-reduced} are proven in \cref{Sec:Convergence-metric-gradient-flow}. We provide a general convergence result for curves of maximal slope and afterwards check the assumptions for our setting. The main focus is here on establishing a lower bound statement for the local slopes.

\subsection{Notation}

Throughout this paper, we use the following notation which we recall from \cite{NR26}.

\begin{itemize}
	\item (Decomposition):
	$e_1, e_2, e_3$ denotes the canonical basis of $\R^3$. We denote by $\compose{A_1}{\bar{A}}$ the $n \times 3$ matrix with $A_1 \in \R^n$ as the first column and $\bar{A} \in \R^{n \times 2}$ as the two last columns. Vise versa, given $A \in \R^{n \times 3}$ we denote the first and last two columns as $A_1 \in \R^n$ and $\bar{A} \in \R^{n \times 2}$, respectively. Especially, we denote $\grad = \compose{\partial_1}{\bar{\grad}}$. For $x \in \R^3$ we also use the notation $x = (x_1, \bar{x})$.
	
	\item (Tangential and normal part):
	We denote the mapping of vectors in $\R^3$ to their normal part by $\mathbf{\bar{x}} : \R^3 \to \R^3$, $x \mapsto x_2e_2 + x_3e_3$ and identify it with the map $\mathbf{\bar{x}} : \R^2 \to \R^3$, $\bar{x} \mapsto x_2e_2 + x_3e_3$. We also define the map $\mathbf{x}_1: \R^3 \to \R$, $x \mapsto x_1$.
	
	\item (Subsets):
	We introduce notations for certain subsets of $\R$ and $\R^{3\times 3}$. We denote $\R_+ := [0,\infty)$ and $\bar{\R}_+ := [0,\infty] := \R_+ \cup \set*{\infty}$. Moreover, we set 
	\begin{alignat*}{2}
		\Gl(3) &:= \set*{F \in \R^{3\times 3} \given \det F \neq 0}, 
		&\quad \Gl_+(3) &:= \set*{F \in \R^{3\times 3} \given \det F > 0}, \\
		\SO(3) &:= \set*{F \in \Gl_+(3) \given F^\transp F = FF^\transp = I},
		&\quad \R^{3\times 3}_{\sym} &:= \set*{F \in \R^{3\times 3} \given F^\transp = F}, \\
		\R^{3\times 3}_{\skewMat} &:= \set*{F \in \R^{3\times 3} \given F^\transp = -F}.
	\end{alignat*}
	
	\item (Strain gradients):
	Let $R \in \R^{3\times 3}$, $\mathbb F \in \R^{3\times 3 \times 3}$ and $\mathscr F \in \R^{3\times 3 \times 3 \times 3}$. We define the multiplications $(R\mathbb F)_{ijk} := \sum_{n=1}^3 R_{in}\mathbb F_{njk}$ and $(R\mathscr F)_{ijkl} := \sum_{n=1}^3 R_{in}\mathscr F_{njkl}$. Furthermore, given $F \in \SobW^{2,1}_{\rm loc}(\R^3, \R^{3\times 3})$, we define the gradient $(\grad F)_{ijk} := \partial_k F_{ij}$ and second gradient $(\grad^2 F)_{ijkl} := \partial_{kl} F_{ij}$. With these definitions the multiplications above can be understand component wise for the partial derivatives, that is $(R\grad F)_{ijk} = (R\partial_k F)_{ij}$ and $(R\grad^2 F)_{ijkl} = (R\partial_{kl} F)_{ij}$.
	
	\item (Polar decomposition):
	Given $h > 0$ and $y \in \SobH^1(\Omega, \R^3)$ with $\det\grad[h]y > 0$ a.e., we denote the components of the polar decomposition of $\grad[h] y$ by $R_y$ and $\sqrt{C_y}$, i.e.\ $C_y := \grad[h]y^\transp \grad[h]y$ and $R_y := \grad[h]y (\sqrt{C_y})^{-1}$. We suppress the dependence of $h$ in this notation, but when used, the choice of $h$ can be taken from the context. It is common knowledge that with this definition, $R_y$ minimizes the distance of $\grad[h]y$ to $\SO(3)$, i.e.\ $R_y \in \SO(3)$ a.e.\ and
	\begin{equation}
		\dist\of[\big]{\grad[h]y, \SO(3)} = \abs*{\grad[h]y - R_y} = \abs*{\sqrt{C_y} - I}.
	\end{equation}
	
	\item (Rotation fields):
	We define the sets $\SobW^{k,p}_{\SO(3)}(\omega) := \SobW^{k,p}(\omega, \R^{3\times 3}) \cap \Leb^\infty(\omega, \SO(3))$ and $\SobH^k_{\SO(3)}(\omega) := \SobW^{k,2}_{\SO(3)}(\omega)$. Given $R \in \SobH^1_{\SO(3)}(\omega)$, we introduce $K_R := R^\transp \partial_1 R \in \Leb^2(\omega, \R^{3\times 3}_{\skewMat})$.
	
	\item We recall $\SobH^1_{\rm av}(S, \R^3) := \set*{\varphi \in \SobH^1(S, \R^3) \given {\textstyle \fint_S} \varphi = {\textstyle \fint_S} \partial_3\varphi_2 - \partial_2\varphi_3 = 0}$.
\end{itemize}

\section{The viscoelastic rod models}
\label{Sec:Rod-models}

In this section we discuss the two viscoelastic rod models and establish existence, uniqueness, time-regularity and a couple of principles that we frequently use in subsequent sections. 
The two models are formulated as (non-autonomous) metric gradient flows with a state space $\mathcal Q^\limSymbol_{\rm bc}$, energy functional $\mathcal E^\limSymbol: [0,T] \times \mathcal Q^\limSymbol_{\rm bc}\to\R$ and dissipation distance $\mathcal D^\limSymbol:\mathcal Q^\limSymbol_{\rm bc}\times\mathcal Q^\limSymbol_{\rm bc}\to\R_+$, see \cref{Sec:Introduction}. 
Here, the superscript $\limSymbol\in\{0,{\rm red}\}$ serves as a placeholder to allow for a simultaneous reference of the two systems. The case $\limSymbol={\rm red}$ corresponds to a ``reduced model'', which has a state space with reduced degrees of freedom; the case $\limSymbol=0$ corresponds to the ``extended'' model. In the later sections, we justify these models as evolutionary $\Gamma$-convergence limits of nonlinear three-dimensional viscoelasticity. However, for the purpose of this section, we may forget about this justification and only assume that $Q_\el, Q_\el^{\rm red}, Q_\vi, Q_\vi^{\rm red}$ are positive definite quadratic forms on the respective spaces. In fact, we shall see that it's justified to assume that $Q \in \set{Q_\el, Q_\vi}$ and $Q^{\rm red} \in \set{Q_\el^{\rm red}, Q_\vi^{\rm red}}$ satisfy
\begin{subequations}
\begin{alignat}{2}
	\tfrac1c \abs{\sym G}^2 &\leq Q(G) \leq c \abs{\sym G}^2, 
	&& \qquad\text{for all } G \in \R^{3\times 3}, 
	\\
	\tfrac1c \abs{K}^2 &\leq Q^{\rm red}(K) \leq c \abs{K}^2, 
	&& \qquad\text{for all } K \in \R^{3\times 3}_{\skewMat}.
	\label{Eq:Rod-models:Non-degneracy-Qred}
\end{alignat}
\end{subequations}

The viscoelastic rod models contain a lot of structure, which becomes apparent, when we change the variables from the deformation (plus orientation frame) $(v,R)$ to the bending-torsion strain $K$ as the elastic part of the energy, as well as the dissipation distances are quadratic w.r.t.\ $K$. The change of variables is achieved via the bijective map 
\begin{equation} \label{Eq:Deformation-to-strain-map}
	\Phi: \Leb^2(\omega, \R^{3\times 3}_{\skewMat}) \to \mathcal A^{\rm rod}_{\rm bc} \subset \SobH^2(\omega, \R^3) \times \SobH^1(\omega, \R^{3\times 3}),
\end{equation}
defined as $\Phi(K) = (\Phi_v(K),\Phi_R(K)) = (v_K, R_K)$, where $R_K$ is the unique solution to the initial value problem
\begin{equation}
	\partial_1 R = R K, \qquad
	R(0) = R_{\rm bc},
\end{equation}
and $v_K(x_1) := v_{\rm bc} + \int_0^{x_1} R_K e_1$. We study this map in \cref{Sec:Change-of-variables} and shall see that this is in fact a homeomorphism and twice Frech{\'e}t differentiable. This means that we identify $(v,R) \in \mathcal Q^{\rm red}_{\rm bc}$ with $K = \Phi^{-1}(v,R)$ and $q = (v,R,a,\varphi) \in \mathcal Q^0_{\rm bc}$ with $\kappa = (K,a,\varphi) = (\Phi^{-1}(v,R),a,\varphi)$. We observe that the transformed state spaces are in fact Hilbert spaces, which we denote by 
\begin{subequations}
\begin{alignat}{2}
	\mathcal H^0 &:= \Leb^2(\omega, \R^{3\times 3}_{\skewMat}) \times \Leb^2(\omega) \times \Leb^2(\omega, \SobH^1_{\rm av}(S, \R^3)) 
	&& \ni \kappa = (K,a,\varphi),
	\\
	\mathcal H^{\rm red} &:= \Leb^2(\omega, \R^{3\times 3}_{\skewMat})
	&& \ni K.
\end{alignat}
\end{subequations}
The distances $\mathcal D^\limSymbol$ correspond to the respective scalar products
\begin{subequations}
\begin{align}
	\sprod*{\kappa}{\hat \kappa}_{\mathcal H^0} &:= \int_\Omega \sym\compose*{K \mathbf{\bar x} + a e_1}{\bar\grad\varphi} : \mathbb D_\vi \sym\compose*{\hat K \mathbf{\bar x} + \hat a e_1}{\bar\grad\hat\varphi},
	\\
	\sprod*{K}{\hat K}_{\mathcal H^{\rm red}} &:= \int_\omega K : \mathbb D_\vi^{\rm red} K,
\end{align}
\end{subequations}
where $\mathbb D_\vi, \mathbb D_\vi^{\rm red}$ denote the bilinear forms associated to the quadratic forms $Q_\vi$ and $Q_\vi^{\rm red}$ respectively,
and we identify $\mathbb D_\vi$ with a fourth order tensor $\mathbb D_\vi \in \R^{3\times 3\times 3\times 3}$. These indeed define scalar products, as the following lemma shows.

\begin{lemma}[{cf.\ \cite[Lem.~3.2~(a)]{BNS20}}] \label{Lem:Non-degeneracy-of-Qvi}
There exists a constant $c > 0$, such that for all $\kappa = (K,a,\varphi) \in \mathcal H^0$, it holds,
\begin{equation}
	\norm{K}_{\Leb^2(\omega)} + \norm{a}_{\Leb^2(\omega)} + \norm{\varphi}_{\Leb^2(\omega,\SobH^1(S))}
	\leq c\, \norm[\big]{\sym\compose[\big]{K\mathbf{\bar x} + ae_1}{\bar\grad \varphi}}_{\Leb^2(\Omega)}.
\end{equation}
\end{lemma}

The transformed energies read,
\begin{subequations}
\begin{align}
	\mathcal I^0(t, \kappa) &:= \int_\Omega Q_\el\compose*{K\mathbf{\bar x} + ae_1}{\bar\grad\varphi} - \sprod*{l(t)}{v_K}, \\
	\mathcal I^{\rm red}(t,K) &:= \int_\omega Q_\el^{\rm red}(K) - \sprod*{l(t)}{v_K}.
\end{align}
\end{subequations}
We observe from \cref{Lem:Change-of-variables:homeomorphism} that for fixed $t \in [0,T]$, $\mathcal I^\limSymbol(t,\placeholder)$ is twice continuously Frech{\'e}t-differentiable with Lipschitz continuous derivative.
We are able to identify the viscoelastic rod models with the gradient flow of $\mathcal I^\limSymbol$ in the Hilbert space $\mathcal H^\limSymbol$.

\begin{proposition} \label{Prop:MetricGF-HilbertGF-correspondence}
Assume $l \in \AC([0,T]; \SobW^{1,\infty}(\omega, \R^3)')$.
\begin{enumerate}[(a)]
	\item $q= (v,R,a,\varphi): [0,T] \to \mathcal Q^0_{\rm bc}$ is a curve of maximal slope for $(\mathcal{Q}^0_{\rm bc}, \mathcal{D}^0, \mathcal{E}^0)$, if and only if $\kappa = (K,a,\varphi) = (\Phi^{-1}(v,R), a,\varphi) \in \Cont([0,T]; \mathcal H^0) \cap \SobH^1(0,T; \mathcal H^0)$ satisfies
	\begin{equation} \label{Eq:HilbertGF-0}
		\dot \kappa(t) = -\grad_{\mathcal H^0} \mathcal I^0\of[\big]{t,\kappa(t)}, \qquad
		\text{for a.e.\ } t \in [0,T]. 
	\end{equation}
	
	\item $(v,R): [0,T] \to \mathcal Q^{\rm red}_{\rm bc}$ is a curve of maximal slope for $(\mathcal{Q}^{\rm red}_{\rm bc}, \mathcal{D}^{\rm red}, \mathcal{E}^{\rm red})$, if and only if $K = \Phi^{-1}(v,R) \in \Cont([0,T]; \mathcal H^{\rm red}) \cap \SobH^1(0,T; \mathcal H^{\rm red})$ satisfies
		\begin{equation} \label{Eq:HilbertGF-red}
			\dot K(t) = -\grad_{\mathcal H^{\rm red}} \mathcal I^{\rm red}\of[\big]{t,K(t)}, \qquad
			\text{for a.e.\ } t \in [0,T]. 
		\end{equation}
\end{enumerate}
\end{proposition}

In view of the Lipschitz continuity of $\grad_{\mathcal H^\limSymbol} \mathcal I^\limSymbol(t,\placeholder)$, existence and uniqueness of (continuous) solutions to \eqref{Eq:HilbertGF-0} and \eqref{Eq:HilbertGF-red} can be obtained even classically from the Picard-Lindel{\"o}f theorem, in the case that merely $l \in \Leb^p(0,T; \SobW^{1,\infty}(\omega, \R^3)')$, $p \in [1,\infty]$, where we obtain $\SobW^{1,p}(0,T; \mathcal H^\limSymbol)$-regularity for the solutions. If $l \in \Cont([0,T], \SobW^{1,\infty}(\omega, \R^3)')$, then, the solution is even contained in $\Cont^1([0,T]; \mathcal H^\limSymbol)$.
Furthermore, in the situation of \cref{Prop:MetricGF-HilbertGF-correspondence}, existence and uniqueness can also be established by the theory of $\lambda$-convex (metric) gradient flows, see e.g.\ \cite{RMS08,FG17}. Indeed, the energies are $\lambda$-convex for some $\lambda \in \R$ (but not necessarily convex, as potentially $\lambda < 0$) in view of the Lipschitz continuity of the gradient. Especially, this approach shows that the solution can be approximated by the time-discretized scheme,
\begin{equation*}
	\left\{\begin{aligned}
		t_i &:= t_{i - 1} + \tau, \\
		\kappa_i &:= \operatorname{arg\,min}_{\kappa \in \mathcal H^\limSymbol} \left(\mathcal I^\limSymbol(t_i,\kappa) + \tfrac{1}{2\tau} \norm*{\kappa - \kappa_{i - 1}}_{\mathcal H^\limSymbol}^2\right).
	\end{aligned}\right.
\end{equation*}

In order to justify \cref{Prop:MetricGF-HilbertGF-correspondence}, let us first note that the local slopes and metric derivatives in $\mathcal Q^\limSymbol_{\rm bc}$ and $\mathcal H^\limSymbol$ coincide, which follows straight-forwardly from the definition. Moreover, as $\mathcal H^\limSymbol$ is a Hilbert space and the energy is $\lambda$-convex, we get further representations from \cite[Rem.~1.1.3, Cor.~1.4.5 and Thm.~2.4.9]{AGS05}. We present these statements in the following lemma but skip the proof.

\begin{lemma} \label{Lem:Rod-models:Characterization-slope}
Let $\limSymbol \in \set{0, \rm red}$. The following statements hold:
\begin{enumerate}[(a)]
	\item A map $q: [0,T] \to \mathcal Q^\limSymbol_{\rm bc}$ is absolutely continuous w.r.t.\ $\mathcal D^\limSymbol$, if and only if $\kappa: [0,T] \to \mathcal H^\limSymbol$ related to $q$ via $\Phi$ is absolutely continuous w.r.t.\ $\norm{\placeholder}_{\mathcal H^\limSymbol}$ (and hence contained in $\SobW^{1,1}(0,T; \mathcal H^\limSymbol)$) and then,
	\begin{equation}
		\abs*{q'}_{\mathcal D^\limSymbol}
		= \abs*{\kappa'}_{\norm{\placeholder}_{\mathcal H^\limSymbol}} 
		= \norm*{\dot\kappa(\placeholder)}_{\mathcal H^\limSymbol}
		\qquad 
		\text{a.e.\ in } [0,T].
	\end{equation}
	
	\item \label{Item:Rod-models:Characterization-slope} There exists $\lambda \in \R$, such that for all $q \in \mathcal Q^\limSymbol_{\rm bc}$, $t \in [0,T]$ and $\kappa \in \mathcal H^\limSymbol$ related to $q$ via $\Phi$, we have
	\begin{equation}
	\begin{aligned}
		\abs*{\partial \mathcal I^\limSymbol(t)}_{\norm{\placeholder}_{\mathcal H^\limSymbol}}(\kappa)
		&= \norm*{\grad_{\mathcal H^\limSymbol} \mathcal I^\limSymbol(t,\kappa)}_{\mathcal H^\limSymbol}
		\\
		&= \abs*{\partial \mathcal E^\limSymbol(t)}_{\mathcal D^\limSymbol}(q) 
		= \sup_{\hat q \neq q} \left(\frac{\mathcal E^\limSymbol(t,q) - \mathcal E^\limSymbol(t,\hat q) + \tfrac\lambda2 \mathcal D^\limSymbol(q,\hat q)^2}{\mathcal D^\limSymbol(q,\hat q)}\right)^+.
	\end{aligned}
	\end{equation}
	Moreover, the local slope $\abs*{\partial \mathcal E^\limSymbol(t)}_{\mathcal D^\limSymbol}$ is lower semi-continuous w.r.t.\ the metric $\mathcal D^\limSymbol$.
\end{enumerate}
\end{lemma}

Hence, $q: [0,T] \to \mathcal Q^\limSymbol_{\rm bc}$ is a curve of maximal slope for $(\mathcal Q^\limSymbol_{\rm bc}, \mathcal D^\limSymbol, \mathcal E^\limSymbol)$ if and only if the by $\Phi$ related map $\kappa:[0,T] \to \mathcal H^\limSymbol_{\rm bc}$ is a curve of maximal slope for $(\mathcal H^\limSymbol, \norm{\placeholder}_{\mathcal H^\limSymbol}, \mathcal I^\limSymbol)$. Then, \cref{Prop:MetricGF-HilbertGF-correspondence} is a classical result for gradient flows in Hilbert spaces, see e.g.\ \cite{AGS05,RMS08}. For the readers convenience, we provide the main arguments.

\begin{proof}[Proof of \cref{Prop:MetricGF-HilbertGF-correspondence}]
Provided $l \in \AC([0,T]; \SobW^{1,\infty}(\omega, \R^3)')$, we observe that for any absolutely continuous map $\kappa: [0,T] \to \mathcal H^\limSymbol$, $t \mapsto \mathcal I^\limSymbol\of[\big]{t,\kappa(t)}$ is absolutely continuous as well. Furthermore, at any point $t \in [0,T]$ where this map and $\kappa$ are differentiable, we find that by the Cauchy-Schwarz inequality,
\begin{align*}
	&\sprod*{\grad_{\mathcal H^\limSymbol} \mathcal I^\limSymbol(t,\kappa(t))}{\dot\kappa(t)}
	= \tfrac{\dd}{\dd t} \mathcal I^\limSymbol\of[\big]{t, \kappa(t)} - \partial_t \mathcal I^\limSymbol\of[\big]{t,\kappa(t)}
	\leq -\tfrac12 \norm*{\grad_{\mathcal H^\limSymbol} \mathcal I^\limSymbol\of[\big]{t,\kappa(t)}}_{\mathcal H^\limSymbol}^2 - \tfrac12 \norm*{\dot\kappa(t)}_{\mathcal H^\limSymbol}^2,
	\\
	&\qquad \Leftrightarrow \quad
	\dot \kappa(t) = -\grad_{\mathcal H^\limSymbol} \mathcal I^\limSymbol\of[\big]{t,\kappa(t)}.
\end{align*}
The claim follows, as the left-hand side of this equivalence characterizes curves of maximal slope and the right-hand side is \eqref{Eq:HilbertGF-0} or \eqref{Eq:HilbertGF-red}.
\end{proof}

\begin{remark} \label{Rem:Rod-models:Local-slope-maximizing-sequence}
The $\lambda$-convexity of the energy and the last characterization of the local slope in \cref{Lem:Rod-models:Characterization-slope}~(\ref{Item:Rod-models:Characterization-slope}) show that maximizing sequences can be chosen of the form $\kappa_s := \kappa + s (\hat \kappa - \kappa)$, where $\hat\kappa$ maximizes the right-hand side in \cref{Lem:Rod-models:Characterization-slope}~(\ref{Item:Rod-models:Characterization-slope}), cf.\ \cite[Thm.~2.4.9]{AGS05}. Indeed, we observe,
\begin{align*}
	\limsup_{s \to 0} \frac{\left( \mathcal I^\limSymbol(t,\kappa) - \mathcal I^\limSymbol(t,\kappa_s) \right)^+}{\norm*{\kappa_s - \kappa}_{\mathcal H^\limSymbol}}
	\geq \limsup_{s \to 0} \frac{s\left( \mathcal I^\limSymbol(t,\kappa) - \mathcal I^\limSymbol(t,\hat \kappa) + \frac{\lambda}{2}(1-s) \norm*{\hat \kappa - \kappa}_{\mathcal H^\limSymbol}^2 \right)^+}{s\norm*{\hat \kappa - \kappa}_{\mathcal H^\limSymbol}}&
	\\
	= \left(\frac{\mathcal E^\limSymbol(t,q) - \mathcal E^\limSymbol(t,\hat q) + \tfrac\lambda2 \mathcal D^\limSymbol(q,\hat q)^2}{\mathcal D^\limSymbol(q,\hat q)}\right)^+&.
\end{align*}
Pulled back to the Cosserat rod configurations, we obtain maximizing sequences, $(v_s, R_s) = \Phi(K_s)$ together with the correctors $a_s = a + s(\hat a - a)$, $\varphi_s = \varphi + s(\hat\varphi - \varphi)$. We use this fact later in \cref{Sec:Convergence-metric-gradient-flow} to construct a recovery sequence for the local slope when passing from the finite thickness to the rod models using evolutionary $\Gamma$-convergence, see \cref{Prop:LSC-local-slope}. For this purpose let us introduce,
\begin{equation}
	G_s := \compose*{K_{R_s}\mathbf{\bar{x}} + a_s e_1}{\bar{\grad}\varphi_s}, \qquad
	G := \compose*{K_R\mathbf{\bar{x}} + a e_1}{\bar{\grad}\varphi}, \qquad 
	\hat{G} := \compose*{K_{\hat R}\mathbf{\bar{x}} + \hat{a} e_1}{\bar{\grad}\hat{\varphi}}
\end{equation}
and note that,
\begin{equation}
\label{Eq:Rod-models:Local-slope-maximizing-sequence:Strain}
	G_s = G + s(\hat{G} - G), \qquad K_{R_sR^\transp} = R(K_{R_s} - K_R)R^\transp = s K_{\hat{R}R^\transp}.
\end{equation}
\end{remark}

\begin{remark} \label{Rem:Manifold-structure-on-bending-deformations}
We may view $\mathcal A^{\rm rod}_{\rm bc}$ as an (infinite dimensional) manifold modeled on the Hilbert space $\Leb^2(\omega, \R^{3\times 3}_{\skewMat})$, given by the single chart $(\Leb^2(\omega, \R^{3\times 3}_{\skewMat}), \Phi)$. Using this structure allows us to ``undo'' the change of variables in \cref{Prop:MetricGF-HilbertGF-correspondence} and consider the viscoelastic rod models as gradient flows w.r.t.\ this manifold structure, that is, curves of maximal slope for $\mathcal E^\limSymbol$ w.r.t.\ $\abs*{\partial\mathcal E^\limSymbol}_{\mathcal D^\limSymbol}$ are characterized by solutions to
\begin{equation*}
	\dif_q \mathcal E^\limSymbol(t,q) + \dif_{\dot q} \mathcal R^\limSymbol(q,\dot q) = 0,
\end{equation*}
where as in the finite thickness model we can recover from the dissipation distances $\mathcal D^\limSymbol$ the dissipation potential $\mathcal R^\limSymbol: \operatorname{T}\mathcal Q^\limSymbol_{\rm bc} \to [0,\infty)$ as
\begin{equation*}
	\mathcal R^\limSymbol(q,\dot q) := \tfrac{1}{4}\dif_{q_1}^2 \mathcal D^\limSymbol(q,q)^2[\dot q, \dot q] = \tfrac12 \norm*{\dot q}_{\mathcal H^\limSymbol}^2,
\end{equation*}
identifying the tangential spaces $\operatorname{T}_q\mathcal Q^\limSymbol_{\rm bc}$ with $\mathcal H^\limSymbol$.
\end{remark}

\section{Assumptions on the model}
\label{Sec:Assumptions}

In this section we provide in detail our assumptions on the viscoelasticity model \eqref{Eq:Intro:PDE} that we study in this paper. The model is determined by an energy $\mathcal E^h$ defined in \eqref{Eq:Intro:Energy} and a dissipation distance $\mathcal D^h$ given in \eqref{Eq:Intro:Dissipation}. These depend on the following quantities, which we introduce in this section:
\begin{itemize}
	\item the reference domain $\Omega \subset \R^3$,
	\item the elastic energy density $W_\el: \Gl_+(3) \to [0,\infty]$ and strain gradient energy densities $H_\el^C, H_\el^F: \R^{3\times 3\times 3} \to \R$ and $H_\el^\mathbb{F}: \R^{3\times 3\times 3 \times 3} \to \R$,
	\item the loading term $\sprod{l}{\placeholder}: [0,T] \times \mathcal Q^h \to \R$,
	\item the dissipation distance density $\mathscr{D}_\vi: \Gl_+(3) \times \Gl_+(3) \to [0,\infty)$.
\end{itemize}
We also discuss the quantities that we obtain for the limiting model \eqref{Eq:Intro:PDE-limit} and the reduced model \eqref{Eq:Intro:PDE-reduced-limit} which we obtain from the quantities above. Those are:
\begin{itemize}
	\item the linearized elastic energy density $Q_\el: \R^{3\times 3} \to \R_+$ and the reduced linearized elastic energy density $Q_\el^{\rm red}: \R^{3 \times 3}_{\skewMat} \to \R_+$, together with the associated fourth-order tensors $\mathbb C_\el$ and $\mathbb C_\el^{\rm red}$,
	\item the linearized dissipation distance density $Q_\vi: \R^{3\times 3} \to \R_+$ and the reduced linearized dissipation distance density $Q^{\rm red}_\vi: \R^{3\times 3}_{\skewMat} \to \R_+$, together with the associated fourth-order tensors $\mathbb D_\vi$ and $\mathbb D_\vi^{\rm red}$.
\end{itemize}
Finally, we define a compatibility condition for the elastic and viscous material law, which we only require for the reduced limit, cf.\ \cref{Thm:Main-result:convergence-to-reduced}.

\paragraph{Reference domain.} 
Let $S \subset \R^2$ be an open, bounded, connected Lipschitz domain and $\omega := (0,\ell)$ for some $\ell > 0$. We set $\Omega := \omega \times S$. We assume that $S$ is centered in $0$ the sense that
\begin{equation}
	\int_S x_2 \,\dd \bar{x} = \int_S x_3 \,\dd \bar{x} = \int_S x_2 x_3 \,\dd \bar{x} = 0.
\end{equation}

\paragraph{Loading term.}
We suppose $l \in \AC([0,T]; \SobH^1(\Omega, \R^3)')$.

\paragraph{Elastic energy.}
Let us introduce the following assumptions for the elastic energy density $W_\el: \Gl_+(3) \to [0,\infty]$:

\begin{enumerate}
\renewcommand{\labelenumi}{(WE\arabic{enumi})}
\renewcommand{\theenumi}{WE\arabic{enumi}}
\makeatletter
\def\@currentlabelname{\upshape (WE1) -- (WE4)}
\makeatother

	\item \label{Ass:WE}\label{Ass:W_FrameIndifference}(Frame indifference):
	\begin{equation*}
		W_\el(RF) = W_\el(F) \qquad \text{for all } F \in \Gl_+(3), R \in \SO(3).
	\end{equation*}
	
	\item \label{Ass:NaturalState}(Natural state):
		\begin{equation*}
			0 = W_\el(I) \leq W_\el(F) \qquad \text{for all } F \in \Gl_+(3).
		\end{equation*}
	
	\item \label{Ass:W_regularity}(Regularity): $W_\el$ is Borel measurable and there exists a neighborhood $U$ of $\SO(3)$ such that $W_\el \in \Cont^2(U)$ and $\dif^2 W_\el$ is Lipschitz continuous in $U$.
	
	\item \label{Ass:W_NonDegeneracy}(Non-degeneracy): There exists some $c > 0$, such that
	\begin{alignat*}{2}
		W_\el(F) &\geq c \dist^2\of[\big]{F, \SO(3)} &&\qquad \text{for all } F \in \Gl_+(3).
	\end{alignat*}
\end{enumerate}

The assumptions \nameref{Ass:WE} are standard. The frame indifference \eqref{Ass:W_FrameIndifference} is necessary to comply with physics and due to \eqref{Ass:NaturalState}, \eqref{Ass:W_regularity} and \eqref{Ass:W_NonDegeneracy} we can linearize $W_\el$ at the identity and define the linearized elastic energy density as the quadratic form
\begin{equation}
	Q_\el(G) := \tfrac{1}{2}\dif_F^2 W_\el(I)[G,G], 
	\qquad G \in \R^{3\times 3}.
\end{equation}
We introduce the associated symmetric fourth-order tensor $\mathbb{C}_\el \in \R^{3\times 3\times 3 \times 3}$ by
\begin{equation}
	G_2 : \mathbb{C}_\el G_1 = \tfrac{1}{2}\dif_F^2 W_\el(I)[G_1,G_2], 
	\qquad G_1, G_2 \in \R^{3\times 3}.
\end{equation}
The regularity assumption \eqref{Ass:W_regularity} is not a strong restriction but for many results weaker assumptions suffice to pass to a bending limit, see e.g.\ \cite{Neu12,NR26}. However, here we require for the results in \Cref{Sec:Lsc-slope} a cubic growth of the rest term in the quadratic Taylor expansion. We follow here the assumptions in \cite[Sec.~2.1]{FK18}.
\medskip

In our assumptions on the strain gradient term $\mathcal{E}^h_\mathrm{sg}$, we follow \cite{NR26} (see especially \cite[Sec.~A.1]{NR26}). It was shown there that this special form of the strain gradients is well-suitable for the bending regime, as it implies sufficient regularity while we can still expect that the terms vanish in the limit $h \to 0$. Additionally, we introduce here a third gradient term which ensures that the energy is convex in the highest order. This term is only relevant for the existence of solutions to the bulk models and not for the dimension reduction. Let us introduce the assumptions on the energy densities $H_\el^C, H_\el^F: \R^{3\times 3\times 3} \to \R$ and $H_\el^\mathbb{F}: \R^{3\times 3\times 3 \times 3} \to \R$:

\begin{enumerate}
\renewcommand{\labelenumi}{(HE\arabic{enumi})}
\renewcommand{\theenumi}{HE\arabic{enumi}}
\makeatletter
\def\@currentlabelname{\upshape (HE1) -- (HE4)}
\makeatother

	\item \label{Ass:HE}\label{Ass:HE_FrameIndifference}(Frame indifference):
	\begin{alignat*}{2}
		H_\el^F(R\mathbb{F}) &= H_\el^F(\mathbb{F}) &&\qquad \text{for all } \mathbb{F} \in \R^{3 \times 3 \times 3}, R \in \SO(3), \\
		H_\el^\mathbb{F}(R\mathscr{F}) &= H_\el^\mathbb{F}(\mathscr{F}) &&\qquad \text{for all } \mathscr{F} \in \R^{3 \times 3 \times 3 \times 3}, R \in \SO(3).
	\end{alignat*}
	
	\item \label{Ass:HE_Growth}($p_*$-growth): There exist $0 < c \leq C$, such that for each $H_\el \in \set*{H_\el^C, H_\el^F}$,
	\begin{alignat*}{3}
		c \abs*{\mathbb{F}}^{p_*} &\leq H_\el(\mathbb{F}) &&\leq C \left( \abs*{\mathbb{F}}^{p_*} + 1\right) &&\qquad \text{for all } \mathbb{F} \in \R^{3 \times 3 \times 3}, \\
		c \abs*{\mathscr{F}}^{p_*} &\leq H_\el^{\mathbb{F}}(\mathscr{F}) &&\leq C \left( \abs*{\mathscr{F}}^{p_*} + 1\right) &&\qquad \text{for all } \mathscr{F} \in \R^{3 \times 3 \times 3 \times 3}.
	\end{alignat*}
	
	\item \label{Ass:HE_Continuity}(Regularity): $H_\el^F$, $H_\el^C$ and $H_\el^{\mathbb F}$ are continuously differentiable and	there exist constants $L_C, L_F, L_\mathbb{F} \in \R_+$, such that
	\begin{align*}
		\abs*{H_\el^C(\mathbb{F}) - H_\el^C(\hat{\mathbb{F}})} &\leq L_C \left(\abs*{\mathbb{F}}^{p_*-1} + \abs*{\hat{\mathbb{F}}}^{p_*-1}\right) \abs*{\mathbb{F} - \hat{\mathbb{F}}} && \text{for all } \mathbb{F}, \hat{\mathbb{F}} \in \R^{3 \times 3 \times 3}, \\
		\abs*{H_\el^F(\mathbb{F}) - H_\el^F(\hat{\mathbb{F}})} &\leq L_F \left(\abs*{\mathbb{F}}^{p_*-1} + \abs*{\hat{\mathbb{F}}}^{p_*-1} + 1\right) \abs*{\mathbb{F} - \hat{\mathbb{F}}}  && \text{for all } \mathbb{F}, \hat{\mathbb{F}} \in \R^{3 \times 3 \times 3}, \\
		\abs*{H_\el^\mathbb{F}(\mathscr{F}) - H_\el^{\mathbb F}(\hat{\mathscr{F}})} &\leq L_\mathbb{F} \left(\abs*{\mathscr{F}}^{p_*-1} + \abs*{\hat{\mathscr{F}}}^{p_*-1} + 1\right) \abs*{\mathscr{F} - \hat{\mathscr{F}}}  && \text{for all } \mathscr{F}, \hat{\mathscr{F}} \in \R^{3 \times 3 \times 3 \times 3}.
	\end{align*}
	
	\item \label{Ass:HE_Convexity}(Convexity): $H_\el^F$, $H_\el^C$ and $H_\el^{\mathbb F}$ are convex.
\end{enumerate}

\begin{assumption}[Elastic material law] \label{Ass:ElasticMaterialLaw}
We assume that $W_\el$ satisfies \nameref{Ass:WE} and $H_\el^C, H_\el^F, H_\el^\mathbb{F}$ satisfy \nameref{Ass:HE} for some $p_* > 3$. Moreover, we fix $0 < \alpha_C, \alpha_F < 1$ and $\alpha_\mathbb{F} > 0$ with $\alpha_F < \frac{2}{3}(1- \alpha_C)$ and $\alpha_\mathbb{F} > \max\set*{\alpha_F, \frac{p_* - 2}{p_*}(1-\alpha_C)}$.
\end{assumption}

\begin{remark}~
\begin{enumerate}[(i)]
	\item Using a Taylor expansion of $W_\el$ at $I$ (see \cite[Lem.~4.1]{FK20}), it is not hard to show that \nameref{Ass:WE} imply that there exists a constant $c > 0$, such that
	\begin{equation} \label{Eq:TaylorEnergy}
		\abs[\big]{W_\el(I+G_2) - W_\el(I+G_1) - Q_\el(G_2) + Q_\el(G_1)} \leq c\,(\abs{G_1}^2 + \abs{G_2}^2)\abs*{G_2 - G_1},
	\end{equation}
	whenever $I + sG_1, I + sG_2 \in U$ for all $s \in [0,1]$.
	
	\item It is standard that $Q_\el$ is convex and non-degenerate on symmetric matrices in the sense that for some constant $c > 0$,
	\begin{equation} \label{Eq:SymmetryD2WatI}
		\tfrac{1}{c} \abs{\sym G}^2 \leq Q_\el(G) \leq c \abs{\sym G}^2
		\qquad \text{for all } G \in \R^{3\times 3}.
	\end{equation}
	Especially, $Q_\el(G)$ and $\mathbb{C}_\el G$ only depend on the symmetric part of $G$.
	
	\item More generally, \nameref{Ass:WE} yield similar to \cite[Lem.~A.1]{FM23} that in a potentially smaller neighborhood $\tilde U \subset U$ of $\SO(3)$, we find a constant $c > 0$, such that
	\begin{equation} \label{Eq:SymmetryD2WatF}
		\tfrac{1}{c} \abs*{\sym\of*{F^\transp G}}^2 \leq \dif_F^2 W_\el(F)[G,G] \leq c \abs*{\sym\of*{F^\transp G}}^2,
		\qquad \text{for all } F \in \tilde U, G \in \R^{3\times 3}.
	\end{equation}
\end{enumerate}
\end{remark}

\paragraph{Viscous dissipation distance.}
As proposed in \cite{MOS14} we model the dissipative effects via a dissipation distance. For our assumptions on the dissipation distance density $\mathscr{D}_\vi: \Gl_+(3) \times \Gl_+(3) \to [0,\infty)$ we follow \cite[Sec.~2.1]{FK18}:

\begin{enumerate}
\renewcommand{\labelenumi}{(DV\arabic{enumi})}
\renewcommand{\theenumi}{DV\arabic{enumi}}
\makeatletter
\def\@currentlabelname{\upshape (DV1) -- (DV4)}
\makeatother

	\item \label{Ass:DV}\label{Ass:DV_FrameIndifference}(Separate frame indifference):
	\begin{equation*}
		\mathscr{D}_\vi(R_1F_1, R_2F_2) = \mathscr{D}_\vi(F_1, F_2)
		\qquad \text{for all } F_1, F_2 \in \Gl_+(3), R_1, R_2 \in \SO(3).
	\end{equation*}
	
	\item \label{Ass:DV_Pseudometric} (Pseudometric): For all $F_1, F_2, F_3 \in \Gl_+(3)$, we have
	\begin{alignat*}{2}
		\mathscr{D}_\vi(F_1, F_2) &= \mathscr{D}_\vi(F_2, F_1), \\
		\mathscr{D}_\vi(F_1, F_2) &= 0 &&\qquad \text{if } F_1^\transp F_1 = F_2^\transp F_2, \\
		\mathscr{D}_\vi(F_1, F_2) &> 0 &&\qquad \text{if } F_1^\transp F_1 \neq F_2^\transp F_2, \\
		\mathscr{D}_\vi(F_1, F_3) &\leq \mathscr{D}_\vi(F_1, F_2) + \mathscr{D}_\vi(F_2, F_3).
	\end{alignat*}
	
	\item \label{Ass:DV_Regularity} (Regularity): $\mathscr D_\vi$ is Borel measurable and there exists a neighborhood $U$ of $\SO(3)$ such that $\mathscr{D}_\vi \in \Cont^2(U \times U)$ and $\dif^2 \mathscr{D}_\vi^2$ is Lipschitz continuous in $U\times U$.
	
	\item \label{Ass:DV_NonDegeneracy} (Non-degeneracy): There exists a constant $c > 0$ such that
	\begin{equation*}
		\mathscr{D}_\vi (I,F) \geq c \dist\of[\big]{F, \SO(3)}
		\qquad \text{for all } F \in \Gl_+(3).
	\end{equation*}
\end{enumerate}

\begin{assumption}
\label{Ass:ViscousMaterialLaw}
We assume that $\mathscr{D}_\vi$ satisfies \nameref{Ass:DV}.
\end{assumption}

As for the elastic energy density, in the limiting model a linearized dissipation distance is relevant, which is determined by the quadratic form,
\begin{equation}
	Q_\vi(G) :=  \tfrac{1}{2} \dif_{F_1}^2 \mathscr{D}_\vi^2(I,I)[G,G], \qquad G \in \R^{3\times 3},
\end{equation}
where $\dif^2 \mathscr{D}_\vi^2$ denotes the second derivative of $\mathscr{D}_\vi^2$ w.r.t.\ the first component. The associated fourth-order tensor $\mathbb{D}_\vi \in \R^{3\times 3 \times 3 \times 3}$ is given by
\begin{equation}
	G_2 : \mathbb{D}_\vi G_1 :=  \tfrac{1}{2} \dif_{F_1}^2 \mathscr{D}_\vi^2(I,I)[G_1,G_2], \qquad G_1,G_2 \in \R^{3\times 3}.
\end{equation}

\begin{remark}~ \label{Rem:TaylorDissipation}
\begin{enumerate}[(i)]
	\item It was shown in \cite[Lem.~4.1 and 4.3]{FK18} by a Taylor expansion argument that \nameref{Ass:DV} imply that there exists a constant $c > 0$ such that whenever $I + sG_1, I+sG_2 \in U$ for all $s \in [0,1]$, then
	\begin{equation} \label{Eq:TaylorDissipation}
		\abs[\big]{\mathscr{D}_\vi^2(I + G_1, I + G_2) - Q_\vi(G_2 - G_1)} \leq c\, (\abs{G_1} + \abs{G_2})\abs{G_2 - G_1}^2.
	\end{equation}
	
	\item Furthermore, \cite[Lem.~4.3]{FK18} shows that there exists a constant $c > 0$ such that whenever $F + s(\hat{F} - F) \in U$ for all $s \in [0,1]$, then
	\begin{equation} \label{Eq:TaylorDissipation2}
		\abs[\big]{\mathscr{D}_\vi^2(F,\hat F) - \tfrac{1}{2} \dif_{F_1}^2 \mathscr{D}_\vi^2(F,F)[\hat F - F,\hat F - F]} \leq c\, \abs{\hat F - F}^3.
	\end{equation}
	
	\item As for $Q_\el$ we observe from \nameref{Ass:DV} that $Q_\vi(G) = Q_\vi(\sym G)$ for all $G \in \R^{3\times 3}$ and for some constant $c > 0$,
	\begin{equation} \label{Eq:NonDegeneracy_Qvi}
		\tfrac{1}{c} \abs*{\sym G}^2 \leq Q_\vi(G) \leq c \abs*{\sym G}^2.
	\end{equation}
		
	\item More generally, as shown in \cite[Lem.~A.1]{FM23} we obtain that in a potentially smaller neighborhood $\tilde U \subset U$ of $\SO(3)$, we find a constant $c > 0$, such that
	\begin{equation} \label{Eq:NonDegeneracy_QviAtF}
		\tfrac{1}{c} \abs[\big]{\sym\of[\big]{F^\transp G}}^2 \leq \dif_{F_1}^2 \mathscr{D}_\vi^2(F,F)[G,G] \leq c \abs[\big]{\sym\of[\big]{F^\transp G}}^2
		\quad \text{for all } F \in \tilde U, G \in \R^{3\times 3}.
	\end{equation}
\end{enumerate}
\end{remark}

\paragraph{Compatibility condition.}
It has been shown for static elasticity \cite{Neu12,BNS20} and elastoplasticity \cite{NR26} that the bending limit can be stated in terms of a reduced linearized elastic energy density given by
\begin{equation} \label{Eq:Reduced-energy}
	Q_\el^{\rm red}(K) := \min_{\chi \in \mathbb{H}_{\rm rel}} \int_S Q_\el(\compose*{K\mathbf{\bar{x}}}{0} + \chi), \qquad
	K \in \R^{3\times 3}_{\skewMat},
\end{equation}
where we define the relaxation space as
\begin{equation}
	\mathbb{H}_{\rm rel} := \set*{\sym\compose*{a e_1}{\bar{\grad}\varphi} \given a \in \R, \varphi \in \SobH^1_{\rm av}(S, \R^3)} \subset \Leb^2(S, \R^{3\times 3}_{\sym}).
\end{equation}
Similarly, we define the reduced linearized dissipation distance as
\begin{equation} \label{Eq:Reduced-dissipation}
	Q_\vi^{\rm red}(K) := \min_{\chi \in \mathbb{H}_{\rm rel}} \int_S Q_\vi(\compose*{K\mathbf{\bar{x}}}{0} + \chi), \qquad
	K \in \R^{3\times 3}_{\skewMat}.
\end{equation}
Moreover, we denote the associated correctors by $\chi_\el(K), \chi_\vi(K) \in \mathbb{H}_{\rm rel}$, i.e.\ the unique minimizers satisfying
\begin{equation*}
	Q_\el^{\rm red}(K) = \int_S Q_\el(\compose*{K\mathbf{\bar{x}}}{0} + \chi_\el(K)), \qquad
	Q_\vi^{\rm red}(K) = \int_S Q_\vi(\compose*{K\mathbf{\bar{x}}}{0} + \chi_\vi(K)).
\end{equation*}
Since the correctors show up in both the energy and dissipation part, in the case of viscoelasticity it can in general not be expected that the energy and dissipation reduce to the effective quantities above. However, in this paper we show that this is the case under the following compatibility condition.

\begin{enumerate}
\renewcommand{\labelenumi}{(C\arabic{enumi})}
\renewcommand{\theenumi}{C\arabic{enumi}}
\makeatletter
\def\@currentlabelname{\upshape (C1)}
\makeatother

	\item \label{Ass:C}\label{Ass:C_Compatibility}(Compatibility): For all $K \in \R^{3\times 3}_{\skewMat}$, we have $\chi_\el(K) = \chi_\vi(K)$.
\end{enumerate}

\begin{remark}
A similar compatibility condition is assumed in \cite{FK20} for von~K{\'a}rm{\'a}n plates by the assumption that the limit model has zero Poisson ratio in $e_3$ direction, see \cite[Rem.~5.11]{FK20}.
This condition translates essentially to our setting to $\chi_\el(K) = \chi_\vi(K) = 0$ so that \eqref{Ass:C_Compatibility} is more general.
\end{remark}

\begin{remark}
It follows from an orthogonality argument shown in \cite[Prop.~2.13]{Neu12}, see also \cite[Lem.~2.10]{BNS20} and \cite[Lem.~2.8]{BGKN22-PP} that $Q_\vi^{\rm red}$ and $Q_\el^{\rm red}$ are positive definite, i.e.\ \eqref{Eq:Rod-models:Non-degneracy-Qred} is satisfied.
\end{remark}

\section{Compactness and geometric rigidity estimate}
\label{Sec:Compactness}

As a preparation for the dimension reduction analysis, we provide in this section compactness statements for deformations with bounded energy. The statements are similar to those provided in \cite{NR26}. As usual the compactness statements are due to the geometric rigidity estimate of \cite[Thm.~3.1]{FJM02} and a corollary of it for rod configurations $(v,R) \in \mathcal{A}_{\rm rod}$ due to \cite[Prop.~3.6]{Neu12}.
Here, we use a slight adaptation of the versions of \cite{NR26}, which we shall prove in the appendix.

\begin{proposition}[cf.\ {\cite[Prop.~4.1]{NR26}}] \label{Prop:RigidityEstimate}
Let $U \subset \R^3$ a bounded Lipschitz domain.
Let $1 < p < \infty$. There exists a constant $c > 0$, such that for all $\mathsf{y} \in \SobW^{1,p}(U, \R^3)$, we find a constant rotation $R \in \SO(3)$ (independent of $p$), such that
\begin{equation} \label{Eq:RigidityEstimate}
	\norm*{\grad \mathsf{y} - R}_{\Leb^p(U)} \leq c \norm*{\dist(\grad \mathsf{y}, \SO(3))}_{\Leb^p(U)}.
\end{equation}
\end{proposition}

\begin{proposition}[cf.\ {\cite[Prop.~4.2]{NR26}}] \label{Prop:RigidityRods} 
Let $h_0 > 0$.
For all $1 < p < \infty$, there exists a constant $c > 0$, such that for all $y \in \SobW^{1,p}(\Omega, \R^3)$ and any $0 < h \leq h_0$, we find a rotation field $R \in \SobW^{2,\infty}_{\SO(3)}(\omega)$ (independent of $p$), such that
\begin{subequations} \label{Eq:RigidityRods}
\begin{alignat}{4}
	\norm*{\grad[h] y - R}_{\Leb^p(\Omega)} + &\norm*{h \partial_1 R}_{\Leb^p(\omega)} &+& \norm*{h^2 \partial_1^2 R}_{\Leb^p(\omega)} &&\leq c \norm*{\dist\of[\big]{\grad[h]y, \SO(3)}}_{\Leb^p(\Omega)}, \label{Eq:RigidityRods:Lq}\\
	& \norm*{h \partial_1 R}_{\Leb^\infty(\omega)} &+& \norm*{h^2 \partial_1^2 R}_{\Leb^\infty(\omega)} && \leq c \norm*{\dist\of[\big]{\grad[h]y, \SO(3)}}_{\Leb^\infty(\Omega)}. \label{Eq:RigidityRods:Linf}
\end{alignat}
\end{subequations}
Moreover, if $y \in \mathcal{Q}^h_{\rm bc}$, then we may choose $R$, such that $R(0) = R_\mathrm{bc}$. \\
(For the proof see \cref{Sec:GRE-Proofs}.)
\end{proposition}

\begin{remark} \label{Rem:CosseratAnsatz}
The geometric rigidity estimate \cref{Prop:RigidityRods} motivates to represent $y \in \mathcal{Q}^h_{\rm bc}$ as
\begin{equation}
	y(x) := v(x_1) + hR(x_1) \mathbf{\bar{x}}(x) + h\phi(x),
\end{equation}
where $R^h[y] := R$ as in \cref{Prop:RigidityRods} and $v^h[y] := v = v_{\rm bc} + \int_0^{x_1} R e_1$, such that $(v, R) \in \mathcal{A}_\mathrm{rod}$ and $\phi^h[y] := \phi \in \SobW^{2,\infty}(\Omega, \R^3)$ with $v(0) = v_{\rm bc}$, $R(0) = R_{\rm bc}$ and $\phi(0,\placeholder) = 0$ a.e.\ in $S$. Indeed, $\phi$ is a first order correction of the Cosserat ansatz $v + hR\mathbf{\bar x}$ and relates to the correctors $(a,\varphi)$ for $q = (v, R, a,\varphi) \in  \mathcal Q^0_{\rm bc}$. We shall discuss this connection more thoroughly after the next proposition.
\end{remark}

From the preceding propositions we obtain the following compactness result for deformations $y^h$ with bounded energy. Up to a slight refinement, the result is analogous to \cite[Prop.~4.4]{NR26} and thus, we omit the proof here. Note that later in \cref{Lem:Regularity-energy} we show that the internal energy $\mathcal E^h_\el + \mathcal E^h_{\rm sg}$ controls the total energy $\mathcal{E}^h$.

\begin{proposition}[{cf.\ \cite[Prop.~4.4]{NR26}}] \label{Prop:Compactness}
Let $M > 0$, $h_0 > 0$, $p \in [1,\infty]$ and $\tilde{p} \in [1,\infty)$. There exists $\tilde{M} > 0$ such that for all sequences $(y^h) \subset \mathcal{Q}^h_{\rm bc}$ with $\mathcal{E}^h_\el(y^h) + \mathcal{E}^h_{\rm sg}(y^h) \leq M$ and $(v^h,R^h,\phi^h)$ as in \cref{Rem:CosseratAnsatz} the following statements hold:
\begin{enumerate}[(a)]
	\item We have for all $0 < h < h_0$,
	\begin{subequations}
	\begin{alignat}{10}
		&&& h^{-\alpha_C} && \norm[\Big]{\grad C_{y^h}}_{\Leb^{p_*}(\Omega)} 
		&&+ h^{\alpha_F} && \norm*{\grad \grad[h] y^h}_{\Leb^{p_*}(\Omega)} 
		&&+ h^{\alpha_\mathbb{F}} && \norm*{\grad^2 \grad[h] y^h}_{\Leb^{p_*}(\Omega)} 	\label{Eq:CompactnessStrainGradientsBound}\\
		&+ &&h^{-\alpha_{\tilde{p}}} && \norm*{\grad[h]y^h - R^h}_{\Leb^{\tilde{p}}(\Omega)} 
		&&+ h^{1-\alpha_p} &&\norm*{\partial_1 R^h}_{\Leb^p(\omega)}
		&&+ h^{2-\alpha_p} && \norm*{\partial_1^2 R^h}_{\Leb^p(\omega)} \label{Eq:CompactnessStrainBound} \\
		&+ && h^{1-\alpha_p^*} && \norm*{\phi^h}_{\Leb^p(\Omega)} 
		&&+ h^{1-\alpha_{\tilde{p}}} && \norm*{\grad[h]\phi^h}_{\Leb^{\tilde{p}}(\Omega)}
		&&+ (h^{2-\alpha_{p_*}} + h^{1+\alpha_F}) && \norm*{\grad\grad[h]\phi^h}_{\Leb^{p_*}(\Omega)} \label{Eq:CompactnessRestBound} \\
		&+ && h^{-\alpha_p} && \norm[\Big]{\sqrt{C_{y^h}} - I}_{\Leb^p(\Omega)} && \leq \tilde{M},
	\end{alignat}
	\end{subequations}
	where 
	$\alpha_p := \left\{\begin{smallmatrix*}[l] 1 & p \leq 2, \\ \alpha_C + 2\frac{1 - \alpha_C}{p} & p \in (2,\infty), \\ \alpha_C & p = \infty\end{smallmatrix*}\right.$ and 
	$\alpha_p^* = 1$ for $p < \infty$ and $\alpha_p^* < \alpha_3$ arbitrary for $p = \infty$.
	
	\item We find $q = (v,R,a,\varphi) \in \mathcal{Q}^0_{\rm bc}$, such that up to a subsequence (not relabeled), $(y^h)$ converges to $q$ as a rod configuration (recall \cref{Def:Intro:Convergence-rod-configurations}). Moreover, the subsequence can be chosen such that we additionally find $\phi \in \SobH^1(\omega, \R^3)$ and $\psi \in \Leb^2(\omega, \SobH^1(S,\R^3))$, such that
	\begin{subequations} \label{Eq:ConvergencesFromCompactnessAdvanced}
	\begin{align}
		v^h & \wto v && \text{in } \SobH^2(\omega, \R^3), \\
		R^h & \wto R && \text{in } \SobH^1(\omega, \R^{3\times 3}), \\
		\phi^h &\wto \phi && \text{in } \SobH^1(\Omega, \R^3), \label{Eq:CompactnessRestConvergence}\\
		\grad[h] \phi^h &\wto \compose*{\partial_1\phi}{\bar{\grad} \psi} && \text{in } \Leb^2(\Omega, \R^{3\times 3}), \\
		\sym\of*{(R^h)^\transp\grad[h]\phi^h} &\wto \sym\compose*{a e_1}{\bar{\grad} \varphi} && \text{in } \Leb^2(\Omega, \R^{3\times 3}),
	\end{align}
	\end{subequations}
	and $a = R^\transp\partial_1\phi \cdot e_1$ and $\varphi \in \Leb^2(\omega, \SobH^1_{\rm av}(S,\R^3))$ is such that
	$\sym\compose{0}{\bar\grad\varphi} = \sym\compose{0}{\bar\grad\tilde\varphi}$, where $\tilde\varphi = R^\transp \psi + (R^\transp \partial_1 \phi \cdot \mathbf{\bar{x}})e_1$.
\end{enumerate}
\end{proposition}

\begin{remark} \label{Rem:Third-gradient-term-compactness}
The only estimate in \cref{Prop:Compactness} where the strain gradient term $H_\el^\mathbb{F}$ is relevant is the estimate $h^{\alpha_\mathbb{F}} \norm*{\grad^2\grad[h]y^h}_{\Leb^{p_*}(\Omega)} \leq \tilde M$. Using the Sobolev inequality, we conclude from this estimate that $(\grad\grad[h]y^h)$ and $(\grad C_{y^h})$ are bounded in $\Leb^\infty(\Omega)$ for fixed $h > 0$.
\end{remark}

\begin{remark} \label{Rem:CosseratAnsatzRefined}
Let $y \in \mathcal{Q}^h_{\rm bc}$. One can associate to $y$ a rod configuration $q = q^h[y] \in \mathcal{Q}^0_{\rm bc}$ by letting $(v,R,\phi) := (v^h[y], R^h[y], \phi^h[y])$ as in \Cref{Rem:CosseratAnsatz} and defining
\begin{equation}
	a^h[y] := \big( R^\transp{\textstyle\fint_S} \partial_1 \phi \big) \cdot e_1, \qquad
	\varphi^h[y] := h^{-1} R^\transp\tilde\phi + \left(\big( R^\transp{\textstyle\fint_S} \partial_1 \phi \big) \cdot \mathbf{\bar x}\right) e_1,
\end{equation}
where $\tilde\phi(x_1,\placeholder) \in \SobH^1_{\rm av}(S,\R^3)$ is such that $\sym\compose{0}{\bar\grad\phi} = \sym\compose{0}{\bar\grad\tilde\phi}$, that is $\tilde\phi := \phi - \fint_S \phi - \frac12\left( \fint_S \partial_3\phi_2 - \partial_2\phi_3 \right) (0,x_3,-x_2)^\transp$. We note this construction depends on the choice of the rotation field $R$. However, this choice is only up to an error that vanishes as $h \to 0$ as \cref{Prop:Compactness} shows.
This definition is motivated by the observation that in the situation of \Cref{Prop:Compactness},
\begin{equation}
	a^h[y^h] \wto a \text{ in } \Leb^2(\omega), \qquad
	\varphi^h[y^h] \wto \varphi \text{ in } \Leb^2(\omega, \SobH^1(S, \R^3)).
\end{equation}
On the other hand, one may associate with $q = (v,R,a,\varphi) \in \mathcal{Q}^0_{\rm bc}$ a deformation $y^h[q]$ via the corrected Cosserat ansatz
\begin{equation} \label{Eq:CorrectedCosseratAnsatz}
	y^h[q] := v + hR\mathbf{\bar x} + h\int_0^{x_1} a(s) R(s)e_1\,\dd s + h^2 R\varphi.
\end{equation}
We observe that in general we do not have $q^h[y^h[q]] = q$ or $y^h[q^h[y]] = y$ but there is the following connection. Recall that in elasticity an important property when passing to the limit is the strain tensor $E^h := h^{-1} (\sqrt{C_{y^h}} - I)$, which can be approximated by $G^h := h^{-1} ((R^h)^\transp\grad[h]y^h - I)$. In fact, it is well-known that limits of $(E^h)$ coincide with limits of $(\sym G^h)$ (see \cite[Lem.~4.5]{Neu12}) and are characterized by limiting states $q$ via
\begin{equation}
	G[q] := \compose[\big]{K_R \mathbf{\bar x} + ae_1}{\bar\grad\varphi}, \qquad
	q \in \mathcal Q^0_{\rm bc},
\end{equation}
cf.\ e.g.\ \cite[Prop.~4.6]{NR26}. This is captured by our construction, as a direct calculation shows
\begin{equation}
	\sym G^h[y^h] = \sym G\big[ q^h[y^h] \big] + \sym\compose[\big]{(R^h)^\transp( \partial_1 \phi^h - {\textstyle\fint_S}\partial_1\phi^h)}{0}
\end{equation}
and the latter term vanishes as $h \to 0$.
\end{remark}

\section{Curves of maximal slope in the three dimensional model}
\label{Sec:3D-model}

In this section we discuss the finite thickness model \eqref{Eq:Intro:PDE} and prove the following existence result for curves of maximal slope for $(\mathcal{Q}^h_\mathrm{bc}, \mathcal{D}^h, \mathcal{E}^h)$:

\begin{theorem}[Existence of curves of maximal slope]
\label{Thm:existence-bulk}
Let $(y^h_0) \subset \mathcal{Q}^h_\mathrm{bc}$ be a sequence of initial values such that
\begin{equation*}
	\limsup_{h \to 0} \mathcal{E}^h(0, y^h_0) < \infty.
\end{equation*}
Then, there exists $h_0 > 0$ such that for all $0 < h \leq h_0$ there exists a curve of maximal slope $y^h: [0,T] \to \mathcal{Q}^h_\mathrm{bc}$ for $(\mathcal{Q}^h_\mathrm{bc}, \mathcal{D}^h, \mathcal{E}^h)$ with $y^h(0) = y^h_0$.
\end{theorem}

Below, in \cref{Sec:Existence-of-solutions} we prove this existence result. Furthermore, in \cref{Sec:Characterization-3D} we prove that the weak formulation \eqref{Eq:Intro:PDE} and the metric gradient flow formulation are equivalent.
But before, let us discuss some basic properties of the metric gradient flow. We start with the fact that $\mathcal{D}^h$ is indeed a metric on $\mathcal{Q}^h_{\rm bc}$.

\begin{lemma}
For all $h > 0$, $\mathcal{D}^h$ is a metric on $\mathcal{Q}^h_\mathrm{bc}$.
\end{lemma}

\begin{proof}
Due to the assumption \eqref{Ass:DV_Pseudometric}, symmetry of $\mathcal{D}^h$ and the triangle inequality are trivial. It remains to show positive definiteness of $\mathcal{D}^h$. We show this as a consequence of the version of the geometric rigidity inequality or nonlinear Korn inequality in \cite{CM15}. 
Let $y, \hat{y} \in \mathcal{Q}^h_\mathrm{bc}$ with $\mathcal{D}^h(y, \hat{y}) = 0$. Then, by \eqref{Ass:DV_Pseudometric} we have
\begin{equation*}
	\grad[h]y^\transp \grad[h]y = \grad[h]\hat y^\transp \grad[h]\hat y \qquad \text{a.e.\ in } \Omega.
\end{equation*}
We use that $\mathcal{Q}^h_\mathrm{bc} \subset \set{\tilde{y} \in \Cont^1(\overline{\Omega}, \R^3) \given \det\grad[h]\tilde{y} > 0 \text{ in } \Omega}$. However, in general we do not have $\det\grad[h]y, \det\grad[h]\hat{y} > 0$ in $\overline{\Omega}$ so we cannot apply the results of \cite{CM15} directly. Covering the open, connected set $\Omega$ with overlapping balls that are compactly contained in $\Omega$ we obtain from \cite[Thm.~1]{CM15} applied on these balls that there exists a fixed rotation $R \in \SO(3)$ and a constant $c \in \R^3$, such that
\begin{equation*}
	\hat{y} = Ry + c
	\quad\text{in } \overline{\Omega}.
\end{equation*}
The boundary conditions imply $R_{\rm bc}\mathbf{\bar x} = R R_{\rm bc}\mathbf{\bar x}$ in $S$ which yields $R = I$ and then, $c = (I - R)v_{\rm bc} = 0$.
\end{proof}

Let us further establish a control of the power $\partial_t \mathcal{E}^h$ that is essential to derive bounds for the energy along curves of maximal slope.

\begin{lemma}\label{Lem:Regularity-energy}
For all $h_0 > 0$ the following statements hold:
\begin{enumerate}[(a)]
	\item\label{Item:Regularity-energy:compactness}
	There exists a constant $c > 0$, such that for all $0 < h \leq h_0$ and $y \in \mathcal Q^h_{\rm bc}$, we have
	\begin{equation}
		\max\set*{ \sup_{t \in [0,T]} \mathcal E^h(t,y),\; \mathcal E^h_\el(y) + \mathcal E^h_{\rm sg}(y) }
		\leq c \left( \inf_{t \in [0,T]} \mathcal E^h(t,y) + c \right).
	\end{equation}
	Especially, the energy is bounded from below.
	
	\item\label{Item:Regularity-energy:time-regularity}
	For all $0 < h \leq h_0$ the map $t \mapsto \mathcal E^h(t,y)$ is absolutely continuous for all $y \in \mathcal Q^h_{\rm bc}$ and the set $\Lambda \subset [0,T]$ where $\partial_t \mathcal E^h(t,y)$ exists is independent of $y$ and $h$. Moreover, there exists a non-negative map $c \in \Leb^1([0,T])$ (uniform for $0 < h \leq h_0$), such that for all $t \in \Lambda$ and $y \in \mathcal Q^h_{\rm bc}$,
	\begin{equation}
		\abs*{\partial_t \mathcal E^h(t,y)} \leq c(t) \left( \mathcal E^h_\el(y) + 1 \right).
	\end{equation}
\end{enumerate}
\end{lemma}

\begin{proof}
The proof follows analogously to \cite[Lem.~A.3]{FM23} from the estimate
\begin{equation*}
	\norm{y}_{\SobH^1(\Omega)} \leq c_1 \bigl( \mathcal{E}^h_\el(y)^{1/2} + 1 \bigr)
\end{equation*}
that holds due to \eqref{Ass:W_NonDegeneracy} and the boundary conditions. For the readers convenience we repeat the arguments.
$l \in \AC([0,T]; \SobH^1(\Omega, \R^3)')$ implies that for all $t \in [0,T]$,
\begin{equation*}
	\mathcal E^h_\el(y) + \mathcal E^h_{\rm sg}(y)
	= \mathcal{E}^h(t,y) + \sprod{l(t)}{y}
	\leq \mathcal{E}^h(t,y) + c_2 \norm{y}_{\SobH^1(\Omega)}
	\leq \mathcal{E}^h(t,y) + c_3 \bigl( \mathcal{E}^h_\el(y)^{1/2} + 1 \bigr).
\end{equation*}
The term $c_3 \mathcal{E}^h_\el(y)^{1/2}$ on the right-hand side can be absorbed into the left-hand side by an $\eps$-Young inequality to show $\mathcal E^h_\el(y) + \mathcal E^h_{\rm sg}(y) \leq c_4 (\mathcal{E}^h(t,y) + c_4)$. Furthermore, for all $s \in [0,T]$, we then obtain
\begin{equation*}
	\mathcal{E}^h(s,y)
	= \mathcal{E}^h(t,y) + \sprod{l(t) - l(s)}{y}
	\leq \mathcal{E}^h(t,y) + c_5 \norm{y}_{\SobH^1(\Omega)}
	\leq c_6 (\mathcal{E}^h(t,y) + c_6).
\end{equation*}
From this we conclude (\ref{Item:Regularity-energy:compactness}). Finally, we conclude (\ref{Item:Regularity-energy:time-regularity}) from
\begin{equation*}
	\abs{\partial_t \mathcal{E}^h(t,y)}
	= \abs{\sprod{\dot{l}(t)}{y}}
	\leq \norm{\dot{l}(t)}_{\SobH^1(\Omega)'} \norm{y}_{\SobH^1(\Omega)}
	\leq c_7\norm{\dot{l}(t)}_{\SobH^1(\Omega)'} \bigl( \mathcal{E}^h_\el(y) + 1 \bigr).
\qedhere
\end{equation*}
\end{proof}

\subsection{Existence of curves of maximal slope}
\label{Sec:Existence-of-solutions}

\Cref{Thm:existence-bulk} follows essentially from the abstract existence result in \cite{RMS08} for which we check the requirements in this section. The existence proof is based on the related time-incremental scheme given by the initial state $y^h_0 \in \mathcal{Q}^h_{\rm bc}$ and the iterative definition,
\begin{equation}
\label{Eq:time-incremental-bulk}
	\left\{\begin{aligned}
		t_i &:= t_{i - 1} + \tau, \\
		y^h_i &:\in \operatorname{Arg\,min}_{y \in \mathcal{Q}^h_{\rm bc}} \bigl(\mathcal E^h(t_i,y) + \tfrac{1}{2\tau} \mathcal{D}^h(y^h_{i-1}, y)^2\bigr).
	\end{aligned}\right.
\end{equation}
Important for the analysis in \cite{RMS08} is also the de~Giorgi interpolant $\tilde{y}^h: [0,T] \to \mathcal{Q}^h_{\rm bc}$ given for $t \in (t_{i-1}, t_i)$ by
\begin{equation}
\label{Eq:de-Giorgi-interpolants-bulk}
	\tilde{y}^h(t) :\in \operatorname{Arg\,min}_{y \in \mathcal{Q}^h_{\rm bc}} \bigl(\mathcal E^h(t,y) + \tfrac{1}{2(t - t_{i-1})} \mathcal{D}^h(y^h_{i-1}, y)^2\bigr).
\end{equation}
The control on the power obtained in \cref{Lem:Regularity-energy} is essential to derive uniform boundedness of the energies $\mathcal{E}^h(t, \tilde{y}^h)$ which allows us to use the a~priori estimates from \cref{Prop:Compactness} to establish compactness properties for the states as $\tau \to 0$. The control obtained from \Cref{Lem:Regularity-energy} is a slightly more general than the property \cite[Eq.~(2.19d)]{RMS08} but it is not hard to check that this assumption suffices.
We note that in \cite{FG17} existence of curves of maximal slope is discussed for time-dependent energies with time-regularity assumptions that are in line with \cref{Lem:Regularity-energy} but the convexity assumptions are too strong for our setting.
\medskip

As \cref{Lem:Regularity-energy} establishes estimates that are uniform for $h \ll 1$ and in \Cref{Thm:existence-bulk} we assume that the initial energy is bounded, it is not hard to show that (approximate) minimizers for \eqref{Eq:time-incremental-bulk} and \eqref{Eq:de-Giorgi-interpolants-bulk} are contained in a sublevel $\mathcal{Q}^{h,M}_{\rm bc}$ (recall \eqref{Eq:Intro:state-space-bounded-energy}) where $M > 0$ is independent of $h$ (for $0 < h \leq h_0$). Any further consideration in \cite{RMS08} can be restricted to such an energy sublevel. This observation is crucial in order to deal with the delicate smallness assumption on $h$ in \Cref{Thm:existence-bulk} which is due to the fact that we only assume regularity assumptions on the energy and dissipation distance in a neighborhood of $\SO(3)$.
\medskip

Incorporating these observations, specified to our setting, the result in \cite{RMS08} requires the following properties to be fullfiled (with $M > 0$ as discussed above):

\begin{itemize}
	\item (Compactness, cf.\ \cite[Eq.~(2.38)]{RMS08}):
	$\mathcal{Q}^{h,M}_{\rm bc}$ is (sequentially) compact with respect to the weak topology of $\SobW^{3,p_*}(\Omega, \R^3)$.
	
	\item (Lower semi-continuity, cf.\ \cite[Eq.~(2.34)]{RMS08}):
	Whenever $(y_k), (\hat{y}_k) \subset \mathcal{Q}^{h,M}_{\rm bc}$ satisfy $y_k \wto y$ and $\hat{y}_k \wto \hat{y}$ in $\SobW^{3,p_*}(\Omega, \R^3)$, then
	\begin{equation*}
		\liminf_{k \to \infty} \mathcal{E}^h(t,y_k) \geq \mathcal{E}^h(t,y), \qquad
		\liminf_{k \to \infty} \mathcal{D}^h(y_k, \hat{y}_k) \geq \mathcal{D}^h(y,\hat{y}).
	\end{equation*}
	
	\item (Upper semi-continuity of the power, cf.\ \cite[Eq.~(2.35)]{RMS08}):
	For all $t \in \Lambda$ (with $\Lambda$ as in \cref{Lem:Regularity-energy}) and $(y_k) \subset \mathcal{Q}^{h,M}_{\rm bc}$ such that $y_k \wto y$ in $\SobW^{3,p_*}(\Omega, \R^3)$, we have
	\begin{equation*}
		\limsup_{k \to \infty} \partial_t\mathcal{E}^h(t,y_k) \leq \partial_t\mathcal{E}^h(t,y).
	\end{equation*}
	
	\item (Lower semi-continuity of the local slope, cf.\ \cite[Eq.~(2.18)]{RMS08})
	Whenever $(y_k) \subset \mathcal{Q}^{h,M}_{\rm bc}$ satisfies $y_k \wto y$ in $\SobW^{3,p_*}(\Omega, \R^3)$, then
	\begin{equation*}
		\liminf_{k \to \infty} \abs[\big]{\partial\mathcal{E}^h(t)}_{\mathcal{D}^h}(y_k) \geq \abs[\big]{\partial\mathcal{E}^h(t)}_{\mathcal{D}^h}(y).
	\end{equation*}
	
	\item (Chain rule, cf.\ \cite[Def.~2.5]{RMS08}):
	For all $y \in \AC([0,T]; \mathcal{Q}^h_{\rm bc})$ with $y(t) \in \mathcal{Q}^{h,M}_{\rm bc}$ for all $t \in [0,T]$ and $\int_0^T \abs{y'}_{\mathcal D^h}(t) \abs{\partial\mathcal{E}^h(t)}_{\mathcal{D}^h}(y(t)) \,\dd t < \infty$, the map $t \mapsto \mathcal{E}^h(t,y(t))$ is absolutely continuous and
	\begin{equation*}
		\tfrac{\dd}{\dd t} \mathcal E^h\of[\big]{t,y(t)} \geq \partial_t \mathcal E^h\of[\big]{t,y(t)} - \abs*{y'}_{\mathcal D^h}(t)\abs*{\partial \mathcal E^h(t)}_{\mathcal D^h}\of[\big]{y(t)}
		\qquad \text{for a.e.\ } t \in [0,T].
	\end{equation*}
\end{itemize}

Let us establish the first three of these properties in the following lemma and prepare for the remaining two. For convenience let us introduce two further version of the energy sublevels:
\begin{subequations}
\begin{align}
	\mathcal Q^{h,M-}_{\rm bc} &:= \set[\Big]{y \in \mathcal Q^h_{\rm bc} \given \inf_{t \in [0,T]} \mathcal E^h(t,y) \leq M},
	\\
	\mathcal Q^{h,<M}_{\rm bc} &:= \set[\Big]{y \in \mathcal Q^h_{\rm bc} \given \sup_{t \in [0,T]} \mathcal E^h(t,y) < M}.
\end{align}
\end{subequations}

\begin{lemma}[Properties of the metric gradient flow]
\label{Lem:Properties-MGF-bulk}
Let $M > 0$. There exists $h_0 > 0$, such that for all $0 < h \leq h_0$ the following statements hold.
\begin{enumerate}[(a)]
	\item (Compactness):
	\label{Item:Properties-MFG-bulk:compactness}%
	$\mathcal{Q}^{h,M-}_{\rm bc}$ and $\mathcal{Q}^{h,M}_{\rm bc}$ are (sequentially) compact with respect to the weak topology of $\SobW^{3,p_*}(\Omega, \R^3)$.
	
	\item (Lower semi-continuity of the energy):
	\label{Item:Properties-MFG-bulk:LSC-energy}%
	For all $t \in [0,T]$, $y \in \mathcal{Q}^h_{\rm bc}$ and sequences $(y_k) \subset \mathcal{Q}^{h,M-}_{\rm bc}$, we have
	\begin{equation}
		y_k \wto y \text{ in } \SobW^{3,p_*}(\Omega, \R^3)
		\quad\Rightarrow\quad
		\liminf_{k \to \infty} \mathcal{E}^h(t,y_k) \geq \mathcal{E}^h(t,y).
	\end{equation}
	
	\item (Continuity of the dissipation distance):
	\label{Item:Properties-MFG-bulk:cont-distance}%
	For all $y, \hat{y} \in \mathcal{Q}^h_{\rm bc}$ and sequences $(y_k), (\hat{y}_k) \subset \mathcal{Q}^{h,M-}_{\rm bc}$, we have
	\begin{subequations}
	\begin{align}
		y_k \wto y,\, \hat{y}_k \wto \hat{y} \text{ in } \SobW^{3,p_*}(\Omega, \R^3)
		&\quad\Rightarrow\quad
		\lim_{k \to \infty} \mathcal{D}^h(y_k, \hat{y}_k) = \mathcal{D}^h(y,\hat{y}),
		\\
		\lim_{k \to \infty} \mathcal{D}^h(y_k, y) = 0
		&\quad\Rightarrow\quad
		y_k \wto y \text{ in } \SobW^{3,p_*}(\Omega, \R^3).
	\end{align}
	\end{subequations}
	Especially, we obtain that the topology induced by $\mathcal{D}^h$ coincides on $\mathcal{Q}^{h,M-}_{\rm bc}$ with the weak topology of $\SobW^{3,p_*}(\Omega, \R^3)$.
	
	\item (Continuity of the power):
	\label{Item:Properties-MFG-bulk:cont-power}%
	For all $t \in \Lambda]$, $y \in \mathcal{Q}^h_{\rm bc}$ and sequences $(y_k) \subset \mathcal{Q}^{h,M-}_{\rm bc}$, we have
	\begin{equation}
		y_k \wto y \text{ in } \SobW^{3,p_*}(\Omega, \R^3)
		\quad\Rightarrow\quad
		\lim_{k \to \infty} \partial_t \mathcal{E}^h(t,y_k) = \partial_t \mathcal{E}^h(t,y).
	\end{equation}
	
	\item (Estimate for the dissipation distance):
	\label{Item:Properties-MFG-bulk:equivalence-distance}%
	There exist constants $\delta^h,c^h > 0$, such that for all $y, \hat{y} \in \mathcal{Q}^{h,M-}_{\rm bc}$ with $\norm{\grad[h]\hat y - \grad[h]y}_{\Leb^\infty(\Omega)} \leq \delta^h$, we have
	\begin{equation}
		\tfrac{1}{c^h}\norm{\grad[h]\hat y - \grad[h]y}_{\Leb^2(\Omega)}
		\leq \mathcal{D}^h\of{y, \hat y}
		\leq c^h\norm{\grad[h]\hat y - \grad[h]y}_{\Leb^2(\Omega)}.
	\end{equation}
	
	\item (Semiconvexity):
	\label{Item:Properties-MFG-bulk:convexity}%
	There exist $\delta^h, \lambda^h, c^h > 0$, such that for all $t \in [0,T]$ and $y,\hat y \in \mathcal{Q}^{h,M-}_{\rm bc}$ with $0 < \norm{\grad[h]\hat y - \grad[h]y}_{\Leb^\infty(\Omega)} \leq \delta^h$ and $s \in [0,1]$,
	$y_s := y + s(\hat y - y)$, satisfies
	\begin{subequations}
	\begin{align}
		\mathcal E^h(t, y_s)
		&\leq \mathcal E^h(t,y) + s \left( \mathcal E^h(t,\hat y) - \mathcal E^h(t,y) \right)
		+ \lambda^h s(1-s) \mathcal{D}^h(y,\hat y)\norm*{\grad[h]y - \grad[h]\hat y}_{\SobW^{1,\infty}(\Omega)},
		\label{Eq:Lambda-convexity-energy:Energy}
		\\
		\mathcal{D}^h(y,y_s) 
		&\leq s\mathcal{D}^h(y,\hat y)(1 + c^h\norm*{\grad[h]\hat y - \grad[h]y}_{\Leb^\infty(\Omega)})^{1/2}.
		\label{Eq:Lambda-convexity-energy:Dissipation}
	\end{align}
	\end{subequations}
\end{enumerate}
\end{lemma}

\begin{proof}
\stepemph{Step 1 -- Proof of (\ref{Item:Properties-MFG-bulk:compactness})--(\ref{Item:Properties-MFG-bulk:cont-power}):}
The proof of (\ref{Item:Properties-MFG-bulk:compactness})--(\ref{Item:Properties-MFG-bulk:cont-power}) follow analogously to the procedure in \cite[Sec.~4]{FK18}, \cite[Sec.~4]{FK20}, \cite[Sec.~5]{FM23} and \cite[Sec.~4]{Mac25} so that we refrain from a detailed proof here and only comment on the main ideas.
\Cref{Prop:Compactness,Lem:Regularity-energy} imply that $\mathcal{Q}^{h,M-}_{\rm bc}$ is a bounded subset of $\SobW^{3,p_*}(\Omega, \R^3)$, so that $\mathcal{Q}^{h,M-}_{\rm bc}$ and $\mathcal{Q}^{h,M}_{\rm bc}$ are precompact w.r.t.\ the weak topology of $\SobW^{3,p_*}(\Omega, \R^3)$. Furthermore, these sets are closed which follows from \Cref{Lem:Regularity-energy} and the lower semi-continuity of the energy in (\ref{Item:Properties-MFG-bulk:LSC-energy}) which itself is a consequence of the convexity \eqref{Ass:HE_Convexity} and $p_*$-growth \eqref{Ass:HE_Growth} of the third gradient term $H^{\mathbb F}_\el$. Note that we also use here that $h_0 > 0$ can be chosen sufficiently small in view of \Cref{Prop:Compactness}, so that $\grad[h]y \in U$ for all $y \in \mathcal{Q}^{h,M-}_{\rm bc}$ where $U$ is a neighborhood of $\SO(3)$ as in \eqref{Ass:W_regularity}. From this observation also (\ref{Item:Properties-MFG-bulk:cont-distance}) follows as a consequence of \eqref{Ass:DV_Regularity} and the compactness (\ref{Item:Properties-MFG-bulk:compactness}). Finally, (\ref{Item:Properties-MFG-bulk:cont-power}) follows easily from $\partial_t \mathcal{E}^h(t,y) = -\sprod{l(t)}{y}$.
\medskip

\stepemph{Step 2 -- Proof of (\ref{Item:Properties-MFG-bulk:equivalence-distance}):}
Let $y, \hat{y} \in \mathcal{Q}^{h,M-}_{\rm bc}$ with $\norm*{\grad[h]\hat y - \grad[h]y}_{\Leb^\infty(\Omega)} \leq \delta^h$, where $0 < \delta^h \ll 1$ is chosen below. Then, by \eqref{Eq:CompactnessStrainBound}, $\grad[h]y, \grad[h]\hat y$ are contained in a uniform neighborhood of $\SO(3)$ for $h_0 > 0$ small enough. Recall the definition of $\mathcal D^h$ in \eqref{Eq:Intro:Dissipation}. Thus, we may apply \eqref{Eq:TaylorDissipation2} and obtain
\begin{align*}
	\mathcal{D}^h(y,\hat y)^2
	&\leq h^{-2}\int_\Omega \tfrac{1}{2} \dif_{F_1}^2 \mathscr{D}_\vi^2(\grad[h]y,\grad[h]y)[\grad[h]\hat y - \grad[h]y,\grad[h]\hat y - \grad[h]y] + c_1\, \abs*{\grad[h] \hat y - \grad[h]y}^3 \\
	&\leq c_2 h^{-2} \norm*{\grad[h]\hat y - \grad[h]y}_{\Leb^2(\Omega)}^2,
\end{align*}
where the second inequality follows from the boundedness of $\dif_{F_1}^2\mathscr D_\vi^2$.
On the other hand, \eqref{Eq:NonDegeneracy_QviAtF} yields
\begin{align*}
	\mathcal{D}^h(y,\hat y)^2
	&\geq h^{-2}\int_\Omega \tfrac{1}{2} \dif_{F_1}^2 \mathscr{D}_\vi^2(\grad[h]y,\grad[h]y)[\grad[h]\hat y - \grad[h]y,\grad[h]\hat y - \grad[h]y] - c_1\, \abs*{\grad[h] \hat y - \grad[h]y}^3 \\
	&\geq h^{-2} \left(c_3\norm*{\sym\of*{(\grad[h]y)^\transp(\grad[h]\hat y - \grad[h]y)}}_{\Leb^2(\Omega)}^2 - c_1\delta^h \norm*{\grad[h]\hat y - \grad[h]y}_{\Leb^2(\Omega)}^2 \right) \\
	&\geq h^{-2}(c_4 - c_1\delta^h)\norm*{\grad[h]\hat y - \grad[h]y}_{\Leb^2(\Omega)}^2.
\end{align*}
Here, the third inequality follows from the generalized Korn inequality of \cite{Pom03} (see also \cite[Thm.~3.3]{MR20} and \cite[Thm.~5.3]{FM23}) which holds as $\grad[h]y^\transp$ is uniformly close to $\SO(3)$ with a constant that does not depend on $y$. Choosing $\delta^h < \frac{c_4}{c_1}$, we obtain (\ref{Item:Properties-MFG-bulk:equivalence-distance}). Note that the constant $c_4$ may depend on $h$ as observed in \cite[Lem.~5.2]{FM23}.
\medskip

\stepemph{Step 3 -- Proof of (\ref{Item:Properties-MFG-bulk:convexity}):}
Finally, we show the semiconvexity statement (\ref{Item:Properties-MFG-bulk:convexity}). Let $\delta^h > 0$ be as in Step 2. Let us first consider the strain gradient energy $\mathcal E^h_{\rm sg}$. As $H_\el^F$ and $H_\el^{\mathbb F}$ are convex, it suffices to estimate the term $\mathcal I(y) := h^{-\alpha_C p_*} \int_\Omega H_\el^C(\grad C_y)$. It is straight-forward to derive the formula,
\begin{equation*}
	C_{y_s} = C_y + s(C_{\hat y} - C_y) - s(1-s)(\grad[h]\hat y - \grad[h]y)^\transp(\grad[h]\hat y - \grad[h]y).
\end{equation*}
Then, using \eqref{Ass:HE_Continuity}, we obtain,
\begin{align*}
	&\abs*{H_\el^C\of[\big]{\grad C_{y_s}} - H_\el^C\of[\big]{\grad C_y + s(\grad C_{\hat y} - \grad C_y)}} \\
	&\qquad\leq c_5 s(1-s)\left( \abs[\big]{\grad C_{y_s}}^{p_*-1} + \abs[\big]{\grad C_y}^{p_*-1} + \abs[\big]{\grad C_{\hat y}}^{p_*-1} \right) \abs[\big]{\grad[h]\hat y - \grad[h]y}\abs[\big]{\grad\grad[h]\hat y - \grad\grad[h]y} \\
	&\qquad\leq c_6 s(1-s) \abs[\big]{\grad[h]\hat y - \grad[h]y}\abs[\big]{\grad\grad[h]\hat y - \grad\grad[h]y},
\end{align*}
where in the second inequality, we used that due to \cref{Prop:Compactness}, $\grad[h]y$, $\grad[h]\hat y$, $\grad\grad[h]y$ and $\grad\grad[h]\hat y$ are bounded. Note that the constant $c_6$ is allowed to depend on $h$. Hence, we may apply Step 2 and use the convexity of $H_\el^C$ to obtain,
\begin{align*}
	\mathcal I(y_s) 
	&\leq \mathcal I(y) + s \left(\mathcal I(\hat y) - \mathcal I(y)\right)
	+ c_7 s(1-s) \norm*{\grad[h]\hat y - \grad[h]y}_{\Leb^2(\Omega)}\norm*{\grad\grad[h]\hat y - \grad\grad[h]y}_{\Leb^2(\Omega)}
	\\
	&\leq \mathcal I(y) + s \left(\mathcal I(\hat y) - \mathcal I(y)\right)
	+ c_8 s(1-s) \mathcal D^h(y,\hat y)\norm*{\grad[h]\hat y - \grad[h]y}_{\SobW^{1,\infty}(\Omega)}.
\end{align*}
Hence, by convexity of $H_\el^F$ and $H_\el^{\mathbb F}$ we obtain,
\begin{equation*}
	\mathcal E^h_{\rm sg}(y_s)
	\leq \mathcal E^h_{\rm sg}(y) + s \left( \mathcal E^h_{\rm sg}(\hat y) - \mathcal E^h_{\rm sg}(y)\right)
	+ c_8 s(1-s) \mathcal D^h(y,\hat y)\norm*{\grad[h]\hat y - \grad[h]y}_{\SobW^{1,\infty}(\Omega)}.
\end{equation*}
In order to proceed with the estimates for the elastic energy, we first note that for any map $f \in \Cont^2(U)$, $U \subset \R^n$ open, a standard Taylor expansion argument shows, that whenever $x_s := x + s(\hat x - x) \in U$ for all $s \in [0,1]$,
\begin{equation*}
	f(x_s)
	= f(x) + s\left( f(\hat x) - f(x)\right) - s \int_0^1 \int_{st}^t \dif^2 f(x_r)[\hat x - x, \hat x - x] \,\dd r \dd t.
\end{equation*}
Recall \eqref{Ass:W_regularity}. For $0< h_0 \ll 1$ and $0 < \delta^h \ll 1$ small enough, $\grad[h]y$, $\grad[h]\hat y$ and $\grad[h]y_s$ are contained in the uniform neighborhood $U$ of $\SO(3)$. Thus, boundedness of $\dif^2_F W_\el$ in $U$ implies the semi-convexity,
\begin{equation*}
	\mathcal E^h_\el(t,y_s) \leq \mathcal E^h_\el(t,y) + s\left(\mathcal E^h_\el(t,\hat y) - \mathcal E^h_\el(t,y)\right) + c_9 s(1-s) \norm{\grad[h]\hat y - \grad[h]y}_{\Leb^2(\Omega)}^2.
\end{equation*}
Hence, from the linearity of the loading term $y \mapsto \sprod{l(t)}{y}$, we conclude \eqref{Eq:Lambda-convexity-energy:Energy}.
Finally, for the dissipation distance, we infer from \eqref{Eq:TaylorDissipation2},
\begin{align*}
	\mathscr{D}_\vi^2\of[\big]{\grad[h]y, \grad[h]y_s}
	&\leq \tfrac{s^2}{2} \dif_{F_1}^2\mathscr{D}_\vi^2\of[\big]{\grad[h]y, \grad[h]y}[\grad[h]\hat y - \grad[h]y, \grad[h]\hat y - \grad[h]y]  + c_{10}s^3 \abs*{\grad[h]\hat y - \grad[h]y}^3,
	\\
	\mathscr{D}_\vi^2\of[\big]{\grad[h]y, \grad[h]y}
	&\geq \tfrac{\mathrlap{\;1}\phantom{s^2}}{2} \dif_{F_1}^2\mathscr{D}_\vi^2\of[\big]{\grad[h]y, \grad[h]y}[\grad[h]\hat y - \grad[h]y, \grad[h]\hat y - \grad[h]y] - c_{11} \abs{\grad[h]\hat y - \grad[h]y}^3.
\end{align*}
Combining these statements, we obtain,
\begin{equation*}
	\mathcal D^h(y, y_s)^2
	\leq s^2\mathcal D^h(y,\hat y) + c_{12} s^2 \norm{\grad[h]\hat y - \grad[h]y}_{\Leb^2(\Omega)}^2 \norm*{\grad[h]\hat y - \grad[h]y}_{\Leb^\infty(\Omega)}.
\end{equation*}
Then, we infer \eqref{Eq:Lambda-convexity-energy:Dissipation} from (\ref{Item:Properties-MFG-bulk:equivalence-distance}).
\end{proof}

The semi-convexity statement of \Cref{Lem:Properties-MGF-bulk} that is in the spirit of \cite[Sec.~2.4]{AGS05} is a very useful property in order to study the local slope. One consequence is a characterization of the local slope as in \cite[Thm.~2.4.9]{AGS05}, see also \cite[Lem.~4.9]{FK18}, \cite[Lem.~4.8]{FK20} and \cite[Lem.~5.4]{FM23}. Eventually, we obtain the lower semi-continuity of the local slope and the chain rule property from this characterization. We note that the chain rule property is related to the local slope being a strong upper gradient in the terminology of \cite{AGS05}.

\begin{corollary}
\label{Cor:Local-slope-3D:Characterization}
Let $M > 0$. There exists $h_0 > 0$ such that for all $0 < h < h_0$ the following statement holds. There exist $\delta^h, \lambda^h, c^h > 0$ such that for all $t \in [0,T]$ $y \in \mathcal Q^{h,M}_{\rm bc}$, we have the representation
\begin{equation} \label{Eq:Local-slope-3D:Characterization}
	\abs*{\partial \mathcal{E}^h(t)}_{\mathcal{D}^h}(y) 
	= \sup_{\hat y \in \mathcal A(h,M,y,\delta^h)} \frac{\left( \mathcal{E}^h(t,y) - \mathcal{E}^h(t,\hat y) - \lambda^h\mathcal{D}^h(y,\hat y)\norm*{\grad[h]y - \grad[h]\hat y}_{\SobW^{1,\infty}(\Omega)}\right)^+}{\mathcal{D}^h(y,\hat y)(1 + c^h\norm*{\grad[h]\hat y - \grad[h]y}_{\Leb^\infty(\Omega)})^{1/2}},
\end{equation}
where $\mathcal A(h,M,y,\delta^h) := \set*{\hat y \in \mathcal{Q}^{h,M-}_{\rm bc} \given \hat y \neq y,\, 0 < \norm*{\grad[h]\hat y - \grad[h]y}_{\Leb^\infty(\Omega)} < \delta^h}$.
\end{corollary}

\begin{proof}
Let $h_0 > 0$ be as in \Cref{Lem:Properties-MGF-bulk} and $\delta^h, \lambda^h, c^h > 0$ as in \Cref{Lem:Properties-MGF-bulk}~(\ref{Item:Properties-MFG-bulk:convexity}). Fix $y \in \mathcal{Q}^{h,M}_{\rm bc}$ and $t \in [0,T]$.
\medskip

\stepemph{Step 1 -- Upper bound for \eqref{Eq:Local-slope-3D:Characterization}:} 
We first show that the right-hand side is an upper bound for the local slope for any choice of $\delta^h, \lambda^h, c^h > 0$. Let therefore $(y_k) \subset \mathcal{Q}^h_{\rm bc}$ be a sequence, such that $\mathcal{D}^h(y,y_k) \to 0$ and
\begin{equation*}
	\abs*{\partial \mathcal{E}^h(t)}_{\mathcal{D}^h}(y)
	= \lim_{k \to \infty} \frac{\left( \mathcal{E}^h(t,y) - \mathcal{E}^h(t,y_k) \right)^+}{\mathcal{D}^h(y, y_k)}
\end{equation*}
and assume without loss of generality that $\abs*{\partial \mathcal{E}^h(t)}_{\mathcal{D}^h}(y) > 0$. Then, we observe that necessarily $\mathcal E^h(t,y_k) < \mathcal E^h(t,y)$ and therefore $y_k \in \mathcal Q^{h,M-}_{\rm bc}$ for $k \gg 1$, as otherwise $(\mathcal{E}^h(t,y) - \mathcal{E}^h(t,y_k))^+ = 0$. Hence, by \Cref{Lem:Properties-MGF-bulk}~(\ref{Item:Properties-MFG-bulk:cont-distance}), $y_k \wto y$ in $\SobW^{3,p_*}(\Omega, \R^3)$ and thus,
\begin{align*}
	\abs*{\partial \mathcal{E}^h(t)}_{\mathcal{D}^h}(y)
	&= \lim_{k \to \infty} \frac{\left( \mathcal{E}^h(t,y) - \mathcal{E}^h(t,y_k) \right)^+}{\mathcal{D}^h(y,y_k)} \\
	&= \lim_{k \to \infty} \frac{\left( \mathcal{E}^h(t,y) - \mathcal{E}^h(t,y_k) - \lambda^h \mathcal{D}^h(y,y_k)\norm*{\grad[h]y - \grad[h]y_k}_{\SobW^{1,\infty}(\Omega)} \right)^+}{\mathcal{D}^h(y,y_k)(1 + c^h\norm*{\grad[h]y_k - \grad[h]y}_{\Leb^\infty(\Omega)})^{1/2}}.
\end{align*}
Clearly, the right-hand side of the second equation is bounded from above by the right-hand side of \eqref{Eq:Local-slope-3D:Characterization} since $y_k \in \mathcal{A}(h,M,y,\delta^h)$ for $k \gg 1$.
\medskip

\stepemph{Step 2 -- Lower bound for \eqref{Eq:Local-slope-3D:Characterization}:} 
For the reverse estimate consider $\hat y \in \mathcal A(h,M,y,\delta^h) \subset \mathcal Q^{h,M-}_{\rm bc}$. Let $y_s := y + s(\hat y - y)$. Then, we obtain from \eqref{Eq:Lambda-convexity-energy:Dissipation} that $\mathcal{D}^h(y,y_s) \leq c_1 s \to 0$ as $s \to 0$, so that the definition of the local slope and the semi-convexity established in \Cref{Lem:Properties-MGF-bulk}~(\ref{Item:Properties-MFG-bulk:convexity}) show
\begin{align*}
	\abs*{\partial \mathcal{E}^h(t)}_{\mathcal{D}^h}(y)
	&\geq \limsup_{s \to 0} \frac{\left( \mathcal{E}^h(t,y) - \mathcal{E}^h(t,y_s) \right)^+}{\mathcal{D}^h(y,y_s)} \\
	&\geq \limsup_{s \to 0} \frac{s\left( \mathcal{E}^h(t,y) - \mathcal{E}^h(t,\hat y) - \lambda^h(1-s) \mathcal{D}^h(y,\hat y)\norm*{\grad[h]\hat y - \grad[h]y}_{\SobW^{1,\infty}(\Omega)} \right)^+}{s\mathcal{D}^h(y,\hat y)(1 + c^h\norm*{\grad[h]\hat y - \grad[h]y}_{\Leb^\infty(\Omega)})^{1/2}} \\
	&= \frac{\left( \mathcal{E}^h(t,y) - \mathcal{E}^h(t,\hat y) - \lambda^h \mathcal{D}^h(y,\hat y)\norm*{\grad[h]\hat y - \grad[h]y}_{\SobW^{1,\infty}(\Omega)} \right)^+}{\mathcal{D}^h(y,\hat y)(1 + c^h\norm*{\grad[h]\hat y - \grad[h]y}_{\Leb^\infty(\Omega)})^{1/2}}.
\end{align*}
Since the conclusion holds for all $\hat y \in \mathcal A(h,M,y,\delta^h)$, we can pass to the supremum and conclude \eqref{Eq:Local-slope-3D:Characterization}.
\end{proof}

We use an adaptation of \cite[Cor.~2.4.10]{AGS05} to derive the chain rule and the lower semi-continuity of the local slope from the representation formula in \cref{Cor:Local-slope-3D:Characterization}. In fact, the lower semi-continuity follows almost immediately as point-wise suprema of lower semi-continuous functions are again lower semi-continuous. This procedure has already been successfully used in settings similar to ours, see e.g.\ \cite[Lem.~4-9]{FK18} and \cite[Lem.~5.4]{FM23}. For a not explicitly time-dependent energy, the chain rule is equivalent to the local slope being a strong upper gradient in the terminology of \cite{AGS05}. In our setting, we use an argument in the spirit of \cite[Prop.~A.1]{MR23} to obtain the chain rule from the strong upper gradient property.

\begin{corollary}
\label{Cor:Local-slope-3d-properties}
Let $M > 0$. There exists $h_0 > 0$ such that for all $0 < h < h_0$ the following statements hold:
\begin{enumerate}[(a)]
	\item \label{Item:Local-slope-3d-properties:lsc}
	For all $t \in [0,T]$, $y \in \mathcal{Q}^h_{\rm bc}$ and sequences $(y_k) \subset \mathcal{Q}^{h,M}_{\rm bc}$, we have
	\begin{equation*}
		y_k \wto y \text{ in } \SobW^{3,p_*}(\Omega, \R^3)
		\quad\Rightarrow\quad
		\liminf_{k \to \infty} \abs[\big]{\partial \mathcal{E}^h(t)}_{\mathcal{D}^h}(y_k)
		\geq \abs[\big]{\partial \mathcal{E}^h(t)}_{\mathcal{D}^h}(y).
	\end{equation*}
	
	\item \label{Item:Local-slope-3d-properties:chain-rule}
	If $y:[0,T] \to \mathcal{Q}^{h,M}_{\rm bc}$ is absolutely continuous w.r.t.\ $\mathcal{D}^h$ and satisfies
	\begin{equation*}
		\int_0^T \abs[\big]{\partial \mathcal{E}^h(t)}_{\mathcal{D}^h}\of{y(t)}\abs[\big]{y'}_{\mathcal D^h}(t) \,\dd t < \infty,
	\end{equation*}
	then, $t \mapsto \mathcal E^h\of{t,y(t)}$ is absolutely continuous and for a.e.\ $t \in [0,T]$,
	\begin{equation}
	\label{Eq:chain-rule-inequality}
		\tfrac{\dd}{\dd t} \mathcal E^h\of{t,y(t)}
		\geq \partial_t \mathcal E^h\of{t,y(t)} - \abs[\big]{\partial \mathcal{E}^h(t)}_{\mathcal{D}^h}\of{y(t)} \abs[\big]{y'}_{\mathcal D^h}(t).
	\end{equation}
\end{enumerate}
\end{corollary}

\begin{proof}
Consider the situation of \Cref{Cor:Local-slope-3D:Characterization}.
\medskip

\stepemph{Step 1 -- Proof of (\ref{Item:Local-slope-3d-properties:lsc}):}
From \Cref{Lem:Properties-MGF-bulk}~(\ref{Item:Properties-MFG-bulk:LSC-energy}) and (\ref{Item:Properties-MFG-bulk:cont-distance}) we obtain that for all $\hat y \in \mathcal{Q}^{h,M-}_{\rm bc}$ the map
\begin{equation*}
	\mathcal Q^{h,M}_{\rm bc}\ni y \mapsto \frac{\left( \mathcal{E}^h(t,y) - \mathcal{E}^h(t,\hat y) - \lambda^h\mathcal{D}^h(y,\hat y)\norm*{\grad[h]y - \grad[h]\hat y}_{\SobW^{1,\infty}(\Omega)}\right)^+}{\mathcal{D}^h(y,\hat y)(1 + c^h\norm*{\grad[h]\hat y - \grad[h]y}_{\Leb^\infty(\Omega)})^{1/2}}
\end{equation*}
is lower semi-continuous. Thus, (\ref{Item:Local-slope-3d-properties:lsc}) follows from the characterization in \Cref{Cor:Local-slope-3D:Characterization}, the fact that point-wise suprema of lower semi-continuous functions are again lower semi-continuous and the fact that for all sequences $(y_k) \subset \mathcal{Q}^{h,M}_{\rm bc}$ with $y_k \wto y$ in $\SobW^{3,p_*}(\Omega, \R^3)$, we have $\mathcal{A}(h,M,y,\delta^h) \subset \bigcup_{k_0 \in \N} \bigcap_{k \geq k_0} \mathcal{A}(h,M,y_k,\delta^h)$.%
\medskip

\stepemph{Step 2 -- The local slope is a strong upper gradient:}
Let throughout this step $t \in [0,T]$ be arbitrary, but fixed. We claim that for all $z \in \AC([0,T]; \mathcal Q^{h,M}_{\rm bc})$ with
\begin{equation}\label{Pf:Lem:MetricDh1:SUP1}
	\int_0^T \abs*{z'}_{\mathcal D^h}(s) \abs*{\partial \mathcal E^h(t)}_{\mathcal D^h}\of[\big]{z(s)} \,\dd s < \infty,
\end{equation}
the map $s \mapsto \mathcal E^h(t,z(s))$ is absolutely continuous and for a.e.~$s\in[0,T]$ (more precisely for all $s$ for which $\tfrac{\dd}{\dd s} \mathcal E^h(t,z(s))$ and $\abs*{z'}_{\mathcal D^h}(s)$ exist), we have
\begin{equation}\label{Pf:Lem:MetricDh1:SUP2}
	\abs*{\tfrac{\dd}{\dd s} \mathcal E^h(t,z(s))}
	\leq \abs*{z'}_{\mathcal D^h}(s) \abs*{\partial \mathcal E^h(t)}_{\mathcal D^h}\of[\big]{z(s)}.
\end{equation}
The inequality \eqref{Pf:Lem:MetricDh1:SUP2} is a straight-forward consequence of the definition of the local slope and metric derivative in \Cref{Def:Intro:Metric-gradient-flow}, see \cite[Thm.~1.2.5]{AGS05}. Thus, we only show here that $s \mapsto \mathcal E^h(t,z(s))$ is absolutely continuous. 
For the argument we note that by \Cref{Lem:Properties-MGF-bulk}~(\ref{Item:Properties-MFG-bulk:cont-distance}), there exists a partition $0 =: s_0 < \dots < s_N =: T$, such that
\begin{equation*}
	\norm{\grad[h]z(b) - \grad[h]z(a)}_{\SobW^{1,\infty}(\Omega)} < \delta^h, \qquad
	\text{for all } a,b \in [s_{i-1}, s_i], i = 1,\dots, N.
\end{equation*}
Evidently it suffices to show that for any segment $i=1,\ldots, N$, the function $[s_{i-1},s_i]\ni s\mapsto\mathcal E^h(t,z(s))$ is absolutely continuous. To that end set $\mathcal Q_i := z([s_{i-1},s_i])$ and note that it is a closed subset of $(\mathcal{Q}^{h,M}_{\rm bc}, \mathcal D^h)$. Thus, $(\mathcal Q_i,\mathcal D^h)$ is a complete metric space. 
Consider the global slope of the restricted energy $\mathcal E^h\big\vert_{\mathcal Q_i}$, i.e., the functional
\begin{equation*}
	\mathcal{G}_i:\mathcal Q_i\to[0,\infty],\qquad \mathcal{G}_i(\tilde z) 
	:= \sup_{\tilde z \neq \hat z \in \mathcal Q_i}  \frac{\left( \mathcal{E}^h(t,\tilde z) - \mathcal{E}^h(t,\hat z) \right)^+}{\mathcal{D}^h(\tilde z, \hat z)}.
\end{equation*}
We observe from \eqref{Eq:Local-slope-3D:Characterization} that for some constant $c_1 > 0$ we have the bound
\begin{equation*}
	\mathcal{G}_i(\tilde z) 
	\leq c_1\Bigl(\abs*{\partial \mathcal{E}^h(t)}_{\mathcal{D}^h}(\tilde z) + 1\Bigr)\qquad\text{for all }\tilde z\in\mathcal Q_i.
\end{equation*}
Thus, in view of condition \eqref{Pf:Lem:MetricDh1:SUP1}, we conclude $(\mathcal G_i\circ z)\abs*{z'}_{\mathcal D^h}\in \Leb^1([s_{i-1},s_i])$. Since the energy is lower semi-continuous by \Cref{Lem:Properties-MGF-bulk}~(\ref{Item:Properties-MFG-bulk:LSC-energy}), we conclude with \cite[Thm.~1.2.5]{AGS05} that $\mathcal G_i$ is a strong upper gradient for $\mathcal E^h\big\vert_{\mathcal Q_i}$ and thus, $s \mapsto \mathcal{E}^h(t,z(s))$ is absolutely continuous on $[s_{i-1},s_i]$.
\medskip

\stepemph{Step 3 -- Estimate for the local slopes:}
We show that there exists a constant $c_1 >0$ such that for all $s,t \in [0,T]$ and $y \in \mathcal Q^{h,M}_{\rm bc}$, we have the estimate,
\begin{equation*}
	\abs*{\partial \mathcal{E}^h(s)}_{\mathcal{D}^h}(y)
	\leq \abs*{\partial \mathcal{E}^h(t)}_{\mathcal{D}^h}(y) + c_1.
\end{equation*}
From the definition of the local slope, we obtain,
\begin{align*}
	\abs*{\partial \mathcal{E}^h(s)}_{\mathcal{D}^h}(y)
	&= \limsup_{\mathcal D^h(y,\hat y)\to 0} \frac{\left( \mathcal{E}^h(s,y) - \mathcal{E}^h(s,\hat y) \right)^+}{\mathcal{D}^h(y, \hat y)}
	\\
	&= \limsup_{\mathcal D^h(y,\hat y)\to 0} \frac{\left( \mathcal{E}^h(t,y) - \mathcal{E}^h(t,\hat y) + \int_t^s \partial_t \mathcal E^h(r,y) - \partial_t \mathcal E^h(r,\hat y) \,\dd r \right)^+}{\mathcal{D}^h(y, \hat y)}
	\\
	&= \limsup_{\mathcal D^h(y,\hat y)\to 0} \frac{\left( \mathcal{E}^h(t,y) - \mathcal{E}^h(t,\hat y) - \int_t^s \sprod{\dot{l}(r)}{\hat{y} - y} \,\dd r \right)^+}{\mathcal{D}^h(y, \hat y)}
	\\
	&\leq \abs*{\partial \mathcal{E}^h(t)}_{\mathcal{D}^h}(y) + \left(\int_0^T \norm*{\partial_t l(r)}_{\SobH^1(\Omega)'} \,\dd r\right) \limsup_{\mathcal D^h(y,\hat y)\to 0} \frac{\norm*{\hat y - y}_{\SobH^1(\Omega)}}{\mathcal{D}^h(y, \hat y)}.
\end{align*}
Note that we can take without loss of generality $\hat{y} \in \mathcal{Q}^{h,M-}_{\rm bc}$. The second term on the right-hand side is bounded in view of $\partial_t l \in \Leb^1([0,T]; \SobH^1(\Omega,\R^3)')$, \Cref{Lem:Properties-MGF-bulk}~(\ref{Item:Properties-MFG-bulk:equivalence-distance}) and the Poincar{\'e} inequality.
\medskip

\stepemph{Step 4 -- Proof of the chain rule (\ref{Item:Local-slope-3d-properties:chain-rule}):}
Let $y \in \AC([0,T]; \mathcal Q^{h,M}_{\rm bc})$ with 
\begin{equation*}
	\int_0^T \abs*{y'}_{\mathcal D^h}(t) \abs*{\partial \mathcal E^h(t)}_{\mathcal D^h}\of[\big]{y(t)} \,\dd t < \infty.
\end{equation*}
For the sake of clarity we set 
\begin{equation*}
\mathcal I_t(s):=\mathcal E^h(t,y(s)),\qquad\mathcal I(t):=\mathcal I_t(t).
\end{equation*}
We first prove that $t \mapsto \mathcal E^h(t,y(t))=\mathcal I(t)$ is absolutely continuous. 
Indeed, by Step 3 we have
\begin{equation*}
	\forall s,t\in[0,T]\,:\qquad \int_0^T \abs*{y'}_{\mathcal D^h}(s) \abs*{\partial \mathcal E^h(t)}_{\mathcal D^h}\of[\big]{y(s)} \,\dd s < \infty.
\end{equation*}
Thus, Step 2 shows that the map $s \mapsto \mathcal E^h(t,y(s))=\mathcal I_t(s)$ is absolutely continuous and
\begin{equation}\label{Pf:MetrichDh1:33}
	\abs*{\mathcal I_t'(s)}\leq \abs*{y'}_{\mathcal D^h}(s) \abs*{\partial \mathcal E^h(t)}_{\mathcal D^h}\of[\big]{y(s)}\qquad
	\text{for a.e.\ } s \in [0,T].
\end{equation}
We conclude that
\begin{align*}
	\abs*{\mathcal I(t)-\mathcal I(s)}
	&\leq \abs[\big]{\mathcal I_t(t)-\mathcal I_t(s)} + \abs[\big]{\mathcal I_t(s)-\mathcal I_s(s)}
	\leq \int_s^t \abs*{\mathcal I_t'(r)} + \abs*{\partial_r \mathcal E^h\of[\big]{r,y(s)}} \,\dd r.
\end{align*}
With \eqref{Pf:MetrichDh1:33}, Step 3 and \cref{Lem:Regularity-energy} we conclude that the integrand on the right-hand side can be bounded by a function in $\Leb^1([0,T])$ that is independent of $s$ and $t$ and thus $\mathcal I$ is absolutely continuous.
It remains to prove the chain-rule inequality. 
Let $\hat{\Lambda}\subset[0,T]$ denote a set of full-measure such that the derivatives $\mathcal I'(t)$, $\abs*{y'}_{\mathcal D^h}(t)$, and $\dot l(t)$ exist for all $t\in\hat{\Lambda}$. Note that the existence of $\hat{\Lambda}$ follows from the absolute continuity of $\mathcal I$, $y$, and $l$.  
Then for all $t\in\hat{\Lambda}$ we have	
\begin{align*}
	\frac{\mathcal I_t(t)-\mathcal I_t(s)}{t-s}
	&=\frac{\mathcal I(t)-\mathcal I(s)}{t-s}-\frac{\mathcal I_t(s)-\mathcal I_s(s)}{t-s}
	\\
	&=\frac{\mathcal I(t)-\mathcal I(s)}{t-s}-\sprod*{\tfrac{l(s) - l(t)}{t-s}}{y(s)}
	\to \mathcal I'(t)-\partial_t\mathcal E^h(t,y(t))\qquad\text{as }s\to t.
\end{align*}
Hence, the limit of the left-hand side, and thus $\mathcal I_t'(t)=\frac{\dd}{\dd s}\mathcal E^h(t,y(s))\vert_{s=t}$ exists and equals the right-hand side.
Since $\abs*{y'}_{\mathcal D^h}(t)$ exists as well, we may apply \eqref{Pf:Lem:MetricDh1:SUP2} with $s=t$ and we conclude
\begin{equation*}
	\mathcal I_t'(t)\geq -\abs*{y'}_{\mathcal D^h}(t) \abs*{\partial \mathcal E^h(t)}_{\mathcal D^h}\of[\big]{y(t)}\qquad\text{for all }t\in\hat{\Lambda}.
\end{equation*}
Combined with $\mathcal I'_t(t)=\mathcal I'(t)-\partial_t\mathcal E^h(t,y(t))$, the chain rule \eqref{Eq:chain-rule-inequality} follows for all $t \in \hat{\Lambda}$.
\end{proof}

\subsection{Characterization of the metric gradient flow}
\label{Sec:Characterization-3D}

In this section we characterize curves of maximal slope for $(\mathcal Q^h_{\rm bc}, \mathcal D^h, \mathcal E^h)$ rigorously by a partial differential equation in the form of \eqref{Eq:Intro:PDE}. Therefore, fix $y \in \mathcal Q^{h,<M}_{\rm bc}$ for some $M > 0$ and $0 < h \ll 1$ and consider the dissipation potential,
\begin{equation}
	\mathcal R^h(y,\dot y) := \lim_{\eps \to 0} \tfrac{\eps^{-2}}{2} \mathcal D^h(y + \eps \dot{y}, y)^2 = 
	\tfrac{h^{-2}}{4} \int_\Omega \dif_{F_1}^2 \mathscr D_\vi^2(\grad[h]y, \grad[h]y)[\grad[h]\dot y, \grad[h]\dot y], \qquad
	\dot y \in \SobH^1(\Omega, \R^3).
\end{equation}
We note that $\mathcal R^h(y,\placeholder)$ is quadratic and thus especially Fr{\'e}chet differentiable w.r.t.\ $\SobH^1(\Omega, \R^3)$ and denote the derivative by $\dif_{\dot y} \mathcal R^h(y,\placeholder) \in \SobH^1(\Omega, \R^3)'$.
As we consider boundary conditions it suffices to restrict to the space
\begin{equation}
	\SobH^1_{S,0}(\Omega, \R^3) := \set*{\dot y \in \SobH^1(\Omega, \R^3) \given \dot y(0,\placeholder) = 0 \text{ a.e.\ in } S}.
\end{equation}
Analogously, we define $\SobW^{3,p_*}_{S,0}(\Omega, \R^3)$.
For the energy, we notice that $\mathcal E^h(t,y + \placeholder)$ is Fr{\'e}chet differentiable w.r.t.\ $\SobW^{3,p_*}(\Omega, \R^3)$ at $0$ and denote the derivative by $\dif^E \mathcal E^h(t,y)$. 
However, as we consider a gradient flow w.r.t.\ the metric $\mathcal D^h$, the differential $\dif^E \mathcal E^h$ is of no direct use. Instead, since $\mathcal D^h$ is equivalent to the distance induced by the $\SobH^1(\Omega, \R^3)$-norm, a natural choice for the derivative $\dif \mathcal E^h(t,y)$ in \eqref{Eq:Intro:PDE} is to consider the Fr{\'e}chet subdifferential
\begin{equation}
	\dif \mathcal E^h(t,y)
	:= \set*{\xi \in \SobH^1_{S,0}(\Omega, \R^3)' \given \liminf_{\norm{\dot y}_{\SobH^1(\Omega)} \to 0}\tfrac{\mathcal E^h(t,y + \dot y) - \mathcal E^h(t,y) - \sprod{\xi}{\dot y}}{\norm{\dot y}_{\SobH^1(\Omega)}} \geq 0}.
\end{equation}
We obtain the following characterization.

\begin{proposition} \label{Prop:Characterization-solutions-3D}
For all $M > 0$, there exists $h_0 > 0$, such that for all $0 < h \leq h_0$ the following statements hold. Let $y \in \AC([0,T]; \mathcal Q^{h,<M}_{\rm bc})$. Then, $y$ is a curve of maximal slope for $(\mathcal Q^h_{\rm bc}, \mathcal D^h, \mathcal E^h)$, if and only if $y \in \SobH^1(0,T; \SobH^1(\Omega, \R^3))$, $t \mapsto \mathcal E^h(t,y(t))$ is absolutely continuous and $y$ satisfies for a.e.\ $t \in [0,T]$ the differential inclusion,
\begin{equation} \label{Eq:Characterization-3D:PDE}
	\dif \mathcal E^h\of[\big]{t,y(t)} + \dif_{\dot y} \mathcal R^h\of[\big]{y(t), \dot y(t)} \ni 0 \qquad
	\text{in } \SobH^1_{S,0}(\Omega, \R^3)'.
\end{equation}
Moreover, \eqref{Eq:Characterization-3D:PDE} is equivalent to
\begin{equation} \label{Eq:Characterization-3D:PDE2}
	\dif^E \mathcal E^h\of[\big]{t,y(t)} + \dif_{\dot y} \mathcal R^h\of[\big]{y(t), \dot y(t)} = 0 \qquad
	\text{in } \SobW^{3,p_*}_{S,0}(\Omega, \R^3)',
\end{equation}
that is, for all $\nu \in \SobW^{3,p_*}_{S,0}(\Omega, \R^3)$ it holds,
\begin{equation}
\begin{aligned}
	\sprod*{l(t)}{\nu}
	&= h^{-2} \int_\Omega \dif W_\el\of[\big]{\grad[h]y(t)} \big[\grad[h]\nu \big]
	\\
	&\qquad+ h^{-\alpha_C p_*} \int_\Omega \dif H_\el^C\of[\big]{\grad C_{y(t)}}\big[2\grad \sym(\grad[h]\nu^\transp \grad[h]y(t)) \big]
	\\
	&\qquad+ h^{\alpha_F p_*} \int_\Omega \dif H_\el^F\of[\big]{\grad\grad[h]y(t)}\big[\grad\grad[h]\nu \big]
	\\
	&\qquad+ h^{\alpha_{\mathbb F} p_*} \int_\Omega \dif H_\el^{\mathbb F}\of[\big]{\grad^2\grad[h]y(t)}\big[\grad^2\grad[h]\nu \big]
	\\
	&\qquad+ \tfrac{h^{-2}}{2} \int_\Omega \dif_{F_1}^2 \mathscr D_\vi^2\of[\big]{\grad[h]y(t), \grad[h]y(t)}\big[\grad[h]\dot y(t), \grad[h]\nu \big].
\end{aligned}
\end{equation}
\end{proposition}

For the proof of \cref{Prop:Characterization-solutions-3D} we use the fact that $\mathcal R^h(y,\placeholder)$ introduces a scalar product on $\SobH^1(\Omega, \R^3)$ by
\begin{equation*}
	\sprod*{y_1}{y_2}_y :=
	\tfrac{h^{-2}}{2} \int_\Omega \dif_{F_1}^2 \mathscr D_\vi^2(\grad[h]y, \grad[h]y)[\grad[h]y_1, \grad[h]y_2].
\end{equation*}
In fact, a close look at the proof of \Cref{Lem:Properties-MGF-bulk}~(\ref{Item:Properties-MFG-bulk:equivalence-distance}) reveals that the induced norm
\begin{equation}\label{def:diss:norm}
	\norm{\dot y}_y := \sqrt{\sprod{\dot y}{\dot y}} = \sqrt{2\mathcal R^h(y,\dot y)}	
\end{equation}
is equivalent to the $\SobH^1(\Omega, \R^3)$-norm and therefore to the distance $\mathcal{D}^h$:

\begin{lemma}
For all $M > 0$, there exists $h_0 > 0$, such that for all $0 < h \leq h_0$ there exists a constant $c > 0$, such that for all $y \in \mathcal Q^{h,M-}_{\rm bc}$ and $\dot y \in \SobH^1(\Omega, \R^3)$,
\begin{equation}
	\tfrac1c \norm*{\dot y}_{\SobH^1(\Omega)} \leq \norm*{\dot y}_y \leq c \norm*{\dot y}_{\SobH^1(\Omega)}.
\end{equation}
\end{lemma}

In fact, we can sharpen the estimate. Following the \eqref{Eq:Lambda-convexity-energy:Dissipation}, we see that \eqref{Eq:TaylorDissipation2} implies that $\norm{\placeholder}_y$ is even equal to the distance $\mathcal D^h$ up to a small error:

\begin{lemma} \label{Lem:Characterization-3D:Equivalence-Hilbert-space}
For all $M > 0$, there exists $h_0 > 0$, such that for all $0 < h \leq h_0$ there exists a constant $c > 0$, such that for all $y, \hat y \in \mathcal Q^{h,M-}_{\rm bc}$,
\begin{subequations}
\begin{align}
	\mathcal D^h(y,\hat y) 
	&\leq (1 + c\norm*{\grad[h]\hat y - \grad[h]y}_{\Leb^\infty(\Omega)})^{1/2}\, \norm*{\hat y - y}_y,
	\\
	\norm*{\hat y - y}_y
	&\leq (1 + c\norm*{\grad[h]\hat y - \grad[h]y}_{\Leb^\infty(\Omega)})^{1/2}\, \mathcal D^h(y,\hat y).
\end{align}
\end{subequations}
\end{lemma}

\Cref{Lem:Characterization-3D:Equivalence-Hilbert-space} implies that the local slope and metric derivative in w.r.t.\ $\mathcal D^h$ and $\norm{\placeholder}_y$ coincide in $\mathcal{Q}^{h,M-}_{\rm bc}$, see below for the precise statement. As the latter is (the subset of) a Hilbert space, we can characterize curves of maximal slope similarly as in \cref{Sec:Rod-models}. For the characterization of the local slope we use a correspondence of the Fr{\'e}chet subdifferential $\dif \mathcal E^h(t,y)$ and the derivative $\dif^E \mathcal{E}^h(t,y)$ which is due to the semi-convexity \Cref{Lem:Properties-MGF-bulk}~(\ref{Item:Properties-MFG-bulk:convexity}):

\begin{lemma}
For all $M > 0$, there exists $h_0 > 0$, such that for all $0 < h \leq h_0$, $t \in [0,T]$ and $y \in \mathcal{Q}^{h,<M}_{\rm bc}$, we have $\xi \in \dif \mathcal E^h(t,y)$ if and only if
\begin{equation}
	\label{Eq:Characterization-3D:Subdifferential-energy}
	\xi \in \SobH^1_{S,0}(\Omega, \R^3)'
	\quad\text{and}\quad
	\forall\dot y \in \SobW^{3,p_*}_{S,0}(\Omega, \R^3):\:
	\sprod{\xi}{\dot y} = \dif^E\mathcal E^h(t,y)[\dot y].
\end{equation}
Especially, $\dif \mathcal E^h(t,y)$ is at most a singleton.
\end{lemma}

\begin{proof}
Let first $\xi \in \dif \mathcal E^h(t,y)$. From the definition of the Fr{\'e}chet subdifferential, we obtain for all $\dot y \in \SobW^{3,p_*}_{S,0}(\Omega, \R^3)$,
\begin{equation*}
	\frac{\sprod{\xi}{\pm \dot y}}{\norm{\dot y}_{\SobH^1(\Omega)}}
	\leq \liminf_{s \downarrow 0} \frac{\mathcal{E}^h(t,y \pm s\dot y) - \mathcal{E}^h(t,y)}{s\norm{\dot y}_{\SobH^1(\Omega)}}
	= \frac{\dif^E \mathcal{E}^h(t,y)[\pm \dot y]}{\norm{\dot y}_{\SobH^1(\Omega)}} 
\end{equation*}
which implies \eqref{Eq:Characterization-3D:Subdifferential-energy}.
On the other hand, suppose that $\xi \in \SobH^1_{S,0}(\Omega, \R^3)'$ satisfies \eqref{Eq:Characterization-3D:Subdifferential-energy}. Let $(\dot{y}_k) \subset \SobH^1_{S,0}(\Omega)$ be such that $0 \neq \dot{y}_k \to 0$ in $\SobH^1(\Omega, \R^3)$ and
\begin{equation*}
	\liminf_{\norm{\dot y}_{\SobH^1(\Omega)} \to 0}\frac{\mathcal E^h(t,y + \dot y) - \mathcal E^h(t,y) - \sprod{\xi}{\dot y}}{\norm{\dot y}_{\SobH^1(\Omega)}}
	= \lim_{k \to \infty}\frac{\mathcal E^h(t,y + \dot y_k) - \mathcal E^h(t,y) - \sprod{\xi}{\dot y_k}}{\norm{\dot y_k}_{\SobH^1(\Omega)}}.
\end{equation*}
Suppose without loss of generality that this limit is finite. Then, $\limsup_{k \to \infty} \mathcal{E}^h(t,y + \dot{y}_k) \leq \mathcal{E}^h(t,y) < M$, so that $y + \dot{y}_k \in \mathcal{Q}^{h,M-}_{\rm bc}$. Then, the definition of $\dif^E\mathcal{E}^h(t,y)$ and the semi-convexity \Cref{Lem:Properties-MGF-bulk}~(\ref{Item:Properties-MFG-bulk:convexity}) show that for all $k \gg 1$,
\begin{align*}
	0 &= \lim_{s \downarrow 0} \frac{\mathcal{E}^h(t,y + s\dot{y}_k) - \mathcal{E}^h(t,y) - \dif^E\mathcal{E}^h(t,y)[s\dot{y}_k]}{s}
	= \lim_{s \downarrow 0} \frac{\mathcal{E}^h(t,y + s\dot{y}_k) - \mathcal{E}^h(t,y) - \sprod{\xi}{s\dot y_k}}{s}
	\\
	&\leq \mathcal{E}^h(t,y + \dot{y}_k) - \mathcal{E}^h(t,y) + \lambda^h \mathcal{D}^h(y,y+\dot{y}_k)\norm{\dot{y}_k}_{\SobW^{1,\infty}(\Omega)} - \sprod{\xi}{\dot{y}_k}
	\\
	&\leq \mathcal{E}^h(t,y + \dot{y}_k) - \mathcal{E}^h(t,y) + c_1\norm{\dot{y}_k}_{\SobH^1(\Omega)}\norm{\dot{y}_k}_{\SobW^{1,\infty}(\Omega)} - \sprod{\xi}{\dot{y}_k}.
\end{align*}
We conclude,
\begin{align*}
	\lim_{k \to \infty}\frac{\mathcal E^h(t,y + \dot y_k) - \mathcal E^h(t,y) - \sprod{\xi}{\dot y_k}}{\norm{\dot y_k}_{\SobH^1(\Omega)}}
	\geq \lim_{k \to \infty} -c_1\norm{\dot{y}_k}_{\SobW^{1,\infty}(\Omega)} = 0
\end{align*}
and thus, $\xi \in \dif\mathcal{E}^h(t,y)$.
\end{proof}

We proceed to characterize metric derivatives and the local slope using the preceding lemmas.

\begin{lemma} \label{Lem:Characterization-3D:Local-slope}
For all $M > 0$, there exists $h_0 > 0$, such that for all $0 < h \leq h_0$ the following statements hold:
\begin{enumerate}[(a)]
	\item\label{Item:Characterization-3D:metric-derivative}
	A map $y: [0,T] \to \mathcal Q^{h,M-}_{\rm bc}$ is absolutely continuous w.r.t.\ $\mathcal D^h$, if and only if $y$ is absolutely continuous w.r.t.\ the $\SobH^1(\Omega, \R^3)$-norm and then,
	\begin{equation}
		\abs*{y'}_{\mathcal D^h}(t)
		= \norm{\dot y(t)}_{y(t)}
		\qquad 
		\text{for a.e.\ } t \in [0,T].
	\end{equation}
	
	\item\label{Item:Characterization-3D:local-slope}
	For all $y \in \mathcal Q^{h,<M}_{\rm bc}$ and $t \in [0,T]$, the local slope is finite if and only if $\dif \mathcal E^h(t,y)$ is not empty and then, the unique element $\xi \in \dif \mathcal E^h(t,y)$ satisfies
	\begin{equation*}
		\abs*{\partial \mathcal E^h(t)}_{\mathcal D^h}(y)
		= \norm{\xi}_{y,*}.
	\end{equation*}
	where $\norm{\placeholder}_{y,*}$ denotes the dual norm of $\norm{\placeholder}_{y}$.
\end{enumerate}
\end{lemma}

\begin{proof}
The equivalence of the notions of absolute continuity of a map $y: [0,T] \to \mathcal Q^{h,M}_{\rm bc}$ is due to the equivalence of $\mathcal D^h$ and the $\SobH^1(\Omega, \R^3)$-norm obtained in \cref{Lem:Properties-MGF-bulk}~(\ref{Item:Properties-MFG-bulk:equivalence-distance}). Then, as discussed in \cref{Def:Intro:Metric-gradient-flow}, $t\mapsto y(t)\in\SobH^1(\Omega, \R^3)$ is differentiable for a.e.~$t$. Hence, from \cref{Lem:Characterization-3D:Equivalence-Hilbert-space} we obtain for a.e.\ $t \in [0,T]$,
\begin{equation*}
	\abs*{y'}_{\mathcal D^\limSymbol}(t)
	= \lim_{s \to t} \frac{\mathcal D^h\of[\big]{y(s),y(t)}}{\abs*{t-s}}
	= \lim_{s \to t} \frac{\norm*{y(t) - y(s)}_{y(t)}}{\abs*{t-s}}
	= \norm*{\dot y(t)}_{y(t)}.
\end{equation*}
Note that here we use that $\lim_{s \to t} \norm{\grad[h]y(t) - \grad[h]y(s)}_{\Leb^\infty(\Omega)} = 0$ since \cref{Lem:Properties-MGF-bulk}~(\ref{Item:Properties-MFG-bulk:cont-distance}) implies $y(s) \wto y(t)$ in $\SobW^{3,p_*}(\Omega, \R^3)$. This completes the proof of (\ref{Item:Characterization-3D:metric-derivative}).
To prove (\ref{Item:Characterization-3D:local-slope}) let $t \in [0,T]$ and $y \in \mathcal Q^{h,<M}_{\rm bc}$. It is straight-forward to obtain from \cref{Lem:Characterization-3D:Equivalence-Hilbert-space} that $\abs*{\partial \mathcal E^h(t)}_{\mathcal D^h}(y) = \abs*{\partial \mathcal E^h(t)}_{\norm{\placeholder}_y}(y)$. 
Let $\xi \in \dif \mathcal E^h(t,y)$. Then, for all $\dot y \in \SobH^1_{S,0}(\Omega, \R^3)$ we find,
\begin{equation*}
	\abs*{\partial \mathcal E^h(t)}_{\norm{\placeholder}_y}(y)
	= \limsup_{\norm{\dot y}_y \to 0} \frac{\left(\mathcal E^h(t,y) - \mathcal E^h(t, y + \dot y)\right)^+}{\norm{\dot y}_y}
	\leq \limsup_{\norm{\dot y}_y \to 0} \,\sprod*{\xi}{\tfrac{\dot y}{\norm{\dot y}_y}}
	\leq \norm{\xi}_{y,*}.
\end{equation*}
On the other hand, suppose $\abs*{\partial \mathcal E^h(t)}_{\norm{\placeholder}_y}(y) < \infty$. For all $\dot y \in \SobW^{3,p_*}(\Omega, \R^3)$, we have
\begin{align*}
	\frac{\pm \dif^E \mathcal E^h(t,y)[\dot y]}{\norm{\dot y}_y}
	&= \lim_{s \to 0} \frac{\mathcal E^h(t,y) - \mathcal E^h(t,y \mp s\dot y)}{\norm{s\dot y}_y}
	\\
	&\leq \limsup_{\norm{\dot y}_y \to 0} \frac{\left( \mathcal E^h(t,y) - \mathcal E^h(t,y + \dot y) \right)^+}{\norm{\dot y}_y} 
	= \abs*{\partial \mathcal E^h(t)}_{\norm{\placeholder}_y}(y).
\end{align*}
Hence, $\dif^E \mathcal E^h(t,y)$ is bounded w.r.t.\ the $\SobH^1(\Omega, \R^3)$-norm and can be extended uniquely to some $\xi \in \SobH^1_{S,0}(\Omega, \R^3)'$. Then, \eqref{Eq:Characterization-3D:Subdifferential-energy} yields $\xi \in \dif \mathcal E^h(t,y)$ and the inequality above shows,
\begin{equation*}
	\norm{\xi}_{y,*}
	\leq \abs*{\partial \mathcal E^h(t)}_{\norm{\placeholder}_y}(y).
	\qedhere.
\end{equation*}
\end{proof}

Now we argue similarly to \cref{Prop:MetricGF-HilbertGF-correspondence} to prove \cref{Prop:Characterization-solutions-3D}.

\begin{proof}[Proof of \cref{Prop:Characterization-solutions-3D}]
Let $y: [0,T] \to \mathcal Q^{h,<M}_{\rm bc}$ be absolutely continuous w.r.t.\ $\mathcal D^h$ (or equivalently w.r.t.\ the $\SobH^1(\Omega, \R^3)$-norm). 
Since \eqref{Eq:Characterization-3D:PDE} and \eqref{Eq:Characterization-3D:PDE2} are equivalent as a consequence of \eqref{Eq:Characterization-3D:Subdifferential-energy}, it suffices to consider \eqref{Eq:Characterization-3D:PDE}.
\medskip

\stepemph{Step 1 -- Equivalence at points of differentiability:}
We first prove the equivalence under the following regularity assumption: Let $\hat{\Lambda}$ denote the set of all $t\in[0,T]$ such that the derivatives $\frac{\dd}{\dd t}\mathcal E^h(t,y(t))$, $\partial_t\mathcal E^h(t,y(t))$, $\abs*{y'}_{\mathcal D^h}(t)$, and $\dot y(t)$ exist and $\dif\mathcal E^h(t,y(t))\neq\emptyset$. We claim that
\begin{equation}\label{Eq:Characterization-3D:equiv}
	\forall t\in\hat{\Lambda}\,:\quad
	\eqref{Eq:Characterization-3D:PDE}%
	\ \Leftrightarrow\ 
	\tfrac{\dd}{\dd t} \mathcal E^h\of[\big]{t,y(t)} - \partial_t \mathcal E^h\of[\big]{t,y(t)}\leq -\tfrac12\abs*{y'}_{\mathcal D^h}^2(t) - \tfrac12 \abs*{\partial \mathcal E^h(t)}_{\mathcal D^h}^2\of[\big]{y(t)}.
\end{equation}
For the argument let $\xi(t)\in \dif\mathcal E^h(t,y(t))$. 
Then, the definition of the Fr{\'e}chet subdifferential yields $\tfrac{\dd}{\dd t} \mathcal E^h\of[\big]{t,y(t)} = \partial_t \mathcal E^h\of[\big]{t,y(t)} + \sprod*{\xi(t)}{\dot y(t)}$.
With help of Young's inequality and \cref{Lem:Characterization-3D:Local-slope} we get
\begin{equation} \label{Eq:Characterization-3D:slope-inequality}
\begin{aligned}
	&\tfrac{\dd}{\dd t} \mathcal E^h\of[\big]{t,y(t)} - \partial_t \mathcal E^h\of[\big]{t,y(t)}
	= \sprod{\xi(t)}{\dot y(t)}
	\\
	&\qquad\geq -\tfrac12 \norm*{\dot y(t)}_{y(t)}^2 - \tfrac12 \norm*{\xi(t)}_{y(t),*}^2
	= -\tfrac12\abs*{y'}_{\mathcal D^h}^2(t) - \tfrac12 \abs*{\partial \mathcal E^h(t)}_{\mathcal D^h}^2\of[\big]{y(t)}.
\end{aligned}
\end{equation}
Furthermore, using convex analysis we observe that $\mathcal R^h_*(y,\placeholder) := \frac12 \norm*{\placeholder}_{y,*}^2$ is the Legendre transform of $\mathcal R^h(y,\placeholder) = \frac12 \norm{\placeholder}_y^2$, see e.g.\ \cite[Sec.~2.3]{Mie14-PP}. Hence, by Fenchel duality we observe that equality in \eqref{Eq:Characterization-3D:slope-inequality} holds, if and only if 
\begin{equation}
	\label{Eq:Characterization-3D:inclusion}
	-\xi(t) \in \partial_{\dot y} \mathcal R^h\of[\big]{y(t),\dot y(t)},
\end{equation}
where $\partial_{\dot y} \mathcal R^h\of{y,\dot y}$ denotes the convex subdifferential.
Now, \eqref{Eq:Characterization-3D:equiv} follows: On the one hand, by Fr{\'e}chet differentiability of $\mathcal R^h(y,\placeholder)$,  \eqref{Eq:Characterization-3D:PDE} is equivalent to \eqref{Eq:Characterization-3D:inclusion} and thus to equality holding in \eqref{Eq:Characterization-3D:slope-inequality} and on the other hand, the latter is equivalent to the right-hand side of \eqref{Eq:Characterization-3D:equiv}.
\medskip

\stepemph{Step 2 -- Curves of maximal slope satisfy the equation:}
Now, assume that $y$ is a curve of maximal slope for $(\mathcal Q^h_{\rm bc}, \mathcal D^h, \mathcal E^h)$. 
Then, with help of \cref{Lem:Characterization-3D:Local-slope} it is easy to see that $\hat{\Lambda}$ (the set of regular points from the previous step) has full measure. 
Hence, by \eqref{Eq:Intro:Curve-of-max-slope} and \eqref{Eq:Characterization-3D:equiv} the differential equation \eqref{Eq:Characterization-3D:PDE} holds for a.e.\ $t\in[0,T]$ and then $y \in \SobH^1(0,T;\SobH^1(\Omega, \R^3))$ follows from equality holding in \eqref{Eq:Characterization-3D:slope-inequality} and the fact that $t \mapsto \abs{\tfrac{\dd}{\dd t} \mathcal E^h\of{t,y(t)}} + \abs{\partial_t \mathcal E^h\of{t,y(t)}} \in \Leb^1([0,T])$.%
\medskip

\stepemph{Step 3 -- Solutions of the equation are curves of maximal slope:}
On the other hand, assume that $y$ satisfies \eqref{Eq:Characterization-3D:PDE}, $y \in \SobH^1(0,T; \SobH^1(\Omega, \R^3))$ and $t \mapsto \mathcal E^h(t,y(t))$ is absolutely continuous. Then by \eqref{Eq:Characterization-3D:PDE} the Fr{\'e}chet subdifferential $\dif \mathcal E^h(t,y(t))$ is not empty and with help of \cref{Lem:Characterization-3D:Local-slope} we see that $\hat{\Lambda}$ has full-measure. Now, \eqref{Eq:Characterization-3D:equiv}  implies that $y$ is a curve of maximal slope for $(\mathcal Q^h_{\rm bc}, \mathcal D^h, \mathcal E^h)$.
\end{proof}

\section{Convergence of the metric gradient flows}
\label{Sec:Convergence-metric-gradient-flow}

This section is dedicated to the proof of \cref{Thm:Main-result:convergence-to-limit,Thm:Main-result:convergence-to-reduced}. 
We formulate a general convergence result that is in the spirit of \cite{Ser11,FK18,FM23}.
Let therefore $(Q^h, D^h)$ and $(Q^0, D^0)$ be metric spaces and consider measurable maps $E^h: [0,T] \times Q^h \to (-\infty,\infty]$, $E^0: [0,T] \times Q^0 \to (-\infty,\infty]$. Let $(\tilde Q, \sigma)$ be a topological Hausdorff space, and $\pi^h: Q^h \to \tilde Q$ and $\pi^0: Q^0 \to \tilde Q$ maps. Assume the following regularity conditions:

\begin{enumerate}[(i)]
	\item (Domain independent of time):
	For all $\bar h \geq 0$ there exists a non-empty set $\operatorname{dom} E^{\bar h} \subset Q^{\bar h}$, such that $\set{(t,q) \in [0,T] \times Q^{\bar h} \given E^{\bar h}(t,q) < \infty} = [0,T] \times \operatorname{dom} E^h$.
	
	\item (Absolute continuity w.r.t.\ time):
	For all $\bar{h} \geq 0$ and $q \in \operatorname{dom} E^{\bar h}$, the map $E^{\bar h}(\placeholder, q):[0,T] \to \R$ is absolutely continuous.

	\item (Uniform power estimate):
	There exists a set $\Lambda \subset [0,T]$ of full measure, a strictly positive map $\lambda \in \Leb^1([0,T])$ and $c > 0$, such that for all $\bar{h} \geq 0$, $t \in \Lambda$ and $q \in \operatorname{dom} E^{\bar h}$, $\partial_t E^{\bar h}(t,q)$ exists and
	\begin{equation*}
		\abs[\big]{\partial_t E^{\bar h}(t,q)} \leq \lambda(t) \bigl( E^{\bar h}(t,q) + c \bigr).
	\end{equation*}
	
	\item (Chain rule for the limit):
	For all $q \in \AC([0,T]; Q^0)$ with
	\begin{equation*}
		\sup_{t \in [0,T]} E^0(t,q(t)) < \infty
		\quad\text{and}\quad
		\int_0^T \abs*{q'}_D(t) \abs*{\partial E^0(t)}_{D^0}\of[\big]{q(t)} \,\dd t < \infty,
	\end{equation*}
	the map $t \mapsto E(t,q(t))$ is absolutely continuous and we have the chain rule inequality
	\begin{equation}\label{Eq:chainrule:limit}
		\tfrac{\dd}{\dd t} E^0\of[\big]{t,q(t)} \geq \partial_t E^0\of[\big]{t,q(t)} - \abs*{q'}_{D^0}(t)\abs*{\partial E^0(t)}_{D^0}\of[\big]{q(t)}
		\qquad \text{for a.e.\ } t \in [0,T].
	\end{equation}
\end{enumerate}

We note that in our setting where $(Q^h,D^h,E^h)=(\mathcal{Q}^h_{\rm bc}, \mathcal D^h,\mathcal E^h) $ and $(Q^0,D^0,E^0) = (\mathcal Q^0_{\rm bc}, \mathcal D^0, \mathcal E^0)$ (or respectively $(Q^0,D^0,E^0) = (\mathcal Q^{\rm red}_{\rm bc}, \mathcal D^{\rm red}, \mathcal E^{\rm red})$), it is easy to verify these properties: 
The uniform power estimate is a consequence of \Cref{Lem:Regularity-energy} (the estimate for the limits follows analogously) and the chain rule follows because clearly $t \mapsto \mathcal{E}^0(t,q(t))$ is absolutely continuous for any $q \in \AC([0,T]; \mathcal Q^0_{\rm bc})$ and then \eqref{Eq:chainrule:limit} follows from the definition of the local slope.
Note that for the $\Gamma$-convergence result the chain rule property is only required for the limit and not for $h > 0$ but it was necessary for the existence proof in \Cref{Sec:Existence-of-solutions}.

\begin{theorem}[$\Gamma$-convergence]
\label{Thm:General-convergence-result}
For $h > 0$ let $q^h:[0,T] \to Q^h$ be a curve of maximal slope for $(Q^h, D^h, E^h)$ and assume that $q^h_0 := q^h(0)$ satisfies for some $q^0 \in Q^0$,
\begin{equation*}
	\pi^h(q^h_0) \to \pi^0(q^0)
	\quad\text{and}\quad
	E^h(0,q^h_0) \to E^0(0,q^0) < \infty. 
\end{equation*}
Then, there exists $M > 0$ and $h_0 > 0$, such that for all $0 < h \leq h_0$, $\sup_{s,t \in [0,T]} E^h(s,q^h(t)) \leq M$.
Assume the following properties:
\begin{enumerate}[(i)]
	\item \label{Item:General-convergence-result:compactness}(Compactness):
	For all sequences $(q^h) \subset Q^h$ with $\sup_{t \in [0,T]} E^h(t,q^h) \leq M$, any subsequence admits a further subsequence (not relabeled), such that $\pi^h(q^h) \to \pi^0(q)$ for some $q \in Q^0$.
	
	\item \label{Item:General-convergence-result:lsc-ed}(Lower semi-continuity of the energy and dissipation distance):
	For all $(q^h), (\hat q^h) \subset Q^h$ with $\sup_{t \in [0,T]} E^h(t,q^h) \leq M$, $t \in [0,T]$ and $q, \hat q \in Q^0$, we have
	\begin{alignat*}{3}
		\pi^h(q^h) &\to \pi^0(q)  &\quad\Rightarrow\quad&& \liminf_{h \to 0} E^h(t,q^h) &\geq E^0(t,q), \\
		\pi^h(q^h) &\to \pi^0(q),\, \pi^h(\hat q^h) \to \pi^0(\hat q)  &\quad\Rightarrow\quad&& \liminf_{h \to 0} D^h(q^h, \hat q^h) &\geq D^0(q,\hat q).
	\end{alignat*}
		
	\item \label{Item:General-convergence-result:continuity-dt}(Upper semi-continuity of the power):
	For all $(q^h) \subset Q^h$ with $\sup_{t \in [0,T]} E^h(t,q^h) \leq M$ and $q \in Q^0$ and $t \in \Lambda$, we have
	\begin{equation*}
		\pi^h(q^h) \to \pi^0(q)  \quad\Rightarrow\quad \limsup_{h \to 0} \partial_t E^h(t,q^h) \leq \partial_t E^0(t,q).
	\end{equation*}
	
	\item \label{Item:General-convergence-result:lsc-slope}(Lower semi-continuity of the local slope):
	For all $(q^h) \subset Q^h$ with $\sup_{t \in [0,T]} E^h(t,q^h) \leq M$, $t \in [0,T]$ and $q \in Q^0$, we have
		\begin{alignat*}{3}
			\pi^h(q^h) \to \pi^0(q)  \quad\Rightarrow\quad 
			\liminf_{h \to 0}\, \abs*{\partial E^h(t)}_{D^h}(q^h) \geq \abs*{\partial E^0(t)}_{D^0}(q).
		\end{alignat*}
\end{enumerate}
Then, there exists a curve of maximal slope $q: [0,T] \to Q^0$ for $(Q^0,D^0,E^0)$ with $q(0) = q^0$ and a subsequence (not relabeled, independent of $t$), such that for all $t \in [0,T]$,
\begin{equation*}
	\pi^h(q^h(t)) \to \pi^0(q(t)), \qquad
	E^h(t,q^h(t)) \to E^0(t,q(t)). 
\end{equation*}
\end{theorem}

\begin{proof}
By integrating \eqref{Eq:Intro:Curve-of-max-slope}, we obtain that for all $t \in [0,T]$,
\begin{equation}
\label{Eq:Gamma-convergence:energy-inequality}
	E^h\of[\big]{t,q^h(t)} + \int_0^t \tfrac12\abs*{q^h{}'}_{D^h}^2(s) + \tfrac12 \abs*{\partial E^h(s)}_{D^h}^2\of[\big]{q^h(s)} \,\dd s
	\leq E^h(0,q^h_0) + \int_0^t \partial_t E^h\of[\big]{s,q^h(s)} \,\dd s.
\end{equation}
As $E^h(0,q^h_0)$ is bounded, Gronwall's inequality implies in view of the uniform power estimate that the energy is bounded, i.e.\ $\limsup_{h \to 0} \sup_{t \in [0,T]} E^h(t,q^h(t)) < \infty$. Combined with the uniform power estimate, we obtain that there exist $M > 0$ and $h_0 > 0$, such that for all $0 < h \leq h_0$,
\begin{equation*}
	\sup_{s,t \in [0,T]} E^h\of[\big]{s,q^h(t)} \leq M.
\end{equation*}
Furthermore, again due to the uniform power estimate, we observe $\abs{\partial_t E^h(t,q^h(t))} \leq (c+M)\lambda(t)$ for all $t \in \Lambda$, so that the right-hand side of \eqref{Eq:Gamma-convergence:energy-inequality} is bounded.
Hence, there exists $A,S \in \Leb^2([0,T])$ such that $\abs*{q^h{}'}_{D^h} \wto A$ and $\abs*{\partial E^h(\placeholder)}_{D^h}\of{q^h(\placeholder)} \wto S$ in $\Leb^2([0,T])$ up to a subsequence (not relabeled). Since
\begin{equation*}
	\limsup_{h \to 0} D^h\of[\big]{q^h(s), q^h(t)} \leq \limsup_{h \to 0} \int_s^t \abs*{q^h{}'}_{D^h} \leq \int_s^t A \qquad
	\text{for all } 0 \leq s < t \leq T,
\end{equation*}
a generalized version of the Arzel{\`a}--Ascoli theorem (see \cref{Lem:Arzela-Ascoli}) shows that there exists a further subsequence (not relabeled) and a limit $q: [0,T] \to Q^0$ such that $\pi^h(q^h(t)) \to \pi^0(q(t))$ for all $t \in [0,T]$. By lower semi-continuity (\ref{Item:General-convergence-result:lsc-ed}) we have
\begin{equation*}
	D^0\of[\big]{q(s), q(t)} \leq \liminf_{h \to 0} D^h\of[\big]{q^h(s), q^h(t)} \leq \int_s^t A,
\end{equation*}
and thus $q$ is absolutely continuous with $\abs*{q'}_{D^0} \leq A$.
Note that by Fatou's inequality and (\ref{Item:General-convergence-result:lsc-slope}), we obtain $\abs{\partial E^0(\placeholder)}_{D^0}\of{q(\placeholder)} \leq S$.
We show that $q$ is a curve of maximal slope. To obtain the required energy inequality, we pass to the limit in \eqref{Eq:Gamma-convergence:energy-inequality}:
By (\ref{Item:General-convergence-result:lsc-ed}), (\ref{Item:General-convergence-result:continuity-dt}), $\abs{\partial E^0(\placeholder)}_{D^0}\of{q(\placeholder)} \leq S$, and $\abs*{q'}_D \leq A$ we have for all $t \in [0,T]$,
\begin{align*}
	&E^0\of[\big]{t,q(t)} + \int_0^t \tfrac12\abs*{q'}_{D^0}^2(s) + \tfrac12 \abs*{\partial E^0(s)}_{D^0}^2\of[\big]{q(s)} \,\dd s
	\\
	&\qquad\leq E^0\of[\big]{t,q(t)} + \int_0^t \tfrac12 A(s)^2 + \tfrac12 S(s)^2 \,\dd s
	\\
	&\qquad\leq \liminf_{h \to 0} \left( E^h\of[\big]{t,q^h(t)} + \int_0^t \tfrac12\abs*{q^h{}'}_{D^h}^2(s) + \tfrac12 \abs*{\partial E^h(s)}_{D^h}^2\of[\big]{q^h(s)} \,\dd s \right)
	\\
	&\qquad\leq \limsup_{h \to 0} \left( E^h\of[\big]{t,q^h(t)} + \int_0^t \tfrac12\abs*{q^h{}'}_{D^h}^2(s) + \tfrac12 \abs*{\partial E^h(s)}_{D^h}^2\of[\big]{q^h(s)} \,\dd s \right)
	\\
	&\qquad\leq \limsup_{h \to 0}\left( E^h(0,q^h_0) + \int_0^t \partial_t E^h\of[\big]{s,q^h(s)} \,\dd s \right)
	\\
	&\qquad\leq E^0(0,q_0) + \int_0^t \partial_t E^0\of[\big]{s,q(s)} \,\dd s.
\end{align*}
As the right-hand side is finite in view of $E(t,q(t)) \leq \liminf_{h \to 0} E^h(t,q^h(t)) \leq M$, we may apply the chain rule inequality \eqref{Eq:chainrule:limit} and infer that $t \mapsto E^0(t,q(t))$ is absolutely continuous, so that
\begin{equation*}
	\int_0^T \tfrac{\dd}{\dd t} E^0\of[\big]{t,q(t)} + \tfrac12\abs*{q'}_{D^0}^2(t) + \tfrac12 \abs*{\partial E^0(t)}_{D^0}^2\of[\big]{q(t)} - \partial_t E^0\of[\big]{t,q(t)} \,\dd t \leq 0.
\end{equation*}
But as the chain rule implies that the integrand is a.e.\ non-negative, the integrand equals 0 a.e.\ and hence, $q$ is a curve of maximal slope.
Moreover, we obtain that in the chain of inequalities above, the left- and right-hand side are equal, so that all the inequalities must in fact be equalities. Especially, from the second, third and fourth line it is not hard to infer
\begin{equation*}
	E^0\of[\big]{t,q(t)}
	= \liminf_{h \to 0} E^h\of[\big]{t,q^h(t)}
	= \limsup_{h \to 0} E^h\of[\big]{t,q^h(t)}.
\qedhere
\end{equation*}
\end{proof}

The goal of this section is to apply this result to our situation, where we the notion of convergence is given by $y^h \xrightarrow{\rm rod} q$ (respectively, $y^h \xrightarrow{\rm red\, rod} (v,R)$). This can be put into the terminology above by the choice $\tilde Q := \Leb^2(\omega, \R^3) \times \Leb^2(\omega, \R^{3\times 3}) \times \Leb^2(\Omega, \R^{3\times 3})$ and
\begin{equation*}
	\pi^h(y) := \left(y, \grad[h]y, h^{-1}(\sqrt{C_y} - I)\right), \quad
	\pi^0(q) = \pi(v,R,a,\varphi) := \bigl(v,R,\sym\compose*{K_R\mathbf{\bar x} + ae_1}{\bar\grad\varphi}\bigr).
\end{equation*}
We already showed the compactness \cref{Thm:General-convergence-result}~(\ref{Item:General-convergence-result:compactness}) in \cref{Prop:Compactness}. We check the other assumptions in the remainder of this section. The main difficulty here is the lower semi-continuity of the local slope. This is the main achievement of this paper. \cref{Thm:General-convergence-result} establishes convergence up to a subsequence. The convergence for the whole sequence in \cref{Thm:Main-result:convergence-to-limit,Thm:Main-result:convergence-to-reduced} then follows from the fact that the limit is unique, which we showed in \cref{Sec:Rod-models}.

\subsection{Lower semi-continuity of the energy, dissipation and local slope}
\label{Sec:Lsc-slope}

We establish the remaining assumptions of \cref{Thm:General-convergence-result} where our focus is on the lower semi-continuity of the local slope.

\begin{lemma} \label{Lem:Continuity-dtE}
Let $t \in \Lambda$, $(y^h) \subset \mathcal{Q}^h_{\rm bc}$ converge to $q = (v,R,a,\varphi) \in \mathcal{Q}^0_{\rm bc}$ as a rod configuration and assume $\limsup_{h \to 0} \mathcal{E}^h(t,y^h) < \infty$. Then,
\begin{equation}
	\lim_{h \to 0} \partial_t \mathcal{E}^h(t,y^h) = \partial_t \mathcal E^0(t,q) = \partial_t \mathcal E^{\rm red}(t,v,R).
\end{equation}
\end{lemma}

\begin{proof}
The statement follows easily from $\partial_t \mathcal{E}^h(t,y) = -\sprod{\dot l(t)}{y}$ and $\partial_t \mathcal E^0(t,q) = \partial_t \mathcal E^{\rm red}(t,v,R) = -\sprod{\dot l(t)}{v}$.
\end{proof}

\begin{lemma} \label{Lem:LSC-energy}
Let $(y^h) \subset \mathcal{Q}^h_{\rm bc}$ converge to $q = (v,R,a,\varphi) \in \mathcal{Q}^0_{\rm bc}$ as a rod configuration. Then, for all $t \in [0,T]$,
\begin{equation}
	\label{Eq:LSC-energy:total}
	\liminf_{h \to 0} \mathcal{E}^h(t,y^h)
	\geq \mathcal{E}^0(t,q)
	\geq \mathcal{E}^{\rm red}(t,v,R).
\end{equation}
\end{lemma}

\begin{proof}
The proof of \Cref{Lem:LSC-energy} is standard in static elasticity and given essentially in \cite{Neu12}, see also \cite[Prop.~4.8]{NR26}. We repeat it for the convenience of the reader. Without loss of generality, we may restrict to a subsequence, where the energy is bounded. Trivially, we obtain
\begin{equation*}
	\liminf_{h \to 0} \mathcal{E}^h(t,y^h)
	\geq \liminf_{h \to 0} \mathcal{E}^h_\el(y^h) - \lim_{h \to 0} \sprod{l(t)}{y^h}
	= \liminf_{h \to 0} h^{-2} \int_\Omega W_\el(\grad[h]y^h) - \sprod{l(t)}{v}.
\end{equation*}
Therefore, it suffices to establish the lower bound for the elastic energy. Due to the frame-indifference \eqref{Ass:W_FrameIndifference}, we observe $W_\el(\grad[h]y^h) = W_\el(I + hE^h)$ where $E^h := h^{-1}(\sqrt{C_{y^h}} - I)$. As by \Cref{Prop:Compactness} $hE^h \to 0$ in $\Leb^\infty(\Omega, \R^{3\times 3})$ and $(E^h)$ is bounded in $\Leb^2(\Omega, \R^{3\times 3})$, a Taylor expansion of $W_\el$ at $I$, see \eqref{Eq:TaylorEnergy}, then shows
\begin{equation*}
	\liminf_{h \to 0} h^{-2}\int_\Omega W_\el(\grad[h]y^h)
	= \liminf_{h \to 0} h^{-2}\int_\Omega W_\el(I + hE^h)
	= \liminf_{h \to 0} \int_\Omega Q_\el(E^h).
\end{equation*}
Since $y^h \xrightarrow{\rm rod} q$ implies $E^h \wto \sym\compose{K_R\mathbf{\bar x} + ae_1}{\bar\grad\varphi}$ in $\Leb^2(\Omega, \R^{3\times 3})$, we conclude from the lower semi-continuity of convex functionals,
\begin{equation*}
	\liminf_{h \to 0} \int_\Omega Q_\el(E^h)
	\geq \int_\Omega Q_\el\compose{K_R\mathbf{\bar x} + ae_1}{\bar\grad\varphi}
	= \mathcal{E}^0_\el(q).
\end{equation*}
Finally, the definition of $\mathcal{E}^{\rm red}_\el$ shows $\mathcal{E}^0(t,q) \geq \mathcal{E}^{\rm red}(t,v,R)$.
\end{proof}

In view of \eqref{Eq:TaylorDissipation}, using the exact same procedure, it is not hard to establish the lower bound for the dissipation. Thus, we omit the proof of the following statement.

\begin{lemma} \label{Lem:LSC-dissipation}
Let $M > 0$, $(y^h), (\hat{y}^h) \subset \mathcal{Q}^{h,M}_{\rm bc}$ and $q = (v,R,a,\varphi), \hat{q} = (\hat{v},\hat{R},\hat{a},\hat{\varphi}) \in \mathcal{Q}^0_{\rm bc}$ be such that $y^h \xrightarrow{\rm rod} q$ and $\hat y^h \xrightarrow{\rm rod} \hat q$. Then, 
\begin{equation}
	\liminf_{h \to 0} \mathcal{D}^h\of*{y^h, \hat y^h} \geq \mathcal{D}^0(q, \hat{q}) \geq \mathcal{D}^{\rm red}\of*{(v,R), (\hat{v},\hat{R})}.
\end{equation}
\end{lemma}

The remainder of this section is dedicated to establishing the lower bound for the local slope which is the following proposition.

\begin{proposition} \label{Prop:LSC-local-slope}
Let $M > 0$, $t \in [0,T]$ and $(y^h) \subset \mathcal{Q}^{h,M}_{\rm bc}$ converge to $q = (v,R,a,\varphi) \in \mathcal{Q}^0_{\rm bc}$ as a rod configuration. Then,
\begin{equation} \label{Eq:LSC-local-slope:limit}
	\liminf_{h \to 0}\, \abs*{\partial \mathcal{E}^h(t)}_{\mathcal{D}^h}(y^h) \geq \abs*{\partial \mathcal{E}^0(t)}_{\mathcal{D}^0}(q).
\end{equation}
If additionally the compatibility condition \eqref{Ass:C_Compatibility} holds, then also
\begin{equation} \label{Eq:LSC-local-slope:reduced}
	\liminf_{h \to 0}\, \abs*{\partial \mathcal{E}^h(t)}_{\mathcal{D}^h}(y^h) \geq \abs*{\partial \mathcal{E}^{\rm red}(t)}_{\mathcal{D}^{\rm red}}(v,R).
\end{equation}
\end{proposition}

Our goal is to prove \Cref{Prop:LSC-local-slope} by adapting the construction of the mutual recovery sequence in \cite[Prop.~4.11]{NR26}. To see how such a recovery sequence comes into play, we follow the ideas of \cite[Thm.~5.10]{FK20} and \cite[Thm.~6.5]{FM23}. Roughly speaking, the goal is to construct a sequence $(y^h_s)$ of deformations, such that
\begin{equation}
\label{Eq:recovery-sequence-idea}
	\liminf_{h \to 0} \abs{\partial \mathcal E^h(t)}_{\mathcal{D}^h}(y^h)
	\geq \liminf_{h \to 0} \lim_{s \to 0} \frac{\bigl(\mathcal{E}^h(t,y^h) - \mathcal{E}^h(t,y^h_s)\bigr)^+}{\mathcal{D}^h(y^h, y^h_s)}
	\geq \abs{\partial \mathcal E^0(t)}_{\mathcal{D}^0}(q).
\end{equation}
For the construction we use that in \Cref{Rem:Rod-models:Local-slope-maximizing-sequence} we essentially identified the maximizing sequences $(q_s) = (v_s, R_s, a_s, \varphi_s)$ for the local slope in the limit. This motivates to construct $(y^h_s)$ in such a way that $y^h_s \xrightarrow{\mathrm{rod}} q_s$ as $h \to 0$. The limit passage in \eqref{Eq:recovery-sequence-idea} requires uniform estimates with respect to $s$ so that $\mathcal{E}^h(t,y^h) - \mathcal{E}^h(t,y^h_s) \approx \mathcal{E}^0(t,q) - \mathcal{E}^0(t,q_s)$ and $\mathcal{D}^h(y^h, y^h_s) \approx \mathcal{D}^0(q,q_s)$. We provide the following construction:

\begin{lemma} \label{Lem:RecoverySeq}
Let $M > 0$, $h_0 > 0$, $(y^h) \subset \mathcal{Q}^{h,M}_{\rm bc}$ and $q,\hat q \in \mathcal Q^0_{\rm bc}$, $q = (v,R,a,\varphi), \hat q = (\hat v,\hat R,\hat a,\hat\varphi)$. Assume $y^h \xrightarrow{\rm rod} q$. Moreover, let $(v^h, R^h, \phi^h)$ as in \cref{Rem:CosseratAnsatz} and assume that the convergences \eqref{Eq:ConvergencesFromCompactnessAdvanced} hold for some $(\phi, \psi)$. Then, for all $s \in (0,1)$ there exist $y^h_s \in \mathcal{Q}^h_{\rm bc}$ and $R^h_s \in \SobW^{2,\infty}_{\SO(3)}(\omega)$, such that the following statements hold.
Let $q_s = (v_s, R_s, a_s, \varphi_s) \in \mathcal{Q}^0_{\rm bc}$, $G_s$, $G$ and $\hat{G}$ be as in \cref{Rem:Rod-models:Local-slope-maximizing-sequence} and define 
\begin{equation}
	G^h := h^{-1}\left((R^h)^\transp\grad[h]y^h - I\right), \quad
	G^h_s := h^{-1}\left((R^h_s)^\transp\grad[h]y^h_s - I\right) \quad\text{and}\quad
	\Delta^h_s := R^h (G^h_s - G^h).
\end{equation}
Then, we find $\eta: [0,\infty) \to [0,\infty]$ with $\lim_{h \to 0}\eta(h) = 0$, such that for all $s \in (0,1)$ and $h > 0$,
\begin{subequations}
\begin{align}
	\norm*{G^h_s - G^h - G_s + G}_{\Leb^2(\Omega)} &\leq \eta(h)s, \label{Eq:RecoverySeq:EstimateGh}\\
	\norm*{y^h_s - y^h - v_s + v}_{\SobH^1(\Omega)} &\leq \eta(h)s, \label{Eq:RecoverySeq:EstimateYh}\\
	h^{1 - \alpha_C - \alpha_F} \norm*{\Delta^h_s}_{\Leb^\infty(\Omega)} &\leq \eta(h)s, \label{Eq:RecoverySeq:EstimateGhUniform}\\
	h^{1 - \alpha_C} \norm*{\grad \Delta^h_s}_{\Leb^{p_*}(\Omega)} &\leq \eta(h)s, \label{Eq:RecoverySeq:EstimateDelta}\\
	h^{1 + \alpha_\mathbb{F}} \norm*{\grad^2 \Delta^h_s}_{\Leb^{p_*}(\Omega)} &\leq \eta(h)s, \label{Eq:RecoverySeq:EstimateDelta2} \\
	h^{\min\set{\alpha_F,\alpha_{\mathbb F} - \alpha_F}} \norm*{R^h_s (R^h)^\transp - I}_{\SobW^{2,\infty}(\omega)} &\leq \eta(h)s. \label{Eq:RecoverySeq:EstimateR}
\end{align}
\end{subequations}
\end{lemma}

\begin{proof}
\stepemph{Step 1 -- Construction:}
We consider sequences as follows:
\begin{enumerate}[(a)]
	\item $(\phi^h_*) \subset \Cont^\infty(\overline{\omega}, \R^3)$ with $\phi^h_*(0) = 0$, such that
	\begin{equation*}
		\phi^h_* \to \phi \text{ in } \SobH^1(\omega, \R^3)
		\quad\text{and}\quad
		\lim_{h \to 0} h^{1 - \alpha_C - \alpha_F}\norm{\phi^h_*}_{\SobW^{2,\infty}(\omega)} = 0,
	\end{equation*}
	
	\item $(d^h) \subset \Cont^\infty_\mathrm{c}(\omega, \R^3)$ with
	\begin{equation*}
		d^h \to \hat{a}Re_1 - aRe_1 \text{ in } \Leb^2(\omega, \R^3)
		\quad\text{and}\quad
		\lim_{h \to 0} h^{1 - \alpha_C - \alpha_F}\norm{d^h}_{\SobW^{2,\infty}(\omega)} = 0,
	\end{equation*}
	
	\item $(\psi^h) \subset \Cont^\infty_\mathrm{c}(\omega, \Cont^\infty(\overline{S},\R^3))$ with
	\begin{equation*}
		\psi^h \to R\hat\varphi - R\varphi \text{ in } \Leb^2(\omega, \SobH^1(S, \R^3))
		\quad\text{and}\quad
		\lim_{h \to 0} h^{1- \alpha_C - \alpha_F}\norm{\psi^h}_{\SobW^{3,\infty}(\Omega)} = 0,
	\end{equation*}

	\item $(\tilde R^h) \subset \Cont^\infty(\overline{\omega}, \SO(3))$ with $\tilde{R}^h(0) = I$, such that
	\begin{align*}
		&\tilde{R}^h \to \hat{R}R^\transp \text{ in } \SobH^1(\omega, \R^{3\times 3}) \quad\text{and}\quad 
		\lim_{h \to 0} \delta(h) \norm*{\tilde{R}^h}_{\SobW^{3,\infty}(\omega)} = 0,
		\quad\text{where}\\
		&\delta(h) := \begin{aligned}[t]
			\norm{\phi^h - \phi^h_*}_{\Leb^2(\Omega)}
			&+ h^{\frac{2}{p_*}(1-\alpha_C)}
			+ h^{\frac13\alpha_{\mathbb F} - \frac13\max\set*{\frac{p_* - 2}{p_*}(1-\alpha_C), \alpha_F}}
			+ h^{\frac13\alpha_F}
			\\
			&+ h^{1 - \alpha_C - \alpha_F} \bigl(1 + \norm{\phi^h}_{\Leb^\infty(\Omega)} + \norm{\phi^h_*}_{\SobW^{2,\infty}(\Omega)} + \norm{\psi^h}_{\SobW^{2,\infty}(\Omega)} \bigr).
		\end{aligned}
	\end{align*}
\end{enumerate}
Indeed, such sequences exist by standard smoothening arguments. In fact $(\tilde{R}^h)$ can be constructed as in \cite[Lem.~4.10]{NR26} where we note that $\lim_{h \to 0}\delta(h) = 0$ since by \Cref{Prop:Compactness} it holds that $\lim_{h \to 0} h^{1 - \alpha_C - \alpha_F} \norm{\phi^h}_{\Leb^\infty(\Omega)} = 0$.
We introduce $\tilde{R} := \hat{R}R^\transp$, $\hat{R}^h := \tilde{R}^h R^h$, $(v^h_s, R^h_s) = \Phi(K_{R^h} + s(K_{\hat{R}^h} - K_{R^h}))$ (recall \eqref{Eq:Deformation-to-strain-map}) and $\tilde{R}^h_s := R^h_s (R^h)^\transp$.
We show that the claim holds for the deformations
\begin{subequations}
\begin{align}
	y^h_s &:= v^h_s + h R^h_s\mathbf{\bar x} + h \phi^h_s,
\shortintertext{where}
	\phi^h_s &:= \tilde{R}^h_s\bigl(\phi^h - \phi^h_*\bigr) + \int_0^{\mathbf x_1} \tilde{R}^h_s(r) \partial_1\phi^h_*(r) \,\dd r
	+ s\left(h\tilde{R}^h_s \psi^h + \int_0^{\mathbf x_1} \tilde{R}^h_s(r) d^h(r) \,\dd r \right).
\end{align}
\end{subequations}
We remark on this ansatz below in \cref{Rem:RecoverySequence-ansatz}.
\medskip

\stepemph{Step 2 -- Estimates for $G^h_s$:}
A straight-forward calculation shows
\begin{align*}
	G^h_s &= \compose[\big]{K_{R^h_s}\mathbf{\bar x}}{0} + (R^h_s)^\transp \grad[h]\phi^h_s
	\\
	&= \compose[\big]{K_{R^h_s}\mathbf{\bar x}}{0}
		+ (R^h)^\transp\grad[h]\phi^h
		+ (R^h)^\transp \compose[\big]{K_{\tilde{R}^h_s}(\phi^h - \phi^h_* + sh\psi^h) + sd^h + sh\partial_1\psi^h}{s\bar\grad\psi^h}.
\end{align*}
Recalling that $G^h = \compose{K_{R^h}\mathbf{\bar x}}{0} + (R^h)^\transp \grad[h]\phi^h$ and $K_{\tilde{R}^h_s} = R^h (K_{R^h_s} - K_{R^h}) (R^h)^\transp = sK_{\tilde{R}^h}$ (see \eqref{Eq:Rod-models:Local-slope-maximizing-sequence:Strain}), we obtain
\begin{equation}
\label{Eq:RecoverySequence:IdentityDiff}
\begin{aligned}
	G^h_s - G^h &= \compose[\big]{K_{R^h_s}\mathbf{\bar x} - K_{R^h}\mathbf{\bar x}}{0}
		+ (R^h)^\transp \compose[\big]{K_{\tilde{R}^h_s}(\phi^h - \phi^h_* + sh\psi^h) + sd^h + sh\partial_1\psi^h}{s\bar\grad\psi^h}
	\\
	&= s(R^h)^\transp\compose[\big]{K_{\tilde{R}^h}(R^h\mathbf{\bar x} + \phi^h - \phi^h_* + sh\psi^h) + d^h + h\partial_1\psi^h}{\bar\grad\psi^h}.
\end{aligned}
\end{equation}
Furthermore, using \eqref{Eq:Rod-models:Local-slope-maximizing-sequence:Strain}, we get
\begin{equation*}
	G_s - G = s\compose[\big]{R^\transp K_{\tilde R} R\mathbf{\bar x} + \hat{a}e_1 - ae_1}{\bar\grad\hat\varphi - \bar\grad\varphi},
\end{equation*}
so that the construction in Step 1 yields
\begin{align*}
	&s^{-1}\norm{G^h_s - G^h - G_s + G}_{\Leb^2(\Omega)} 
	\\
	&\qquad\leq c_1
	\biggl(\norm[\big]{(R^h)^\transp K_{\tilde{R}^h} R^h - R^\transp K_{\tilde R}R}_{\Leb^2(\omega)}
	\\
	&\qquad\qquad\quad
	+ \norm[\big]{(R^h)^\transp d^h - \hat{a}e_1 + ae_1}_{\Leb^2(\omega)}
	+ \norm[\big]{(R^h)^\transp\bar\grad\psi^h - \bar{\grad}\hat{\varphi} + \bar{\grad}\varphi}_{\Leb^2(\Omega)}
	\\
	&\qquad\qquad\quad
	+ \norm[\big]{\phi^h - \phi^h_*}_{\Leb^2(\Omega)}\norm[\big]{K_{\tilde{R}^h}}_{\Leb^\infty(\omega)}
	+ \norm[\big]{sh\psi^h}_{\Leb^2(\Omega)}\norm[\big]{K_{\tilde{R}^h}}_{\Leb^\infty(\omega)} \biggr) \xrightarrow{h \to 0} 0,
\end{align*}
which establishes \eqref{Eq:RecoverySeq:EstimateGh}.
\medskip

\stepemph{Step 3 -- Estimate for $y^h_s$:}
Recall,
\begin{align*}
	y^h = v^h + hR^h\mathbf{\bar x} + h\phi^h
	\quad\text{and}\quad
	y^h_s = v^h_s + hR^h_s\mathbf{\bar x} + h\phi^h_s.
\end{align*}
By definition of $v^h_s, v^h, v_s, v$ we obtain
\begin{equation*}
	v^h_s - v^h - v_s + v
	= \int_0^{\mathbf{x}_1} (\tilde{R}^h_s - R_sR^\transp)Re_1 + (\tilde{R}^h_s - I)(R^h - R)e_1
\end{equation*}
and thus, \cref{Rem:Estimate_dR_by_K} establishes
\begin{align*}
	&s^{-1}\norm{v^h_s - v^h -v_s + v}_{\SobH^1(\omega)} 
	\\
	&\qquad \leq c_2 s^{-1}\left( \norm{K_{\tilde{R}^h_s} - K_{R_sR^\transp}}_{\Leb^2(\omega)} + \norm{K_{\tilde{R}^h_s}}_{\Leb^\infty(\omega)}\norm{R^h - R}_{\Leb^2(\omega)} \right) \\
	&\qquad = c_2 \left( \norm{K_{\tilde{R}^h} - K_{\tilde R}}_{\Leb^2(\omega)} + \norm{K_{\tilde{R}^h}}_{\Leb^\infty(\omega)}\norm{R^h - R}_{\Leb^2(\omega)} \right) \xrightarrow{h \to 0} 0.
\end{align*}
Moreover, by \cref{Lem:Change-of-variables:homeomorphism} we observe
\begin{equation*}
	hs^{-1}\norm*{R^h_s - R^h}_{\SobH^1(\omega)} 
	\leq c_3 h\norm{K_{\tilde{R}^h}}_{\Leb^2(\omega)} \xrightarrow{h \to 0} 0.
\end{equation*}
Finally, for the rest terms we notice,
\begin{equation*}
	\phi^h_s - \phi^h
	= (\tilde{R}^h_s - I)(\phi^h - \phi^h_*) + \int_0^{\mathbf{x}_1}(\tilde{R}^h_s - I)\partial_1\phi^h_* + s\left(h\tilde{R}^h_s \psi^h + \int_0^{\mathbf x_1} \tilde{R}^h_s(r) d^h(r) \,\dd r \right),
\end{equation*}
so that again by \eqref{Eq:Rod-models:Local-slope-maximizing-sequence:Strain} and \cref{Lem:Change-of-variables:homeomorphism}
\begin{align*}
	hs^{-1}\norm[\big]{\phi^h_s - \phi^h}_{\SobH^1(\Omega)}
	\leq c_4 h \biggl(
		&\norm[\big]{h\psi^h}_{\SobH^1(\Omega)}
		+ \norm[\big]{d^h}_{\Leb^2(\omega)}
		\\
		&+ \norm[\big]{K_{\tilde R^h}}_{\Leb^\infty(\omega)} \Bigl( \norm[\big]{\phi^h - \phi^h_*}_{\SobH^1(\Omega)} + \norm[\big]{\phi^h_*}_{\SobH^1(\omega)} + \norm[\big]{h\psi^h}_{\Leb^2(\Omega)} \Bigr)
	\biggr)
	\xrightarrow{h \to 0} 0.
\end{align*}
Together, these three estimates clearly establish \eqref{Eq:RecoverySeq:EstimateYh}. We conclude also from the Steps 2 and 3 that $y^h_s \xrightarrow{\mathrm{rod}} q_s$.
\medskip

\stepemph{Step 4 -- Estimates for $\Delta^h_s$:}
Equation \eqref{Eq:RecoverySequence:IdentityDiff} shows
\begin{equation*}
	s^{-1}\Delta^h_s
	= \compose[\big]{K_{\tilde{R}^h}(R^h\mathbf{\bar x} + \phi^h - \phi^h_* + sh\psi^h) + d^h + h\partial_1\psi^h}{\bar\grad\psi^h}.
\end{equation*}
We identify the critical term as $K_{\tilde{R}^h} \phi^h$. Since $\lim_{h \to 0} h^{1-\alpha_C - \alpha_F}\norm{K_{\tilde{R}^h}}_{\Leb^\infty(\omega)}\norm{\phi^h}_{\Leb^\infty(\Omega)} = 0$ by construction in Step 1, we conclude \eqref{Eq:RecoverySeq:EstimateGhUniform}.
Calculating the derivatives of $\Delta^h_s$, we obtain for $i,j=1,2,3$
\begin{align*}
	s^{-1}\partial_i \Delta^h_s
	&= \compose*{\begin{aligned}
			& \partial_i K_{\tilde R^h}(R^h\mathbf{\bar x} + \phi^h - \phi^h_* + sh\psi^h) 
			\\
			+ & K_{\tilde R^h}(\partial_i R^h\mathbf{\bar x} +  R^h\partial_i\mathbf{\bar x} + \partial_i\phi^h - \partial_i\phi^h_* + sh\partial_i\psi^h)
			\\
			+ & \partial_i d^h + h\partial_{1i}\psi^h
		\end{aligned}}{\partial_i\bar\grad\psi^h}
\shortintertext{and}
	s^{-1}\partial_{ij} \Delta^h_s
		&= \compose*{\begin{aligned}
				& \partial_{ij} K_{\tilde R^h}(R^h\mathbf{\bar x} + \phi^h - \phi^h_* + sh\psi^h) 
				\\
				+ & \partial_j K_{\tilde R^h}(\partial_i R^h\mathbf{\bar x} +  R^h\partial_i\mathbf{\bar x} + \partial_i\phi^h - \partial_i\phi^h_* + sh\partial_i\psi^h)
				\\
				+ & \partial_i K_{\tilde R^h}(\partial_j R^h\mathbf{\bar x} +  R^h\partial_i\mathbf{\bar x} + \partial_j\phi^h - \partial_j\phi^h_* + sh\partial_j\psi^h)
				\\
				+ & K_{\tilde R^h}(\partial_{ij} R^h\mathbf{\bar x} + \partial_j R^h \partial_i\mathbf{\bar x} + \partial_i R^h \partial_j\mathbf{\bar x} + \partial_{ij}\phi^h - \partial_{ij}\phi^h_* + sh\partial_{ij}\psi^h)
				\\
				+ & \partial_{ij} d^h + h\partial_{1ij}\psi^h
			\end{aligned}}{\partial_{ij}\bar\grad\psi^h}.
\end{align*}
The critical terms are here the terms involving the correction term $\phi^h$ and derivatives of $R^h$. We obtain from \Cref{Prop:Compactness},
\begin{align*}
	h^{1-\alpha_C} \norm{\phi^h}_{\Leb^{p_*}(\Omega)}
	&\leq c_5 h^{1-\alpha_C},
	\\
	h^{1-\alpha_C} \norm{\grad\phi^h}_{\Leb^{p_*}(\Omega)}
	&\leq c_6 h^{\alpha_{p_*}-\alpha_C} 
	= c_6 h^{\frac{2}{p_*}(1-\alpha_C)},
	\\
	h^{1+\alpha_{\mathbb F}} \norm{\grad^2\phi^h}_{\Leb^{p_*}(\Omega)}
	&\leq c_7 h^{1 + \alpha_{\mathbb F} - \max\set*{2 - \alpha_{p_*}, 1 + \alpha_F}} 
	= c_7 h^{\alpha_{\mathbb F} - \max\set*{\frac{p_* - 2}{p_*}(1-\alpha_C), \alpha_F}},
	\\
	h^{1-\alpha_C} \norm{\partial_1 R^h}_{\Leb^{p_*}(\omega)}
	&\leq c_8 h^{\alpha_{p_*} - \alpha_C} 
	= c_8 h^{\frac{2}{p_*}(1-\alpha_C)},
	\\
	h^{1+\alpha_{\mathbb F}} \norm{\partial_1^2 R^h}_{\Leb^{p_*}(\omega)}
	&\leq c_9 h^{\alpha_{p_*} + \alpha_{\mathbb F} - 1}
	= c_9 h^{\alpha_{\mathbb F} - \frac{p_* - 2}{p_*}(1-\alpha_C)}.
\end{align*}
We conclude \eqref{Eq:RecoverySeq:EstimateDelta} and \eqref{Eq:RecoverySeq:EstimateDelta2} in view of the choice of $\delta$ in Step 1.
\medskip

\stepemph{Step 5 -- Estimate for $\tilde{R}^h_s$:} 
Finally, we recall $\tilde{R}^h_s = R^h_s (R^h)^\transp$. By construction of $R^h_s$, we have
\begin{equation*}
	\tilde{R}^h_s - I = \int_0^{\mathbf{x}_1} \partial_1 \tilde{R}^h_s, \quad
	\partial_1 \tilde{R}^h_s = \tilde{R}^h_s K_{\tilde R^h_s} = s \tilde{R}^h_s K_{\tilde R^h}, \quad
	\partial_1^2 \tilde{R}^h_s = s \tilde{R}^h_s \bigl( s K_{\tilde{R}^h}^2 + \partial_1 K_{\tilde{R}^h} \bigr).
\end{equation*}
Thus, \eqref{Eq:RecoverySeq:EstimateR} follows from $\lim_{h \to 0}h^{\min\set{\alpha_F, \alpha_{\mathbb F} - \alpha_F}/2}\norm{\tilde R^h}_{\SobW^{2,\infty}(\omega)} = 0$.
\end{proof}

\begin{remark} \label{Rem:RecoverySequence-ansatz}
It is plausible to construct the recovery sequence ansatz $y^h_s$ in \cref{Lem:RecoverySeq} again via a (corrected) Cosserat ansatz. The choices for $v^h_s$, $R^h_s$ and $\phi^h_s$ are motivated by the goal to define an ansatz that is (approximately) linear in $G$, i.e.\
\begin{equation}
	G^h_s \approx G^h + s (G[\hat q] - G[q]),
\end{equation}
(recall \cref{Rem:CosseratAnsatzRefined}). In fact, for the correction term $\phi^h_s$, we observe that in the situation of \cref{Prop:Compactness} the relevant correctors $(a,\varphi)$ depend linearly on the limit of $(R^h)^\transp\grad[h]\phi^h$. This leads to an ansatz of the form,
\begin{equation*}
	(R^h_s)^\transp\phi^h_s \approx (R^h)^\transp\phi^h + s(\hat R^\transp\hat \phi - R^\transp\phi),
\end{equation*}
where $\hat \phi$ and $\phi$ belong to the corrected Cosserat ansatz \eqref{Eq:CorrectedCosseratAnsatz}. However, we need a refined treatment of the rotational ansatz, as e.g.\ in
\begin{equation*}
	(R^h_s)^\transp\grad[h]\big( R^h_s (R^h)^\transp\phi^h \big) 
	= (R^h)^\transp \grad[h]\phi^h + (R^h)^\transp\compose*{K_{\tilde R^h_s} \phi^h}{0},
\end{equation*}
the last term does not vanish as $h \to 0$. In contrast with the refined ansatz in \Cref{Lem:RecoverySeq} the corresponding term vanishes in the limit passage for \eqref{Eq:RecoverySeq:EstimateGh}.
\end{remark}

We use this construction to prove \cref{Prop:LSC-local-slope}.

\begin{proof}[Proof of \cref{Prop:LSC-local-slope}]
We start with the proof of \eqref{Eq:LSC-local-slope:limit}. Therefore let $(y^h) \subset \mathcal{Q}^{h,M}_{\rm bc}$ converge to $q = (v,R,a,\varphi) \in \mathcal{Q}^0_{\rm bc}$ as a rod configuration. Let $(v^h, R^h, \phi^h)$ be as in \cref{Rem:CosseratAnsatz} and consider $q \neq \hat{q} \in \mathcal{Q}^0_{\rm bc}$. By restricting to a subsequence (not relabeled) we may assume that the convergences \eqref{Eq:ConvergencesFromCompactnessAdvanced} hold for some $(\phi,\psi)$ and that the subsequence attains $\liminf_{h \to 0}\abs*{\partial \mathcal{E}^h}(y^h)$ as a limit. 
Hence, the assumptions of \Cref{Lem:RecoverySeq} are satisfied. Recalling the definitions and results of \Cref{Lem:RecoverySeq}, we show that there exists some $\hat\eta: [0,\infty) \to [0,\infty]$ with $\lim_{h \to 0} \hat\eta(h) = 0$, such that, for all $h > 0$ and $s \in (0,1)$,
\begin{subequations} \label{Eq:LSC-local-slope:Inequalities}
\begin{align}
	\mathcal{D}^h(y^h,y^h_s) &\leq \mathcal{D}^0(q,q_s) + \hat\eta(h)s, 
		\label{Eq:LSC-local-slope:Inequality-dissipation} \\
	\mathcal{E}^h(t,y^h) - \mathcal{E}^h(t,y^h_s) &\geq \mathcal{E}^0(t,q) - \mathcal{E}^0(t,q_s) - \hat\eta(h)s.
		\label{Eq:LSC-local-slope:Inequality-energy}
\end{align}
\end{subequations}
For an easier presentation we allow for the precise definition of $\hat\eta$ to change throughout the proof from line to line.
Once we established these inequalities, we can easily conclude the statement. Indeed, then especially $\mathcal D^h(y^h, y^h_s) \to 0$ as $s \to 0$, and so, by definition of $\abs{\partial \mathcal{E}^h(t)}_{\mathcal{D}^h}$ and \Cref{Rem:Rod-models:Local-slope-maximizing-sequence},
\begin{align*}
	\abs[\big]{\partial \mathcal{E}^h(t)}_{\mathcal{D}^h}(y^h) 
	&\geq \limsup_{s \to 0} \frac{\left(\mathcal{E}^h(t,y^h) - \mathcal{E}^h(t,y^h_s)\right)^+}{\mathcal{D}^h(y^h, y^h_s)}
	\\
	&\geq \limsup_{s \to 0} \frac{\left(\mathcal{E}^0(t,q) - \mathcal{E}^0(t,q_s)\right)^+ - \hat{\eta}(h)s}{\mathcal{D}^0(q, q_s) + \hat{\eta}(h)s}
	\\
	&\geq \frac{\left(\mathcal{E}^0(t,q) - \mathcal{E}^0(t,\hat q) + \frac{\lambda}{2} \mathcal D^0(q,\hat q)^2\right)^+ - \hat{\eta}(h)}{\mathcal{D}^0(q, \hat q) + \hat{\eta}(h)}.
\end{align*}
Hence, 
\begin{equation*}
	\liminf_{h \to 0}\, \abs[\big]{\partial \mathcal{E}^h(t)}_{\mathcal{D}^h}(y^h)
	\geq \frac{\left(\mathcal{E}^0(t,q) - \mathcal{E}^0(t,\hat q) + \frac{\lambda}{2} \mathcal D^0(q,\hat q)^2\right)^+}{\mathcal{D}^0(q, \hat q)},
\end{equation*}
and since this holds for arbitrary $\hat q \neq q$, we conclude $\liminf_{h \to 0} \abs{\partial \mathcal{E}^h(t)}_{\mathcal{D}^h}(y^h) \geq \abs{\partial \mathcal{E}^0(t)}_{\mathcal{D}^0}(q)$ from \Cref{Lem:Rod-models:Characterization-slope}~(\ref{Item:Rod-models:Characterization-slope}).
In the remainder of this proof we show \eqref{Eq:LSC-local-slope:Inequalities} and afterwards establish \eqref{Eq:LSC-local-slope:reduced}.
\medskip

\stepemph{Step 1 -- Dissipation:}
Let $R^h_* \in \Leb^\infty(\Omega, \SO(3))$ be a rotation field obtained through the polar decomposition, such that
\begin{equation*}
	I + hE^h
	:= \sqrt{C_{y^h}}
	= R^h_* (R^h)^\transp \grad[h]y
	= R^h_* (I + hG^h).
\end{equation*}
Let furthermore $E^h_s := h^{-1}( R^h_*(I + hG^h_s) - I)$. We observe $R^h_* \to I$ in $\Leb^2(\Omega, \R^{3\times 3})$ and due to \eqref{Eq:RecoverySeq:EstimateGh} and \eqref{Eq:RecoverySeq:EstimateGhUniform},
\begin{equation*}
	\norm{E^h - E^h_s - G_s + G}_{\Leb^2(\Omega)} \leq \hat{\eta}(h)s, \qquad
	h\norm{E^h - E^h_s}_{\Leb^\infty(\Omega)} \leq \hat{\eta}(h)s.
\end{equation*}
Since moreover $hE^h \to 0$ in $\Leb^\infty(\Omega, \R^{3\times 3})$ by \Cref{Prop:Compactness}, we may apply the Taylor expansion \eqref{Eq:TaylorDissipation} to show
\begin{align*}
	&\abs*{\mathcal{D}^h(y^h, y^h_s)^2 - \int_\Omega Q_\vi\of{E^h - E^h_s}}
	\\
	&\qquad = h^{-2}\abs*{\int_\Omega \mathscr{D}_\vi^2(I + hE^h, I + hE^h_s) - \int_\Omega Q_\vi\of{hE^h - hE^h_s}} 
	\\
	&\qquad \leq c_1 \bigl( \norm{E^h}_{\Leb^2(\Omega)} + \norm{E^h_s}_{\Leb^2(\Omega)}\bigr) h \norm{\Delta^h_s}_{\Leb^\infty(\Omega)} \norm{E^h - E^h_s}_{\Leb^2(\Omega)}
	\leq \hat{\eta}(h)s^2.
\end{align*}
Let $\mathbb{D}_\vi \in \R^{3\times 3 \times 3\times 3}$ denote the fourth order tensor associated to the quadratic form $Q_\vi$. Then, we obtain
\begin{align*}
	&\abs*{\int_\Omega Q_\vi\of{E^h - E^h_s} - Q_\vi\of{G - G_s}}
	\\
	&\qquad= \abs*{\int_\Omega (E^h - E^h_s + G - G_s): \mathbb{D}_\vi (E^h - E^h_s - G + G_s)}
	\\
	&\qquad \leq c_2 \norm[\big]{E^h - E^h_s + G - G_s}_{\Leb^2(\Omega)} \norm[\big]{E^h - E^h_s - G + G_s}_{\Leb^2(\Omega)}
	\leq \hat{\eta}(h)s^2.
\end{align*}
Since $\int_\Omega Q_\vi\of{G - G_s} = \mathcal{D}^0(q, q_s)^2$, combining the two estimates yields \eqref{Eq:LSC-local-slope:Inequality-dissipation}.
\medskip

\stepemph{Step 2 -- Elastic energy:}
We split \eqref{Eq:LSC-local-slope:Inequality-energy} into several statements and show
\begin{subequations}
\begin{align}
	\mathcal E^h_\el(y^h) - \mathcal E^h_\el(y^h_s) &\geq \mathcal E^0_\el(q) - \mathcal E^0_\el(q_s) - \hat\eta(h)s,
		\label{Eq:LSC-local-slope:Inequality-elastic-energy} \\
	\sprod*{l(t)}{y^h_s - y^h} &\geq \sprod*{l(t)}{v_s - v} - \hat\eta(h)s,
		\label{Eq:LSC-local-slope:Inequality-loading-term} \\
	\mathcal E^h_{\rm sg}(y^h) - \mathcal E^h_{\rm sg}(y^h_s) &\geq  - \hat\eta(h)s.
		\label{Eq:LSC-local-slope:Inequality-strain-gradients}
\end{align}
\end{subequations}
First, we draw our attention to the hyperelastic energy. Analogously as for the dissipation distance in Step 1, we may apply the Taylor expansion \eqref{Eq:TaylorEnergy} and the estimates \eqref{Eq:RecoverySeq:EstimateGh} and \eqref{Eq:RecoverySeq:EstimateGhUniform} to show
\begin{align*}
	&\abs*{\mathcal{E}^h_\el(y^h) - \mathcal{E}^h_\el(y^h_s) - \int_\Omega Q_\el\of{E^h} - Q_\el\of{E^h_s}}
	\\
	&\qquad= h^{-2}\abs*{\int_\Omega W_\el\of{I + hE^h} - W_\el\of{I + hE^h_s} - Q_\el\of{hE^h} + Q_\el\of{hE^h_s}}
	\\
	&\qquad\leq c_3 \bigl( \norm{E^h}_{\Leb^2(\Omega)}^2 + \norm{E^h_s}_{\Leb^2(\Omega)}^2 \bigr) h\norm{\Delta^h_s}_{\Leb^\infty(\Omega)} \leq \hat{\eta}(h)s,
\end{align*}
and use that $Q_\el$ is a quadratic form to establish
\begin{align*}
	&\abs*{\int_\Omega Q_\el(E^h) - Q_\el(E^h_s) - Q_\el(G) + Q_\el(G_s)} \\
	&\qquad= \abs*{\int_\Omega (E^h + E^h_s + G_s - G): \mathbb{C}_\el (E^h - E^h_s - G + G_s)
		+ 2(G - G_s): \mathbb{C}_\el (E^h - G)}
	\\
	&\qquad\leq c_4 \norm{E^h - E^h_s - G + G_s}_{\Leb^2(\Omega)} + 2s\abs*{\int_\Omega (G - \hat{G}):\mathbb{C}_\el (E^h - G)}
	\\
	&\qquad\leq \hat{\eta}(h)s + 2s\abs*{\int_\Omega (G - \hat{G}):\mathbb{C}_\el (E^h - G)}
	\\
	&\qquad= \hat{\eta}(h)s + 2s\abs*{\int_\Omega (G - \hat{G}):\mathbb{C}_\el (E^h - \sym G)}.
\end{align*}
The latter term converges to $0$ as $h \to 0$ by the weak convergence $E^h \wto \sym G$ in $\Leb^2(\Omega, \R^{3\times 3})$. Hence, we conclude \eqref{Eq:LSC-local-slope:Inequality-elastic-energy}.
\medskip

\stepemph{Step 3 -- Loading term:} 
For the loading term we conclude \eqref{Eq:LSC-local-slope:Inequality-loading-term} from the observation that by \eqref{Eq:RecoverySeq:EstimateYh},
\begin{align*}
	&\abs*{\sprod*{l(t)}{y^h} - \sprod{l(t)}{y^h_s} - \sprod*{l(t)}{v} + \sprod*{l(t)}{v_s}} \\
	&\qquad= \abs*{\sprod*{l(t)}{y^h - y^h_s - v + v_s}}
	\leq c_5 \norm*{y^h - y^h_s - v + v_s}_{\SobH^1(\Omega)} \leq \hat{\eta}(h)s.
\end{align*}

\stepemph{Step 4 -- Strain gradient terms:}
It remains to estimate the strain gradient terms. The assumptions \nameref{Ass:HE} and the H{\"o}lder inequality yield
\begin{align*}
		&h^{-\alpha_C p_*}\int_\Omega H_\el^C(\grad C_{y^h}) - H_\el^C(\grad C_{y^h_s})
		\\
		&\qquad\geq h^{-\alpha_C p_*} \int_\Omega \dif H_\el^C(\grad C_{y^h_s})\bigl[\grad C_{y^h} - \grad C_{y^h_s}\bigr]
		\\
		&\qquad\geq -c_6 \norm[\big]{h^{-\alpha_C}\grad C_{y^h_s}}_{\Leb^{p_*}(\Omega)}^{p_* - 1} h^{-\alpha_C}\norm[\big]{\grad C_{y^h_s} - \grad C_{y^h}}_{\Leb^{p_*}(\Omega)},
	\\[0.5\baselineskip]
		&h^{\alpha_F p_*}\int_\Omega H_\el^F(\grad \grad[h] y^h) - H_\el^F(\grad \grad[h]y^h_s)
		\\
		&\qquad= h^{\alpha_F p_*}\int_\Omega H_\el^F(\tilde R^h_s\grad \grad[h] y^h) - H_\el^F(\grad \grad[h]y^h_s)
		\\
		&\qquad\geq h^{\alpha_F p_*}\int_\Omega \dif H_\el^F(\grad \grad[h]y^h_s)\bigl[\tilde R^h_s\grad \grad[h]y^h - \grad\grad[h]y^h_s\bigr] 
		\\
		&\qquad\geq -c_7 \bigl(\norm[\big]{h^{\alpha_F}\grad \grad[h]y^h_s}_{\Leb^{p_*}(\Omega)}^{p_* - 1} + 1\bigr) h^{\alpha_F}\norm[\big]{\grad \grad[h]y^h_s - \tilde R^h_s\grad \grad[h] y^h}_{\Leb^{p_*}(\Omega)},
	\\[0.5\baselineskip]
		&h^{\alpha_{\mathbb F} p_*}\int_\Omega H_\el^{\mathbb F}(\grad^2 \grad[h] y^h) - H_\el^{\mathbb F}(\grad^2 \grad[h]y^h_s)
		\\
		&\qquad= h^{\alpha_{\mathbb F} p_*}\int_\Omega H_\el^{\mathbb F}(\tilde R^h_s\grad^2 \grad[h]y^h) - H_\el^{\mathbb F}(\grad^2 \grad[h]y^h_s)
		\\
		&\qquad\geq h^{\alpha_{\mathbb F} p_*}\int_\Omega \dif H_\el^{\mathbb F}(\grad^2 \grad[h]y^h_s)\bigl[\tilde R^h_s\grad^2 \grad[h]y^h - \grad^2\grad[h]y^h_s\bigr]
		\\
		&\qquad\geq -c_8 \bigl(\norm[\big]{h^{\alpha_{\mathbb F}}\grad^2 \grad[h]y^h_s}_{\Leb^{p_*}(\Omega)}^{p_* - 1} + 1\bigr) h^{\alpha_{\mathbb F}}\norm[\big]{\grad^2 \grad[h]y^h_s - \tilde R^h_s\grad^2 \grad[h] y^h}_{\Leb^{p_*}(\Omega)}.
\end{align*}
To conclude that the right-hand sides vanish, we show the estimates
\begin{subequations}
\label{Eq:recovery-sequence:estimates-strain-gradients}
\begin{align}
	h^{-\alpha_C} \norm[\big]{\grad C_{y^h_s} - \grad C_{y^h}}_{\Leb^{p_*}(\Omega)} 
	& \leq \hat{\eta}(h)s,
	\\
	h^{\alpha_F} \norm[\big]{\grad \grad[h]y^h_s - \tilde R^h_s \grad \grad[h] y^h}_{\Leb^{p_*}(\Omega)}
	& \leq \hat{\eta}(h)s,
	\\
	h^{\alpha_{\mathbb F}} \norm[\big]{\grad^2 \grad[h]y^h_s - \tilde R^h_s\grad^2 \grad[h] y^h}_{\Leb^{p_*}(\Omega)}
	& \leq \hat{\eta}(h)s,
\end{align}
\end{subequations}
Indeed, this suffices, since $h^{-\alpha_C}\grad C_{y^h}$, $h^{\alpha_F}\grad\grad[h]y^h$ and $h^{\alpha_{\mathbb F}}\grad^2\grad[h]y^h$ are bounded in $\Leb^{p_*}(\Omega)$ by \Cref{Prop:Compactness}.
The definition $\Delta^h_s = R^h(G^h_s - G^h)$ yields $\grad[h]y^h_s = \tilde{R}^h_s(\grad[h]y^h + h\Delta^h_s)$. Then, a straight-forward computation shows
\begin{align*}
		\partial_i C_{y^h_s} - \partial_i C_{y^h}
		&= 2h\sym\of[\big]{(\partial_i \Delta^h_s)^\transp \grad[h]y^h + (\Delta^h_s)^\transp \partial_i\grad[h]y^h} + 2h^2 \sym\of[\big]{(\partial_i\Delta^h_s)^\transp \Delta^h_s},
	\\[0.5\baselineskip]
		\partial_i \grad[h]y^h_s - \tilde{R}^h_s\partial_i \grad[h]y^h
		&= (\partial_i \tilde{R}^h_s)(\grad[h]y^h + h\Delta^h_s) + h\tilde{R}^h_s \partial_i \Delta^h_s,
	\\[0.5\baselineskip]
		\partial_{ij}\grad[h]y^h_s - \tilde R^h_s \partial_{ij} \grad[h]y^h
		&=  \partial_j \tilde R^h_s \partial_i \grad[h]y^h +  \partial_i\tilde R^h_s (\partial_j\grad[h]y^h + h\partial_j\Delta^h_s) \\
		&\qquad + \partial_{ij}\tilde R^h_s (\grad[h]y^h + h\Delta^h_s) + h\tilde R^h_s \partial_{ij} \Delta^h_s + h\partial_j\tilde R^h_s \partial_i \Delta^h_s.
\end{align*}
We conclude \eqref{Eq:recovery-sequence:estimates-strain-gradients} from the following estimates which are consequences of \Cref{Prop:Compactness,Lem:RecoverySeq}:
\begin{align*}
	\norm[\big]{\tilde R^h_s}_{\Leb^\infty(\omega)}
	+ \norm[\big]{\grad[h]y^h}_{\Leb^\infty(\Omega)}
	+ h^{\alpha_F}\norm[\big]{\grad\grad[h]y^h}_{\Leb^{p_*}(\Omega)}
	+ h^{\alpha_{\mathbb F}}\norm[\big]{\grad^2\grad[h]y^h}_{\Leb^{p_*}(\Omega)} &\leq c_9,
	\\
	h^{1-\alpha_C - \alpha_F}\norm[\big]{\Delta^h_s}_{\Leb^\infty(\Omega)} 
	+ h^{1- \alpha_C} \norm[\big]{\grad\Delta^h_s}_{\Leb^{p_*}(\Omega)} 
	+ h^{1 + \alpha_{\mathbb F}} \norm[\big]{\grad^2\Delta^h_s}_{\Leb^{p_*}(\Omega)} &\leq \hat{\eta}(h)s,
	\\
	h^{\min\set{\alpha_F, \alpha_{\mathbb F} - \alpha_F}}\norm[\big]{\partial_1 \tilde R^h_s}_{\Leb^\infty(\omega)} + h^{\alpha_{\mathbb F}}\norm[\big]{\partial_1^2 \tilde R^h_s}_{\Leb^\infty(\omega)} &\leq \hat{\eta}(h)s.
\end{align*}

\stepemph{Step 5 -- Lower bound for $\mathcal{E}^{\rm red}$:}
It remains to show that the lower bound also holds for the reduced energy and dissipation, provided \nameref{Ass:C} holds. In fact we show that in this case, we have for all $q = (v,R,a,\varphi) \in \mathcal{Q}^0_{\rm bc}$,
\begin{equation*}
	\abs{\partial \mathcal{E}^0(t)}_{\mathcal{D}^0}(q)
	\geq \abs{\partial \mathcal{E}^\mathrm{red}(t)}_{\mathcal{D}^\mathrm{red}}(v,R).
\end{equation*}
Let therefore $(v_s, R_s) \subset \mathcal{Q}^{\rm red}_{\rm bc}$ be such that $0 \neq \mathcal{D}^\mathrm{red}((v,R),(v_s, R_s)) \to 0$ and
\begin{equation*}
	\abs{\partial \mathcal{E}^\mathrm{red}(t)}_{\mathcal{D}^\mathrm{red}}(v,R)
	= \lim_{s \to 0} \frac{\bigl(\mathcal{E}^\mathrm{red}(t,v,R) - \mathcal{E}^\mathrm{red}(t,v_s, R_s)\bigr)^+}{\mathcal{D}^\mathrm{red}((v,R),(v_s,R_s))}.
\end{equation*}
Let $(a^*, \varphi^*)$, $(a_s^*,\varphi_s^*)$ be the correctors, i.e.,
\begin{equation*}
	\sym\compose*{a^* e_1}{\bar{\grad}\varphi^*} = \chi_\el(K_R) = \chi_\vi(K_R), \quad
	\sym\compose*{a_s^* e_1}{\bar{\grad}\varphi_s^*} = \chi_\el(K_{R_s}) = \chi_\vi(K_{R_s}).
\end{equation*}
Set $a_s := a + a_s^* - a^*$ and $\varphi_s := \varphi + \varphi_s^* - \varphi^*$ and define $q_s := (v_s,R_s,a_s,\varphi_s) \in \mathcal{Q}^0_{\rm bc}$. We observe from the linearity of the correctors for the dissipation,
\begin{equation*}
	\mathcal{D}^0(q,q_s)^2 
	= \int_\Omega Q_\vi\compose[\big]{K_{R_s}\mathbf{\bar{x}} - K_R\mathbf{\bar{x}} + a^*_se_1 - a^*e_1}{\bar{\grad}\varphi_s^* - \bar{\grad}\varphi^*}
	= \mathcal{D}^\mathrm{red}((v,R), (v_s, R_s))^2.
\end{equation*}
Especially, this shows $\mathcal{D}^0(q,q_s) \to 0$.
Moreover, introduce the notation $E[\tilde{K},\tilde{a},\tilde{\varphi}] = \compose{\tilde{K}\mathbf{\bar x} + \tilde{a}e_1}{\bar\grad\tilde\varphi}$. We obtain,
\begin{align*}
	&\mathcal{E}^0_\el(q) - \mathcal{E}^0_\el(q_s)
	= \int_\Omega Q_\el\compose*{K_R\mathbf{\bar{x}} + ae_1}{\bar{\grad}\varphi} - Q_\el\compose*{K_{R_s}\mathbf{\bar{x}} + a_s e_1}{\bar{\grad}\varphi_s}
	\\
	&\quad = \int_\Omega E[K_R - K_{R_s}, a - a_s, \varphi - \varphi_s] : \mathbb{C}_\el E[K_R + K_{R_s}, a + a_s, \varphi + \varphi_s]
	\\
	&\quad = \int_\Omega E[K_R - K_{R_s}, a^* - a_s^*, \varphi^* - \varphi_s^*] : \mathbb{C}_\el E[K_R + K_{R_s}, a + a_s, \varphi + \varphi_s]
	\\
	&\quad = \int_\Omega E[K_R - K_{R_s}, a^* - a_s^*, \varphi^* - \varphi_s^*] : \mathbb{C}_\el E[K_R + K_{R_s}, a^* + a_s^*, \varphi^* + \varphi_s^*]
	\\
	&\quad = \int_\Omega Q_\el\compose*{K_R\mathbf{\bar{x}} + a^*e_1}{\bar{\grad}\varphi^*} - Q_\el\compose*{K_{R_s}\mathbf{\bar{x}} + a_s^*e_1}{\bar{\grad}\varphi_s^*}
	= \mathcal{E}^\mathrm{red}_\el(v,R) - \mathcal{E}^\mathrm{red}_\el(v_s,R_s),
\end{align*}
where the fourth equality is due to the corrector equation for $\chi_\el(K_R - K_{R_s})$. Hence, we conclude
\begin{align*}
	\abs{\partial \mathcal{E}^\mathrm{red}(t)}_{\mathcal{D}^\mathrm{red}}(v,R)
	&= \lim_{s \to 0} \frac{\bigl(\mathcal{E}^\mathrm{red}(t,v,R) - \mathcal{E}^\mathrm{red}(t,v_s, R_s)\bigr)^+}{\mathcal{D}^\mathrm{red}(v,R;v_s,R_s)}
	\\
	&= \lim_{s \to 0} \frac{\bigl(\mathcal{E}^0(t,q) - \mathcal{E}^0(t,q_s)\bigr)^+}{\mathcal{D}^0(q,q_s)}
	\leq \abs{\partial \mathcal{E}^0(t)}_{\mathcal{D}^0}(q),
\end{align*}
where the inequality follows from the definition of the local slope.
\end{proof}

\appendix
\section{Appendix}

\subsection{Characterization of inextensible Cosserat rods by the bending-torsion strain}
\label{Sec:Change-of-variables}

In \cref{Sec:Rod-models} we determined the structure of the viscoelastic rod models by using a change of variables from the inextensible Cosserat rod configurations $(v,R)$ to the associated bending-torsion strain $K$ via the bijective map 
\begin{equation*}
	\Phi: \Leb^2(\omega, \R^{3\times 3}_{\skewMat}) \to \mathcal A^{\rm rod}_{\rm bc} \subset \SobH^2(\omega, \R^3) \times \SobH^1(\omega, \R^{3\times 3}),
\end{equation*}
as defined in \eqref{Eq:Deformation-to-strain-map}, that is $\Phi(K) = (\Phi_v(K),\Phi_R(K)) = (v_K, R_K)$, where $R_K$ is the unique solution to the initial value problem
\begin{equation}
	\partial_1 R = R K, \qquad
	R(0) = R_{\rm bc},
\end{equation}
and $v_K(x_1) := v_{\rm bc} + \int_0^{x_1} R_K e_1$. In this section, we discuss the properties of this map.

\begin{lemma} \label{Lem:Change-of-variables:homeomorphism}
The map $\Phi$ as defined in \eqref{Eq:Deformation-to-strain-map} is a homeomorphism. More precisely, given $M > 0$ there exists $C > 0$, such that for all $K, \hat K \in \Leb^2(\omega, \R^{3\times 3}_{\skewMat})$ with $\norm{K}_{\Leb^2(\omega)} \leq M$,
\begin{equation} \label{Eq:Change-of-variables:Estimate-R-K}
	\tfrac{1}{C} \norm*{K - \hat K}_{\Leb^2(\omega)}
	\leq \norm*{R_K - R_{\hat K}}_{\SobH^1(\omega)} 
	\leq  C \norm*{K - \hat K}_{\Leb^2(\omega)}.
\end{equation}
Furthermore, $\Phi$ is twice continuously Frech{\'e}t differentiable with derivatives
\begin{subequations}
\begin{align}
	\dif \Phi_R(K)[K'](x_1) &= \left( \int_0^{x_1} R_K K' R_K^\transp \right) R_K(x_1), \\
	\dif \Phi_v(K)[K'](x_1) &= \int_0^{x_1} \dif \Phi_R(K)[K'](s) e_1 \,\dd s,
\end{align}
\end{subequations}
and there exists $c > 0$, such that for all $K \in \Leb^2(\omega, \R^{3\times 3}_{\skewMat})$
\begin{subequations}
\begin{alignat}{3}
	&\norm{\dif \Phi_v(K)}_{\mathcal L(\Leb^2,\SobW^{1,\infty})} &&+
	\norm{\dif^2 \Phi_v(K)}_{\mathcal L(\Leb^2, \mathcal L( \Leb^2, \SobW^{1,\infty}))} &&\leq c, 
	\\
	&\norm{\dif \Phi_v(K)}_{\mathcal L(\Leb^2,\SobH^2)} &&+
	\norm{\dif^2 \Phi_v(K)}_{\mathcal L(\Leb^2, \mathcal L( \Leb^2, \SobH^2))} &&\leq c(1 + \norm{K}_{\Leb^2(\omega)}),
\intertext{and}
	&\norm{\dif \Phi_R(K)}_{\mathcal L(\Leb^2,\Leb^\infty)} &&+
	\norm{\dif^2 \Phi_R(K)}_{\mathcal L(\Leb^2, \mathcal L( \Leb^2, \Leb^\infty))} &&\leq c,
	\\
	&\norm{\dif \Phi_R(K)}_{\mathcal L(\Leb^2,\SobH^1)} &&+
	\norm{\dif^2 \Phi_R(K)}_{\mathcal L(\Leb^2, \mathcal L( \Leb^2, \SobH^1))} &&\leq c(1 + \norm{K}_{\Leb^2(\omega)}).
\end{alignat}
\end{subequations}
\end{lemma}

\begin{proof}
It is easy to see that $\Phi$ is bijective, as we can write down the inverse map explicitly as $(v,R) \mapsto K_R := R^\transp \partial_1 R$. We only show the first estimate in \eqref{Eq:Change-of-variables:Estimate-R-K} as the second inequality follows from the bounds on the derivative $\dif \Phi_R$ which we discuss below. Let $K, \hat K \in \Leb^2(\omega, \R^{3\times 3}_{\skewMat})$. Then, using the Sobolev inequality, we obtain,
\begin{align*}
	\norm*{K - \hat K}_{\Leb^2(\omega)}
	= \norm*{R_K^\transp \partial_1 R_K - R_{\hat K}^\transp \partial_1 R_{\hat K}}_{\Leb^2(\omega)}
	&= \norm*{\partial_1 R_K - \partial_1 R_{\hat K} + (R_{\hat K}R_K^\transp - I)\partial_1 R_K}_{\Leb^2(\omega)} \\	
	&\leq c_1 \left(1 + \norm*{\partial_1 R}_{\Leb^2(\omega)}\right) \norm*{R - \hat{R}}_{\SobH^1(\omega)}.
\end{align*}
Furthermore, we formally differentiate the differential equation $\partial_1 \Phi_R(K) = \Phi_R(K)K$ to determine candidates for the derivative of $\Phi_R$ and arrive at,
\begin{alignat*}{3}
	\partial_1 \chi_{K'} &= \chi_{K'}K + R_K K', 
	&& \quad \chi_{K'}(0) = 0,
	&&\qquad \chi_{K'} := \dif \Phi_R(K)[K'],
	\\
	\partial_1 \psi &= \psi K + \chi_{K_1'}K_2' + \chi_{K_2'}K_1', 
	&& \quad\psi(0) = 0,
	&& \qquad \psi := \dif^2 \Phi_R(K)[K_1', K_2'].
\end{alignat*}
These differential equation admit the unique solutions,
\begin{align*}
	\chi_{K'}(x_1) &:= \left( \int_0^{x_1} R_K K' R_K^\transp \right) R_K(x_1), \\
	\psi(x_1) &:= \begin{aligned}[t]
		\chi_{K_1'}(x_1) R_K(x_1)^\transp \chi_{K_2'}(x_1)
		&+ \chi_{K_2'}(x_1) R_K(x_1)^\transp \chi_{K_1'}(x_1) \\
		&- \left( \int_0^{x_1} R_K K_1' R_K^\transp \chi_{K_2'} + R_K K_2' R_K^\transp \chi_{K_1'} \right) R_K(x_1).
	\end{aligned}
\end{align*}
Once we determined the candidates, it is straight-forward to show that indeed $\Phi_R$ is twice Frech{\'e}t differentiable and the bounds can easily be calculated. Finally, as the map $R_K \mapsto v_K$ is affine and bounded, also $\Phi_v$ is twice Frech{\'e}t differentiable with the given formula for the derivative and the bounds.
\end{proof}

\begin{remark} \label{Rem:Estimate_dR_by_K}
The bounds for $\dif \Phi_R$ determined in \cref{Lem:Change-of-variables:homeomorphism} imply that for all $K, \hat K \in \Leb^2(\omega, \R^{3\times 3}_{\skewMat})$,
\begin{equation}
	\norm*{R_K - R_{\hat K}}_{\Leb^\infty(\omega)} 
	\leq  c \norm*{K - \hat K}_{\Leb^2(\omega)}.
\end{equation}
This can also be seen more straight-forwardly by the calculation,
\begin{equation*}
	\norm*{R_K - R_{\hat K}}_{\Leb^\infty(\omega)} 
	= \norm*{\tilde R - I}_{\Leb^\infty(\omega)} 
	\leq c_1\norm*{\partial_1 \tilde R}_{\Leb^2(\omega)} 
	= c_1\norm*{K_{\tilde R}}_{\Leb^2(\omega)} 
	= c_1\norm*{K - \hat K}_{\Leb^2(\omega)},
\end{equation*}
where $\tilde R := R_{\hat K}R_K^\transp$ and an elementary calculation shows, $\hat K - K = R_K^\transp K_{\tilde R} R_K$.
\end{remark}

\subsection{Geometric rigidity estimate}
\label{Sec:GRE-Proofs}

In this section we prove \cref{Prop:RigidityRods}. To establish the statement we refine the construction in \cite[Prop.~3.6]{Neu12}. Therefore, we first prove a slight variation of \cite[Lem.~4.3]{Neu12}.

\begin{lemma}[cf.\ {\cite[Lem.~4.3]{Neu12}}] \label{Lem:RotationInterpolation}
Let $s_1 < s_2$ and $R_1, R_2 \in \SO(3)$. Then, there exists a map $R \in \Cont^\infty([s_1,s_2], \R^{3\times 3}) \cap \Cont([s_1,s_2], \SO(3))$ and a matrix $K \in \R^{3\times 3}_{\skewMat}$, such that
\begin{equation}
	R(s_1) = R_1, \qquad
	R(s_2) = R_2, \qquad
	R(s)^\transp \partial_1 R(s) = (s - s_1)(s_2-s)K \text{ for all } s \in [s_1, s_2],
\end{equation}
as well as
\begin{equation}
	\norm{\partial_1 R}_{\Leb^\infty([s_1,s_2])} \leq \frac{c\abs*{R_2 - R_1}}{s_2 - s_1}, \qquad
	\norm{\partial_1^2 R}_{\Leb^\infty([s_1,s_2])} \leq \frac{c\abs*{R_2 - R_1}}{(s_2 - s_1)^2},
\end{equation}
for some universal constant $c > 0$.
\end{lemma}

\begin{proof}
By a standard rescaling and translation argument, we may without loss of generality assume $s_1 = 0$ and $s_2 = 1$. By Rodrigues' rotation formula we find $N \in \R^{3\times 3}_{\skewMat}$ with $\abs*{N} = 1$ and $\alpha \in [0,\pi]$, such that
\begin{equation*}
	\operatorname{Rot}(N, \alpha) := I + \sin\alpha N + (1-\cos\alpha)N^2 = \exp(\alpha N) = R_1^\transp R_2.
\end{equation*}
Consider $f: \R \to \R$, $f(s) := 3s^2 - 2s^3$. We observe,
\begin{equation*}
	f(0) = 0, \qquad
	f(1) = 1, \qquad
	f'(s) = 6s(1-s), \qquad
	f''(s) = 6-12s.
\end{equation*}
We define
\begin{equation*}
	R(s) := R_1\operatorname{Rot}(N, \alpha f(s)), \qquad s \in [0,1].
\end{equation*}
Clearly, $R(0) = R_1$, $R(1) = R_2$ and we calculate 
\begin{equation*}
	R(s)^\transp\partial_1 R(s) = R(s)^\transp R_1\exp(\alpha f(s)N) \alpha f'(s) N = f'(s) \alpha N =: s(1-s)K.
\end{equation*}
Especially, this shows using $\abs*{N} = 1$,
\begin{equation*}
	\abs*{\partial_1 R(s)} = \abs*{R(s)^\transp\partial_1 R(s)} \leq 6\alpha
\end{equation*}
and \cite[Lem.~4.3]{Neu12} establishes $\alpha \leq c_1 \abs*{R_1^\transp R_2 - I} = c_1 \abs*{R_2 - R_1}$. Finally,
\begin{equation*}
	\partial_1 (R^\transp\partial_1 R)(s) = f''(s)\alpha N,
\end{equation*}
from which it is easy to establish $\abs*{\partial_1^2 R(s)} \leq c_2 \abs*{R_2 - R_1}$.
\end{proof}

The proof of \cref{Prop:RigidityRods} is now very close to \cite[Prop.~3.6]{Neu12}. For the readers convenience we sketch the arguments.

\begin{proof}[Proof of \cref{Prop:RigidityRods}]
Let  $y \in \SobW^{1,p}(\Omega, \R^3)$. All constants throughout this proof are independent of $h$.
\medskip

\stepemph{Step 1 -- Piece-wise constant maps:} \cite[Lem.~4.2]{Neu12} provides the construction of a piece-wise constant, right continuous map $R_{\rm pc}: \overline{\omega} \to \SO(3)$ with jump set contained in $h\Z \cap \omega$ which satisfy $R_{\rm pc}(0) = R_{\rm bc}$, as well as
\begin{equation*}
	\norm{\grad[h]y - R_{\rm pc}}_{\Leb^p(\Omega)}^p + h \sum_{\xi \in h\Z \cap \omega} \abs*{[[R_{\rm pc}(\xi)]]}^p \leq c_1 \norm{\dist(\grad[h]y, \SO(3))}_{\Leb^p(\Omega)}^p,
\end{equation*}
and
\begin{equation*}
	\sup_{\xi \in h\Z \cap \omega} \abs*{[[R_{\rm pc}(\xi)]]} \leq c_2 \norm{\dist(\grad[h]y, \SO(3))}_{\Leb^\infty(\Omega)}.
\end{equation*}
Here $[[R_{\rm pc}(\xi)]] := \lim_{s \to 0} R_{\rm pc}(\xi) - R_{\rm pc}(\xi - s)$ denotes the jump of $R_{\rm pc}$ at $\xi$. The construction uses the geometric rigidity estimate in \cref{Prop:RigidityEstimate} on the domain $U = [\xi,\xi+h) \times hS$. In fact, \cite[Lem.~4.2]{Neu12} provides the arguments for $p = 2$ but the extension to arbitrary $1 < p < \infty$ is straight-forward, see also \cite[Prop.~4.2]{NR26}. \cite[Prop.~4.2]{NR26} shows that we can use the explicit choice $R_{\rm pc}(\xi) = \operatorname{P}_{\SO(3)}(\fint_{[\xi,\xi+h \times hS} \grad[h]y)$. Moreover, we notice that it is possible to circumvent the extension argument (Step 1 in the proof of \cite[Lem.~4.2]{Neu12}) as one may choose $R_{\rm pc}$ (resp. $\hat R_{\rm pc}$) constant on $[\xi^*,l)$ where $\xi^* := \max\set*{\xi \in h\Z \given \xi \leq l - h}$, such that
\begin{equation*}
	\int_{[\xi^*,l) \times S} \abs*{\grad[h]y - R_{\rm pc}}^p \leq c_3 \int_{[\xi^*,l) \times S} \dist^p\of*{\grad[h]y, \SO(3)}.
\end{equation*}
Indeed, this is due to the fact that $c_3$ is independent of $h$, by $h \leq l - \xi^* \leq 2h$ and the uniformity of the constant of the geometric rigidity estimate for Bilipschitz equivalent domains.
\medskip

\stepemph{Step 2 -- Interpolation:} We now construct $R$ by interpolating the piece-wise constant map $R_{\rm pc}$. Let $\xi \in h\Z \cap [0,\xi^*)$. Then, by applying \cref{Lem:RotationInterpolation} to $s_1 = \xi$, $s_2 = \xi+h$, $R_1 = R_{\rm pc}(\xi)$ and $R_2 = R_{\rm pc}(\xi+h)$ we find a smooth map $R:[\xi,\xi+h] \to \SO(3)$ such that
\begin{equation*}
	R(\xi) = R_{\rm pc}(\xi), \qquad
	R(\xi+h) = R_{\rm pc}(\xi+h), \qquad
	\partial_1 R(\xi) = \partial_1 R(\xi+h) = 0,
\end{equation*}
and
\begin{equation*}
	\norm{\partial_1 R}_{\Leb^\infty([\xi,\xi+h])} \leq c_5 h^{-1} \abs*{[[R_{\rm pc}(\xi+h)]]}, \qquad
	\norm{\partial_1^2 R}_{\Leb^\infty([\xi,\xi+h])} \leq c_5 h^{-2} \abs*{[[R_{\rm pc}(\xi+h)]]}.
\end{equation*}
Similar, we can define $R:[\xi^*,l] \to \SO(3)$ and obtain
\begin{equation*}
	R(\xi^*) = R_{\rm pc}(\xi^*), \qquad
	R(l) = R_{\rm pc}(l), \qquad
	\partial_1 R(\xi^*) = 0,
\end{equation*}
and
\begin{equation*}
	\norm{\partial_1 R}_{\Leb^\infty([\xi^*,l])} \leq c_4 h^{-1} \abs*{[[R_{\rm pc}(\xi^*+h)]]}, \qquad
	\norm{\partial_1^2 R}_{\Leb^\infty([\xi^*,l])} \leq c_4 h^{-2} \abs*{[[R_{\rm pc}(\xi^*+h)]]}.
\end{equation*}
This defines a map $R \in \SobW^{2,\infty}_{\SO(3)}(\omega)$ and the estimates above combined with Step 1 yield \eqref{Eq:RigidityRods}. 
\end{proof}

\subsection{An Arzel{\`a}--Ascoli type theorem}

In this section we prove a generalized version of Arzel{\`a}--Ascoli, which we used to prove the general $\Gamma$-convergence theorem for curves of maximal slope in \cref{Thm:General-convergence-result}. The statement and proof is very similar to \cite[Prop.~3.3.1]{AGS05} and \cite[Thm.~3.5]{FK18} and we provide the proof for the readers convenience.

\begin{lemma}[Arzel{\`a}--Ascoli, cf.\ {\cite[Prop.~3.3.1]{AGS05}}] \label{Lem:Arzela-Ascoli}
Let $(Q^k, D^k)$ and $(Q, D)$ be metric spaces for $k \in \N$. Let $(\tilde Q, \sigma)$ be a topological Hausdorff space, and $\pi^k: Q^k \to \tilde Q$ and $\pi: Q \to \tilde Q$ maps. Let $q^k: [0,T] \to \mathcal Q^k$ be curves. Assume that the following statements hold:
\begin{enumerate}[(a)]
	\item\label{Item:Arzela-Ascoli:a} For all $t \in [0,T]$ and for any subsequence of $(q^k(t))$, there exists a further subsequence (not relabeled) and some $q \in Q$, such that $\pi^k(q^k(t)) \to \pi(q)$.
	
	\item\label{Item:Arzela-Ascoli:b} Let $s,t \in [0,T]$, $q,\hat q \in Q$ and consider a subsequence (not relabeled) such that $\pi^k(q^k(s)) \to \pi(q)$ and $\pi^k(q^k(t)) \to \pi(\hat q)$. Then,
	\begin{equation*}
		\liminf_{k \to \infty} D^k\of*{q^k(s), q^k(t)} \geq D(q,\hat q).
	\end{equation*}
	
	\item\label{Item:Arzela-Ascoli:c} There exists an at most countable set $\mathscr C \subset [0,T]$ and a symmetric map $\omega: [0,T]^2 \to [0,\infty)$, such that
		\begin{align*}
			\lim_{(s,t) \to (r,r)} \omega(s,t) &= 0, && \text{for all } r \in [0,T] \setminus \mathscr C, \text{ and} \\
			\limsup_{k \to \infty} D^k\of*{q^k(s),q^k(t)} &\leq \omega(s,t), && \text{for all } s,t \in [0,T].
 		\end{align*}
\end{enumerate}
Then, there exists a subsequence (not relabeled, independent of $t$) and a limiting curve $q: [0,T] \to Q$, such that $q$ is continuous w.r.t.\ $D$ in $[0,T] \setminus \mathscr C$ and
\begin{equation} \label{Eq:Arzela-Ascoli:convergence}
	\pi^k(q^k(t)) \to \pi(q(t)), \qquad\text{for all } t \in [0,T].
\end{equation}
\end{lemma}

\begin{proof}
By a standard diagonal sequence argument and due to (\ref{Item:Arzela-Ascoli:a}), we find a subsequence (not relabeled) and some map $q: ([0,T]\cap \Q) \cup \mathscr C \to Q$, such that
\begin{equation*}
	\pi^k(q^k(t)) \to \pi(q(t)) \qquad\text{for all } t \in ([0,T]\cap \Q) \cup \mathscr C.
\end{equation*}
Let $t \in [0,T] \setminus (\Q \cup \mathscr{C})$. Due to (\ref{Item:Arzela-Ascoli:a}) any subsequence admits a further subsequence (not relabeled) and a limit $q(t) \in Q$, such that $\pi^k(q^k(t)) \to \pi(q(t))$. Let $(t_n) \subset ([0,T]\cap \Q) \cup \mathscr C$ be such that $t_n \to t$. Since, $t \nin \mathscr{C}$, we obtain from (\ref{Item:Arzela-Ascoli:b}) and (\ref{Item:Arzela-Ascoli:c})
\begin{equation*}
	\limsup_{n \to \infty} D\of[\big]{q(t), q(t_n)}
	\leq \limsup_{n \to \infty}\liminf_{k \to \infty} D^k\of[\big]{q^k(t), q^k(t_n)}
	\leq \limsup_{n \to \infty} \omega(t,t_n) = 0.
\end{equation*}
As $D$ is a metric, this defines $q(t)$ uniquely, so that $q(t)$ is independent of the chosen subsequence. Hence, the subsequence principle implies $\pi^k(q^k(t)) \to \pi(q(t))$ for the whole sequence. Finally, we obtain that the so defined map $q: [0,T] \to Q$ satisfies for all $s,t \in [0,T]$,
\begin{equation*}
	D\of[\big]{q(s), q(t)}
	\leq \liminf_{k \to \infty} D^k\of[\big]{q^k(s), q^k(t)}
	\leq \omega(s,t),
\end{equation*}
which shows that $q$ is continuous w.r.t.\ $D$ in $[0,T] \setminus \mathscr{C}$.
\end{proof}

\clearpage
\printbibliography[heading=bibintoc]
\end{document}